\documentclass[10pt,a4paper,twoside]{amsart}
\usepackage{etex}
\usepackage{amsmath,amssymb,amsfonts,enumerate,amsthm,graphicx,bm,fancyhdr,pictexwd,dcpic,mathrsfs,verbatim,aligned-overset}
\usepackage[hidelinks]{hyperref}
\usepackage[dvipsnames]{xcolor}
\usepackage{tikz-cd,hyperref}
\usetikzlibrary{arrows}
\usepackage[all]{xy}
\usepackage{amscd}

\xyoption{all}
\headsep=1cm \numberwithin{equation}{section}
\newtheorem{thm}{Theorem}[section]
    \theoremstyle{Definition}
    
    \theoremstyle{Definition and Remark}
    \newtheorem{defi-rem}[thm]{Definition and Remark}
    \newtheorem{defi-rems}[thm]{Definition and Remarks}
    \newtheorem{defis-rems}[thm]{Definitions and Remarks}
    \newtheorem{defi-Nots}[thm]{Definition and Notations}
    \newtheorem{defi-Not}[thm]{Definition and Notation}
    
    \theoremstyle{Lemma}

    \theoremstyle{Corollary}

    \newtheorem{ques}[thm]{Question}

\newtheorem{subthm}{Theorem}[subsection]
\newtheorem{subdefi}[subthm]{Definition}

\newtheorem{subdefi-rem}[subthm]{Definition and Remark}
\newtheorem{subdefi-rems}[subthm]{Definition and Remarks}
\newtheorem{subdefis-rems}[subthm]{Definitions and Remarks}
\newtheorem{subdefi-Nots}[subthm]{Definition and Notations}
\newtheorem{subdefi-Not}[subthm]{Definition and Notation}

\newtheorem{subnota}[subthm]{Notation}
\newtheorem{sublem}[subthm]{Lemma}
\newtheorem{subrem}[subthm]{Remark}
\newtheorem{subcor}[subthm]{Corollary}

\newtheorem{subprop}[subthm]{Proposition}

\newtheorem{subfact}[subthm]{Fact}
\newtheorem{subconv}[subthm]{Convention}
\DeclareMathOperator{\Hom}{Hom} 

\DeclareMathOperator{\Ker}{Ker}   
 \DeclareMathOperator{\Ass}{Ass}

 \DeclareMathOperator{\Ext}{Ext}
\DeclareMathOperator{\depth}{depth}
 \DeclareMathOperator{\Tot}{Tot}

\newcommand\numberthis{\addtocounter{equation}{1}\tag{\theequation}}

\newcommand{\fm}{\mathfrak{m}}

\newcommand{\fp}{\frak{p}}
\newcommand{\fq}{\frak{q}}
\newcommand{\fa}{\frak{a}}
\newcommand{\fb}{\frak{b}}

\newcommand{\fn}{\frak{n}}

 \newcommand{\mam}{\mathfrak{m}}  \newcommand{\maa}{\mathfrak{a}}
\newcommand{\mab}{\mathfrak{b}}   
  
 \newcommand{\mn}{\mathbb{N}}

\newcommand{\Max}{\text{Max}}
\newcommand{\pd}{\text{pd}} \newcommand{\id}{\text{Id}}

 \newcommand{\im}{\text{im}}

\newcommand{\dsum}{\bigoplus}

 \newcommand{\tor}{\text{Tor}}

 \newcommand{\Frac}{\text{Frac}} 

\DeclareMathAlphabet{\mathcalligra}{T1}{calligra}{g}{f}

\newcommand{\llar}{-\kern-5pt-\kern-5pt\longrightarrow}
\def\restr{{\kern-1pt\restriction\kern-1pt}}

\usepackage[paper=a4paper,left=30mm,right=20mm,top=25mm,bottom=30mm]{geometry}

\begin{document}

	\title{Test modules, weakly regular homomorphisms and complete intersection dimension}
	\author[Ehsan Tavanfar]{Ehsan Tavanfar}
	\address{School of Mathematics, Institute for 	Research in Fundamental Sciences (IPM), P. O. Box: 	19395-5746, Tehran, Iran. }
	\email{tavanfar@ipm.ir and tavanfar@gmail.com}	
	\maketitle
	\markboth{\small \textsc{E. Tavanfar}}{\small  \textsc{Test modules, weakly-regular homomorphisms and complete intersection dimension}}

\date{\today}
\makeatletter{\renewcommand*{\@makefnmark}{}
	\footnotetext{(MSC2010)  13H10, 13D05.}\makeatother}
\makeatletter{\renewcommand*{\@makefnmark}{}
	\makeatother}
\makeatletter{\renewcommand*{\@makefnmark}{}
	\footnotetext{This research was supported by a grant from IPM.}\makeatother}

  \begin{abstract}
 We prove that  if a local ring admits a (pd-)test module of finite complete intersection dimension, then it is  a complete intersection ring. This answers, positively,   a  question proposed by  Celikbas, Dao and Takahashi.  
 To this aim,  we first investigate   another   question raised by  Celikbas and Sather-Wagstaff concerning   ascent properties of (pd-)test complexes   under weakly regular homomorphisms.
\end{abstract}


\tableofcontents

\section{Introduction}

Complete Intersection dimension introduced by Avramov, Gasharov and Peeva in \cite{Avramov}, is a homological invariant lying
 between the classical projective dimension and the	Auslander-Bridger Gorenstein dimension.  Over a complete intersection local ring every finitely generated module has finite complete intersection dimension.  The converse also holds. Namely,  a local ring $(R,\mam)$ is complete intersection if and only if its residue field  $R/\mam$   has finite complete intersection dimension (see \cite{Avramov}). Thus, in view of the fact that $R/\fm$ tests  finiteness of projective dimension of finitely generated modules via   vanishing of positive Tor modules $\tor^R_i(R/\mam,-)$,  the following   question is  proposed by Celikbas, Dao and Takahashi in 2014.

\begin{ques}\label{CDTQuestion} \emph{(\cite[Question 3.5]{Celikbas})
	Let $R$ be a local ring. Let $M$ be a (pd)-test module with $\text{CI-dim}_R(M)<\infty$. Then must $R$ be a complete intersection?}
\end{ques}

The aim of this paper is to answer this question in the affirmative, in general (see Corollary \ref{MainResult}).
In order to accomplish this, we study the  following question proposed by Celikbas and Sather-Wagstaff  whose module version is strongly connected to Question \ref{CDTQuestion}. 

\begin{ques}\label{CWQuestion} 
	\emph{(\cite[Question 3.7]{CelikbasWagstaff}) Let $\varphi:(R,\fm)\rightarrow (S,\fn)$ be  a flat local ring homomorphism  and let $M\in \mathcal{D}^f_b(R)$. Assume that $S/\fm S$ is regular. If  $M$ is a (pd-)test complex over $R$, then must $S\otimes_R^{\mathbf{L}}M$  be  a (pd-)test complex for $S$?}
\end{ques}

 Celikbas and Sather-Wagstaff  answered   Question \ref{CWQuestion} affirmatively in the case where $\varphi:R\rightarrow S$ is a  weakly unramified
  flat local homomorphism inducing finite field extension on residue fields (\cite[Theorem 3.5]{CelikbasWagstaff}).
 Afterwards, a further  step was taken    by Sather-Wagstaff  (\cite[Theorem 4.8]{WagstaffAscent}), where he settled Question \ref{CWQuestion} affirmatively under the assumption that the induced residue field extension by $\varphi$ is algebraic.  We   answer   Question \ref{CWQuestion} positively, in the case where $R$ has uncountable residue field. 
  Moreover, towards answering Question \ref{CDTQuestion},  we were able to relax this ``uncountable condition'' on the residue field of $R$  
  under another assumption that the test complex $M$ (in the statement of Question \ref{CWQuestion}) has finite complete intersection dimension ‌(see   Corollary \ref{TestComplexAndFlatHomomorphism}).

The last section of the paper is devoted to give  alternative proofs for Corollary \ref{TestComplexAndFlatHomomorphism} (which settles affirmatively  Question \ref{CWQuestion} provided $R/\fm$ is uncountable or $\text{CI-dim}_R(M)<\infty$).  The first proof of Corollary \ref{TestComplexAndFlatHomomorphism}, applies the Sather-Wagstaff's recent result \cite[Theorem 4.8]{WagstaffAscent}.
The proof of \cite[Theorem 4.8]{WagstaffAscent} is based on    translating the situation of Question \ref{CWQuestion} to a similar question on finite dimensional  DG-algebras over a field (by passing to the Koszul complex of a minimal generating set of the maximal ideal). In contrast, our alternative proof  uses a different method which is based on reducing the question to the case of coefficient ring base change extensions $R\rightarrow R\widehat{\otimes}_{C_R}C'$ where $C_R\rightarrow C'$ induces  a finitely generated field extension on residue fields (here $C_R$ is a coefficient ring for $R$, and  $C_R\rightarrow C'$ is either a field extension or an extension of complete unramified discrete valuation rings of mixed characteristic). But the desired  reduction to these  coefficient ring base change extensions was not possible for us in general and our alternative proofs are given in two  special cases.
Namely, lets keep in mind the notation of Question \ref{CWQuestion}.  In the case where either $\varphi$ induces (possibly non-algebraic) separable  field extension on residue fields (Corollary \ref{CorollaryAlternativeProofCWQuestionSeparableCase}), or   $R$ is a Cohen-Macaulay ring  containing a field (Corollary \ref{CorollaryAlternativeProofCohenMacaulayCase}), we succeeded to reduce  Question \ref{CWQuestion} to the aforementioned case of coefficient ring base change extensions, by which our  alternative proofs for Corollary \ref{TestComplexAndFlatHomomorphism} are given. 
  Since  Question \ref{CDTQuestion} is basically  a  question on Cohen-Macaulay rings, so our alternative proofs for Corollary \ref{TestComplexAndFlatHomomorphism} provide  also  an alternative proof, avoiding the use of Sather-Wagstaff's Theorem, for the equicharacteristic case of Corollary \ref{MainResult} which answers affirmatively Question \ref{CDTQuestion} (see Corollary \ref{CorollaryAlternativeProofCDTQuestion}).

Many of the results of this paper use  complete tensor product over a coefficient ring of a complete local ring. As far we know, the only reference to  complete tensor products and their properties   is Grothendieck's \cite{GrothendieckIV} (for complete tensor products over a field, or more generally over an absolutely flat ring,  see  \cite{TabaaCompleted}).  In Section 2, we  investigate and prove  some properties of  complete tensor products over arbitrary local rings that are required in our paper. Some of our results concerning complete tensor products are new (See Remark \ref{RemarkSomeResutlsAreNew}).

\section{Complete tensor product} 

In this section, we discuss  complete tensor products over an arbitrary local ring $\Lambda$. But all complete tensor products  in the subsequent sections of the paper are  over a field or over a complete discrete valuation ring.  We begin with  a preparatory  subsection containing  some general notation and conventions that will be used throughout this paper, as well as some elementary facts that will be  used in this section (or also later on).

\subsection{Preliminaries}

All  rings are commutative  with identity. The set of all maximal ideals of a ring $R$ is denoted by $\text{Max}(R)$. The notation $\text{Char}(R)$ stands for the characteristic of the ring $R$. 
The fraction field of an integral domain $R$ is denoted by $\Frac(R)$. If $\mathbf{X}$ is a set of variables over a domain $R$ then $R(\mathbf{X}):=\text{Frac}(R[\mathbf{X}])$.
In this paper by a local ring we mean a Noetherian ring with a unique maximal ideal, while a quasi-local ring is a (not necessarily Noetherian) ring  containing only one maximal ideal.  Given a ring $R$, an $R$-module $M$ and  submodules $N'\subseteq N$ of $M$, by the natural epimorphism   $M/N'\rightarrow M/N$  we   mean the map $m+N'\mapsto m+N$.    

 An $R$-algebra $A$ is said to be of finite type  provided $A$ is a finitely generated $R$-algebra, and it is called essentially of finite type provided $A$ is a localization of a finitely generated  $R$-algebra.   A  field extension $\mathcal{K}/K$ is said to be a finitely generated field extension  provided $\mathcal{K}=\Frac(K[\kappa_1,\ldots,\kappa_n])$ for a finite subset $\{\kappa_1,\ldots,\kappa_n\}$ of  $\mathcal{K}$. A finitely generated field extension   of $K$ can be represented as  a finite algebraic extension of a purely transcendental extension of $K$ with finite transcendence basis. 

\begin{subconv}\label{FirstConvention}
	\emph{For a quasi-local ring, say e.g.   $R$, its maximal ideal   (resp. its residue field) will  be  denoted by $\mathfrak{m}_R$ (resp. by $K_R$), unless otherwise is stated explicitly. So by a quasi-local ring $R$ a triple $(R,\mathfrak{m}_R,K_R)$ is  meant while we often  write only $R$. If $R$ is a (quasi-)local ring then a local $R$-algebra $A$  is an $R$-algebra  $\varphi:R\rightarrow A$  such that $A$ is a local ring and $\varphi$ is a (quasi-)local homomorphism, i.e. $\varphi(\mam_R)\subseteq \mam_A$. If $R$ and $S$ are local rings, then an $(R,S)$-local algebra is a ring $A$ that is   an $R$-local algebra and an $S$-local algebra and it is an $(R,S)$-bimodule with the induced module structures.}
\end{subconv}

\begin{subfact}\label{ElementaryTensorProductFact}
	Suppose that $\Lambda$ is a ring and  that $U$ and $V$ are two $\Lambda$-algebras. Let $\mathfrak{u}$ be an ideal of $U$, and $\mathfrak{v}$ be an ideal of $V$. Finally, let $L$ be an $\Lambda$-module. 
	\begin{enumerate}
		\item[(i)]
		We have the natural isomorphism of $U$-modules (and $U/\mathfrak{u}$-modules), $$(U\otimes_{\Lambda}L)/\big(\mathfrak{u}(U\otimes_{\Lambda}L)\big)\rightarrow (U/\mathfrak{u})\otimes_{\Lambda}L,$$ given by $u\otimes x+\mathfrak{u}(U\otimes_{\Lambda}L)\mapsto (u+\mathfrak{u})\otimes x,\ \ \ (u\in U\text{\ \ and\ \ }x\in L).$
		\item[(ii)]
		The natural ring epimorphism $$\pi:U\otimes_{\Lambda}V\twoheadrightarrow (U/\mathfrak{u})\otimes_{\Lambda}(V/\mathfrak{v}),\ \ \ u\otimes v\mapsto (u+\mathfrak{u})\otimes (v+\mathfrak{v})$$ induces the isomorphism   $$U\otimes_{\Lambda}V/\big(\mathfrak{u}(U\otimes_{\Lambda}V)+\mathfrak{v}(U\otimes_{\Lambda}V)\big)\cong (U/\mathfrak{u})\otimes_{\Lambda}(V/\mathfrak{v}).$$
	\end{enumerate}
	\begin{proof}
		(i) Applying  $-\otimes_\Lambda L$     to the exact sequence $0\rightarrow \mathfrak{u}\overset{g}{\rightarrow} U\rightarrow U/\mathfrak{u}\rightarrow 0$ yields the exact sequence $\mathfrak{u}\otimes_\Lambda L\overset{g\otimes \id_L}{\longrightarrow} U\otimes_\Lambda L\rightarrow (U/\mathfrak{u})\otimes_\Lambda L\rightarrow 0$. Then the statement follows from $$\mathfrak{u}(U\otimes_{\Lambda}L)=(g\otimes \text{Id}_L)\big(\mathfrak{u}\otimes_{\Lambda}L\big).$$ 
		
		(ii) When $\mathfrak{v}=0$, the statement follows from the previous part by setting $L:=V$.  Then the the general case  follows from the Third Isomorphism Theorem in conjunction with the above special case.
	\end{proof}
\end{subfact}

\begin{subfact}\label{ApplyingSharpResult}
	Suppose that $K\subseteq \mathcal{K}$ is an extension of fields, and  that $R$ is a local $K$-algebra such  that  $K_R\otimes_{K}\mathcal{K}$ is a field. Let $\mathcal{L}$ be a subfield of $\mathcal{K}$  containing $K$. Then $\mathfrak{m}_R(R\otimes_{K}\mathcal{L})\in \Max(R\otimes_{K}\mathcal{L})$.
	\begin{proof}
		As $K_R\otimes_K\mathcal{K}$ is a field, at least one of  $K_R,\mathcal{K}$ is algebraic over $K$ by virtue of \cite[Theorem 3.1]{SharpTheDimension}. So $\dim(K_R\otimes_{K}\mathcal{L})=0$ again by \cite[Theorem 3.1]{SharpTheDimension},         while $K_R\otimes_{K}\mathcal{L}$ is a domain (it is a subring  of $K_R\otimes_{K}\mathcal{K}$, as $K_R$ is   flat over $K$). So $K_R\otimes_{K}\mathcal{L}$ is a field, and $\mathfrak{m}_R(R\otimes_{K}\mathcal{L})$ is  maximal  (Fact \ref{ElementaryTensorProductFact}(ii)). 
	\end{proof}
\end{subfact}

\begin{subfact}\label{DirectLimit}
	Suppose that $R$ is a ring and let $\mathfrak{p}\in \text{Spec}(R)$. Assume that   $R=\lim\limits_{\underset{i\in I}{\longrightarrow}}R_i$ for some direct system  $\{R_i,\ \varphi_{i,j}:R_i\rightarrow R_j\}_{{i,j\in I,i\le j}}$ of rings and ring homomorphisms. For each $i\in I$, let $\varphi_i:R_i\rightarrow R$ be the map to the direct limit and set $\fp_i:=\varphi_i^{-1}(\mathfrak{p})$. Then,  $R_{\mathfrak{p}}\cong \lim\limits_{\underset{i\in I}{\longrightarrow}}{R_{i}}_{\mathfrak{p}_i}$.
	\begin{proof}
		For each $i,j\in I$ with $i\le j$ we have $\mathfrak{p}_i=\varphi_{i,j}^{-1}(\mathfrak{p}_j)$. Hence we get the  maps ${R_i}_{\mathfrak{p}_i}\rightarrow {R_j}_{\mathfrak{p}_j}$ induced by $\varphi_{i,j}$ (also ${\varphi_i}_{\mathfrak{p}}:{R_i}_{\mathfrak{p}_i}\rightarrow R_{\mathfrak{p}}$)  in a compatible way such that they induce the  ring homomorphism $\lambda:\lim\limits_{\underset{i\in I}{\longrightarrow}}\ {R_i}_{\mathfrak{p}_i}\rightarrow R_\mathfrak{p}$ (by the universal property of the direct limit). For each fraction $r/s\in R_\mathfrak{p}$ we can consider some $i\in I$, $r_i\in R_i$ and $s_i\in R_i\backslash \mathfrak{p}_i$ such that  $r=\varphi_i(r_i)$ and $s=\varphi_i(s_i)$. Thus $r/s={\varphi_i}_{\mathfrak{p}}(r_i/s_i)$. This shows that $\lambda$ maps the image of, $r_i/s_i\in {R_i}_{\mathfrak{p}_i}$, in $\lim\limits_{\underset{i\in I}{\longrightarrow}}\ {R_i}_{\mathfrak{p}_i}$, to $r/s$. Thus $\lambda$  is surjective. The injectivity also follows similarly, because any relation $rs=0$ ($r\in R$ and $s\in R\backslash \mathfrak{p}$) comes from a similar relation in some $R_i$ $(i\in I)$.
	\end{proof}
\end{subfact}

\subsection{Complete tensor product}

We begin by giving the definition of complete tensor products. We   remind the reader of Convention \ref{FirstConvention}.   

\begin{subdefi-Not}\label{DefinitionCompeleteTensorProduct} \emph{Let $\Lambda$ be a local ring.  The \textit{complete tensor product} of    $\Lambda$-local algebras $A$ and $B$ is   the inverse limit  $$A\widehat{\otimes_{\Lambda}}B:=\lim\limits_{\overset{\longleftarrow}{n\in \mn}} \big((A/\mam_A^n)\otimes_{\Lambda}(B/\mam_B^n)\big),$$ of the inverse system $\Big((A/\mam_A^n)\otimes_{\Lambda}(B/\mam_B^n),\ \xi_n\Big)_{n\in \mathbb{N}}$ with the transition homomorphisms  $$\xi_n:(A/\mam_A^{n+1})\otimes_{\Lambda}(B/\mam_B^{n+1})\rightarrow (A/\mam_A^n)\otimes_{\Lambda}(B/\mam_B^n), \ \ \big((a+\mathfrak{m}_A^{n+1})\otimes (b+\mathfrak{m}_B^{n+1})\big)\mapsto \big((a+\mathfrak{m}_A^{n})\otimes (b+\mathfrak{m}_B^{n})\big).$$  We set $\mathfrak{M}_{A,\Lambda,B}:=\mam_A(A\otimes_{\Lambda}B)+\mam_B(A\otimes_{\Lambda}B)$ and set  $\mathfrak{M}_{A,\Lambda,B}^{e}$ to be the extension of $\mathfrak{M}_{A,\Lambda,B}$ under the natural ring homomorphism $$\eta_{A\otimes_{\Lambda}B}:A\otimes_\Lambda B\rightarrow A\widehat{\otimes}_{\Lambda}B,\ \ \ a\otimes b\mapsto \big((a+\mam_A^n)\otimes(b+\mam_B^n)\big)_{n\in\mn}.$$ When no confusion is likely, we   (often but not always) simplify the notation by writing $\mathfrak{M}$ (resp. $\mathfrak{M}^e$) in place of $\mathfrak{M}_{A,\Lambda,B}$ (resp. $\mathfrak{M}_{A,\Lambda,B}^e$). We consider $A\widehat{\otimes}_\Lambda B$ as an $A$-algebra (resp. $B$-algebra) in the natural way, i.e. via  $$\eta_{A,A\widehat{\otimes}_\Lambda B}:=A\overset{a\mapsto a\otimes 1}{\longrightarrow} A\otimes_{\Lambda}B\overset{\eta_{A\otimes_{\Lambda}B}}{\longrightarrow} A\widehat{\otimes}_{\Lambda}B \text{\ \ \ (resp.\ \ \ }  \eta_{B,A\widehat{\otimes}_{\Lambda}B}:=B\overset{b\mapsto 1\otimes b}{\rightarrow} A\otimes_{\Lambda}B\overset{\eta_{A\otimes_{\Lambda}B}}{\longrightarrow} A\widehat{\otimes}_{\Lambda}B)$$}
\end{subdefi-Not}

\begin{subrem}\label{RemarkDefinitionCompleteTensorProductDoesntAffectCompleteness}
	  $A\widehat{\otimes}_{\Lambda}B\cong \widehat{A}\widehat{\otimes}_{\Lambda}\widehat{B}$ and that $A\widehat{\otimes}_{\Lambda}B\cong B\widehat{\otimes}_{\Lambda}A$, by definition. 
\end{subrem}	

\begin{subrem}\label{RemarkNaturalIsomorphism}
	\emph{Let $\Lambda$ be a local ring with an ideal $\mathfrak{d}$, and $A$ be  a local $\Lambda$-algebra. Then we have the $\Lambda$-algebra isomorphism $$\widehat{A/\mathfrak{d}A}\rightarrow (\Lambda/\mathfrak{d})\widehat{\otimes}_{\Lambda}A, \ \ 
		\big((a_n+\mathfrak{d}A)+(\mathfrak{m}_A^n)\big)_{n\in\mathbb{N}}\mapsto \Big(\big((1+\mathfrak{d})+(\mathfrak{m}_{\Lambda}^n)\big)\otimes(a_n+\mathfrak{m}_A^n)\Big)_{n\in \mathbb{N}}.$$
		Namely, $\mathfrak{m}_\Lambda^nA\subseteq \mathfrak{m}_A^n$ by the local algebra property. Thus   we have the isomorphisms of inverse systems
		\begin{align*}
			\Big((A/\mathfrak{d}A)/(\mathfrak{m}_A^n)
			\Big)_{n\in\mathbb{N}}
			\cong	    
			\Big(A/(\mathfrak{d}A+\mathfrak{m}_A^n+\mathfrak{m}_\Lambda^nA)
			\Big)_{n\in\mathbb{N}}
			\cong
			\Big(\big((\Lambda/\mathfrak{d})/(\mathfrak{m}_{\Lambda}^n)\big)\otimes_{\Lambda}(A/\mathfrak{m}_A^n),\ \xi_n\Big)_{n\in\mathbb{N}},
		\end{align*}
		where the transition homomorphisms of the two left-most ones are the natural epimorphisms.}
\end{subrem}

In  the statement of Lemma \ref{cofinal}(i), the transition maps of the  inverse system corresponding to each inverse limit     are given by the tensor product of the natural epimorphisms.

\begin{sublem}\label{cofinal}
	Let $A$ be a ring, and $M$ and  $N$ be two $A$-modules. Suppose that $\{M_n\}_{n\in \mathbb{N}}$ and $\{M'_n\}_{n\in\mathbb{N}}$ are two decreasing sequences of submodules of $M$ such that for each $n\in \mathbb{N}$ there is some $h_n\ge n$ with $M_{h_n}\subseteq M'_{n}\subseteq M_{n}$ (thus  $\{M_n\}_{n\in \mathbb{N}}$  is cofinal with  $\{M'_n\}_{n\in\mathbb{N}}$). Assume that  $\{N_n\}_{n\in \mathbb{N}}$ is a decreasing sequence of submodules of $N$. Let $\pi_n:M/M'_n\rightarrow M/M_n$ be the natural epimorphism for each $n\in \mathbb{N}$. 
	\begin{enumerate}
		\item[(i)] The natural map $$(\id_{N/N_n}\otimes\pi_n)_{n\in\mathbb{N}}:\lim\limits_{\underset{n\in \mathbb{N}}{\longleftarrow}}\big((N/N_n)\otimes_A(M/M'_n)\big)\rightarrow \lim\limits_{\underset{n\in \mathbb{N}}{\longleftarrow}}\big((N/N_n)\otimes_A(M/M_n)\big)$$ is an isomorphism 
		\item[(ii)] The natural map $(\pi_n)_{n\in\mathbb{N}}:\lim\limits_{\underset{n\in \mathbb{N}}{\longleftarrow}}(M/M'_n)\rightarrow \lim\limits_{\underset{n\in \mathbb{N}}{\longleftarrow}}(M/M_n)$ is an isomorphism.
	\end{enumerate}
	
	\begin{proof} (i) Let $\delta_n:M_n/M'_n\rightarrow M/M'_n$ be the inclusion map and  set $$K_n:=\ker\Big(\big((N/N_{n})\otimes_{A}(M_{n}/M'_{n})\big)\overset{\id_{N/N_n}\otimes \delta_n}{\longrightarrow} \big((N/N_{n})\otimes_{A}(M/M'_{n})\big)\Big)$$ for each $n\in \mathbb{N}$. Then   			
		$$
		0\rightarrow\{\big((N/N_{n})\otimes_{A}(M_{n}/M'_{n})\big)/K_{n}\}\overset{\{\overline{\id_{N/N_n}\otimes \delta_n}\}}{\longrightarrow}\{(N/N_{n})\otimes_{A}(M/M'_{n})\}\overset{\{\id_{N/N_n}\otimes \pi_n\}}{\longrightarrow}\{(N/N_{n})\otimes_{A}(M/M_{n})\}\rightarrow0
		$$
		is an exact sequence	of inverse systems.			 Moreover,  the transition map  
		\begin{equation}
		\label{ML1}
		\big((N/N_{m})\otimes_A (M_{m}/M'_{m})\big)/{K_m}\rightarrow \big((N/N_n)\otimes_{A}(M_n/M'_n)\big)/{K_n}
		\end{equation}
		given by $\big((y+N_{m})\otimes (x+M'_{m})\big)+K_m\mapsto \big((y+N_n)\otimes (x+M'_n)\big)+K_n$ is zero for each $m\ge h_n$ (because  $M_{h_n}\subseteq M'_n$ by our hypothesis).  This, on the one hand shows that the inverse system $$\{\big((N/N_{n})\otimes_{A}(M_{n}/M'_{n})\big)/K_{n}\}$$ satisfies the Mittag-Leffler condition whence by \cite[II, Proposition 9.1(b)]{HartshorneAlgebra} we get the exact sequence,
		\begin{align*}
		\label{MittagLeffler2}
		0
		\numberthis\rightarrow 
		\lim\limits_{\underset{n\in \mathbb{N}}{\longleftarrow}}\Big(\big((N/N_{n})\otimes_{A}(M_{n}/M'_{n})\big)/K_{n}\Big)
		&\rightarrow  
		\lim\limits_{\underset{n\in \mathbb{N}}{\longleftarrow}}\big((N/N_n)\otimes_A(M/M'_n)\big)
		&\\&\hspace{-9mm}\overset{(\id_{N/N_n}\otimes\pi_n)_{n\in\mathbb{N}}}{\rightarrow}  
		\lim\limits_{\underset{n\in \mathbb{N}}{\longleftarrow}}\big((N/N_n)\otimes_A(M/M_n)\big)
		\rightarrow
		0.
		\end{align*}
		On the other hand the vanishing of the maps in (\ref{ML1}) for each $n$ and $m\ge h_n$ shows that $$\lim\limits_{\underset{n\in \mathbb{N}}{\longleftarrow}}\Big(\big((N/N_{n})\otimes_{A}(M_{n}/M'_{n})\big)/K_{n}\Big)=0.$$ Hence the desired conclusion follows from the exact sequence (\ref{MittagLeffler2}).
		
		(ii) The desired conclusion follows from part (i) by setting $N=A$ and $N_n=0$ for each $n\in \mn$.
	\end{proof}
\end{sublem}

\begin{subprop}\label{CompleteTensorProductNoetherianness}Let   $\Lambda$ be a local ring and let  $A,B$ be two  $\Lambda$-local algebras.    Then,
	\begin{enumerate}
		\item[(i)]   $A\widehat{\otimes}_\Lambda B$ (via $\eta_{A\otimes_{\Lambda}B}$) is $(A\otimes_{\Lambda}B)$-algebra isomorphic to the $\mathfrak{M}$-adic completion of $A\otimes_{\Lambda}B$. 
		\item[(ii)] $A\widehat{\otimes}_{\Lambda}B$  is Noetherian and $\mathfrak{M}^e$-adically complete provided  $(A\otimes_{\Lambda}B)/\mathfrak{M}$ is Noetherian.
		
		\item[(iii)] If $\mathfrak{M}\in \text{Max}(A\otimes_{\Lambda}B)$, then $A\widehat{\otimes}_{\Lambda}B$ is $(A\otimes_{\Lambda}B)$-algebra isomorphic to the $\mathfrak{M}(A\otimes_{\Lambda}B)_{\mathfrak{M}}$-adic completion of $(A\otimes_{\Lambda}B)_{\mathfrak{M}}$.
	\end{enumerate} 
	\begin{proof}
		(i) The  decreasing sequence  $\{\mam_A^n(A\otimes_{\Lambda}B)+\mam_B^n(A\otimes_{\Lambda}B)\}_{n\in \mn}$ of ideals of $A\otimes_{\Lambda}B$ is cofinal with the adic filtration $\{\mathfrak{M}^n\}_{n\in\mn}$ (because $\mathfrak{M}$ is finitely generated). Thus the $\mathfrak{M}$-adic completion of $A\otimes_\Lambda B$ is $(A\otimes_{\Lambda}B)$-algebra isomorphic to   $\lim\limits_{\underset{n\in \mn}{\longleftarrow}}\Big((A\otimes_{\Lambda}B)/\big(\mam_A^n(A\otimes_{\Lambda}B)+\mam_B^n(A\otimes_{\Lambda}B)\big)\Big)$  by Lemma \ref{cofinal}(ii). Also, the isomorphism of Fact \ref{ElementaryTensorProductFact}(ii) can be used to construct an isomorphism of inverse systems yielding  the $(A\otimes_{\Lambda}B)$-algebra isomorphism  $\lim\limits_{\underset{n\in \mn}{\longleftarrow}}\Big((A\otimes_{\Lambda}B)/\big(\mam_A^n(A\otimes_{\Lambda}B)+\mam_B^n(A\otimes_{\Lambda}B)\big)\Big)\cong A\widehat{\otimes}_{\Lambda}B.$ 
		
		(ii) This is immediate  in  the light of and art of \cite[Tag 05GH]{Stacks} in conjunction with the previous part.
		
		(iii) This also follows easily from part (i). Namely, since	 $\mathfrak{M}$ is a maximal ideal so the inverse systems $\{A\otimes_{\Lambda}B/\mathfrak{M}^n\}_{n\in \mathbb{N}}$ and $\{(A\otimes_{\Lambda}B)_{\mathfrak{M}}/(\mathfrak{M})^n\}$ are $(A\otimes_{\Lambda}B)$-algebra isomorphic, because  the ring homomorphism $$A\otimes_{\Lambda}B/\mathfrak{M}^n\rightarrow (A\otimes_{\Lambda}B)_{\mathfrak{M}}/(\mathfrak{M})^n,\ \ \ x+\mathfrak{M}^n\mapsto x/1+(\mathfrak{M}^n)\ \ (\forall\ x\in A\otimes_{\Lambda}B)$$ is an  isomorphism for each $n\in \mathbb{N}$.
	\end{proof}
\end{subprop} 

The concept of generalized local rings and some of their properties  will be used in the sequel. 

\begin{subdefi-rem}\label{NoetherianCompletion}\emph{(\cite[Definition, page 56]{CohenOnTheStructure})
		A quasi-local ring $R$ is said to be  a \textit{generalized local ring} if  $\mathfrak{m}_R$ is finitely generated and  $\bigcap\limits_{n\in\mathbb{N}}\mathfrak{m}^n=0$. By \cite[Theorem 2, page 59]{CohenOnTheStructure} the  $\mathfrak{m}_R$-adic completion,   $\widehat{R},$ of a generalized local ring $R$  is a generalized local ring with maximal ideal  $\mathfrak{m}\widehat{R}$. Even more,  $\widehat{R}$ is  local (Noetherian) by \cite[Theorem 3, page 61]{CohenOnTheStructure}.}
\end{subdefi-rem}

\begin{subrem}
	\emph{Note  that neither   $\mathfrak{M}\subseteq A\otimes_{\Lambda}B$ nor $\mathfrak{M}^e\subseteq A\widehat{\otimes}_{\Lambda}B$ is  necessarily a maximal ideal. For instance,  let $\Lambda=\mathbb{F}_2$ be the field with  $2$ elements and $A=B=\mathbb{F}_2[X]/(X^2+X+1)$ be its field extension. Then $A\widehat{\otimes}_{\Lambda}B=A\otimes_{\Lambda}B\cong \big(\mathbb{F}_2[X]/(X^2+X+1)\big)[Y]/(Y^2+Y+1)$ has two distinct maximal ideals, as the polynomial $Y^2+Y+1$ has two distinct roots $X$ and $X+1$ over the field $\mathbb{F}_2[X]/(X^2+X+1)$.}
\end{subrem}

\begin{subprop}\label{IdealIsMaximal}
	Let $\Lambda$ be a local ring and  $A,B$ be $\Lambda$-local algebras. Assume that $\mathfrak{M}\in \Max (A\otimes_{\Lambda}B)$. 
	\begin{enumerate}
		\item[(i)] $A\widehat{\otimes}_{\Lambda} B$ is a complete local ring with  maximal ideal $\mathfrak{M}^e$.
		\item[(ii)] The  natural maps $A\otimes_{\Lambda}B/\mathfrak{M}^n\rightarrow A\widehat{\otimes}_{\Lambda}B/{\mathfrak{M}^e}^n$ (induced by $\eta_{A\otimes_{\Lambda}B}:A\otimes_{\Lambda}B\rightarrow A\widehat{\otimes}_{\Lambda}B$) are  isomorphisms for all $n\in \mathbb{N}$, and they yield the isomorphism of inverse systems $$\{(A\otimes_{\Lambda}B)/\mathfrak{M}^n\}_{n\in\mathbb{N}}\rightarrow \{(A\widehat{\otimes}_{\Lambda}B)/{{\mathfrak{M}^e}^n}\}_{n\in\mathbb{N}}$$ whose  transition maps are natural epimorphisms.
		\item[(iii)] Suppose that $\mathfrak{a}$ and $\mathfrak{b}$ are ideals of $A$ and $B$, respectively.  ‌Then we have an isomorphism of $(A/\mathfrak{a})$-algebras (resp.  $(B/\mathfrak{b})$-algebras)    $$(A\widehat{\otimes}_{\Lambda}B)/\big(\mathfrak{a}(A\widehat{\otimes}_{\Lambda}B)+\mathfrak{b}(A\widehat{\otimes}_{\Lambda}B)\big)\cong (A/\mathfrak{a})\widehat{\otimes}_{\Lambda}(B/\mathfrak{b}),$$ where the algebra structure of the left-hand side is induced by $\eta_{A,A\widehat{\otimes}B}$ (resp. $\eta_{B,A\widehat{\otimes}_{\Lambda}B}$).
	\end{enumerate}
	\begin{proof}   	  (i)  $(A\otimes_{\Lambda}B)_{\mathfrak{M}}$ is a quasi-local ring whose maximal ideal is finitely generated. Thus
		$$R:=(A\otimes_{\Lambda}B)_{\mathfrak{M}}/\bigcap\limits_{n\in\mathbb{N}}\Big(\big(\mathfrak{M}
		(A\otimes_{\Lambda}B)_{\mathfrak{M}}\big)^n\Big)$$ 
		is a generalized local ring.   The localization map yields the natural isomorphism of inverse systems $\{(A\otimes_\Lambda B)/\mathfrak{M}^n\}_{n\in\mn}\cong \{R/\mathfrak{M}^nR\}_{n\in\mn},$ because $\mathfrak{M}\in \Max(A\otimes_{\Lambda}B)$. Thus, $A\widehat{\otimes}_{\Lambda}B$ is  also $(A\otimes_\Lambda B)$-algebra isomorphic to the   $(\mathfrak{M}R)$-adic completion, $\widehat{R}$, of $R$ by Proposition \ref{CompleteTensorProductNoetherianness}(i).
		In view of Definition and Remark \ref{NoetherianCompletion}, $\widehat{R}$ is a local ring with maximal ideal $\mam_R\widehat{R}=\mathfrak{M}\widehat{R}$. Hence, $A\widehat{\otimes}_{\Lambda}B$ is a local ring with  maximal ideal $\mathfrak{M}^e$ and it is complete in view of Proposition \ref{CompleteTensorProductNoetherianness}(ii).

		(ii)  By Proposition \ref{CompleteTensorProductNoetherianness}(i)  we may, and we do, assume that $A\widehat{\otimes}_{\Lambda}B$ is the $\mathfrak{M}$-adic completion of $A\otimes_{\Lambda}B$. Then the statement is well-known, but we include its proof for the reader's convenience.
		
		Fix some $n\in \mathbb{N}$. If we equip $A\otimes_{\Lambda}B$ with the $\mathfrak{M}$-adic topology, then the completion of $\mathfrak{M}^n$ with respect to the induced topology of $A\otimes_{\Lambda}B$ is   $\lim\limits_{\underset{m\in \mathbb{N}}{\longleftarrow}}\ \big(\mathfrak{M}^n/(\mathfrak{M}^{n}\cap \mathfrak{M}^{m})\big)$  which we denote it by $\widehat{{\mathfrak{M}^n}^{\text{in.}}}$. By \cite[Corollary 10.3]{AtiyahIntroduction}, we already have the natural isomorphism $A\widehat{\otimes}_{\Lambda}B/\widehat{{\mathfrak{M}^n}^{\text{in.}}}\overset{}{\rightarrow }\widehat{A\otimes_{\Lambda}B/\mathfrak{M}^n}$, while   the projection map onto the $n$-th component  yields the isomorphism $\widehat{A\otimes_{\Lambda}B/\mathfrak{M}^n}\rightarrow A\otimes_{\Lambda}B/\mathfrak{M}^n$. Let $f_n:A\widehat{\otimes}_{\Lambda}B/\widehat{{\mathfrak{M}^n}^{\text{in.}}}\rightarrow A\otimes_{\Lambda}B/\mathfrak{M}^n$ be  the composition of  these two isomorphisms, that is given by $(x_m+\mathfrak{M}^m)_{m\in \mathbb{N}}+\widehat{{\mathfrak{M}^n}^{\text{in.}}}\mapsto x_n+\mathfrak{M}^n\ \ (x_m\in A\otimes_{\Lambda}B$ for each $m\in \mn$). 
		
		Let $h_n:A\otimes_{\Lambda}B/\mathfrak{M}^n\rightarrow A\widehat{\otimes}_{\Lambda}B/\widehat{{\mathfrak{M}^n}^{\text{in.}}}$ be the induced map by the canonical map to the completion. Clearly $f_n\circ h_n$ is the identity map, hence $h_n$ is isomorphism too. Thus to finish the proof we only need to show that $\mathfrak{M}^n(A\widehat{\otimes}_{\Lambda}B)=\widehat{{\mathfrak{M}^n}^{\text{in.}}}$ for each $n\in \mn$ (it is clear that then the isomorphisms $h_n$ ($n\in\mathbb{N}$) provide us with the desired isomorphism of the inverse systems as in the statement).
		
		As $\mathfrak{M}^n$ is  finitely generated,  the  map $\mathfrak{M}^n\otimes_{A\otimes_{\Lambda}B}(A\widehat{\otimes}_{\Lambda}B)\rightarrow \widehat{\mathfrak{M}^n}$  given by $y\otimes (x_m+\mathfrak{M}^m)_{m\in \mn}\mapsto (yx_m+\mathfrak{M}^{m+n})_{m\in \mn}$ is surjective (see \cite[Proposition 10.13]{AtiyahIntroduction}). Moreover, by Lemma \ref{cofinal}(ii) we have the  isomorphism $\widehat{\mathfrak{M}^n}\rightarrow \widehat{{\mathfrak{M}^n}^{\text{in.}}}$. Let $\psi:\mathfrak{M}^n\otimes_{A\otimes_{\Lambda}B}(A\widehat{\otimes}_{\Lambda}B)\rightarrow \widehat{{\mathfrak{M}^n}^{\text{in.}}}$ be the composition of these two surjections which fits in the commutative diagram
		\begin{center}
			$\begin{CD}
			\mathfrak{M}^n\otimes_{A\otimes_{\Lambda}B}(A\widehat{\otimes}_{\Lambda}B) @>\psi>\text{surjective} > \widehat{{\mathfrak{M}^n}^{\text{in.}}}\\
			@VVV @VVV\\
			\mathfrak{M}^n(A\widehat{\otimes}_{\Lambda}B) @>>\text{inclusion}> A\widehat{\otimes}_{\Lambda}B
			\end{CD}$.
		\end{center}
		The commutativity of this diagram in conjunction with the surjectivity of $\psi$ shows that $\widehat{{\mathfrak{M}^n}^{\text{in.}}}=\mathfrak{M}^n(A\widehat{\otimes}_{\Lambda}B)$, as required.

		(iii) By symmetry, it suffices to prove that they are  $(A/\mathfrak{a})$-algebra isomorphic.     From part (ii) and the Third Isomorphism Theorem we get the isomorphism 
		$$\psi:\lim\limits_{\underset{m\in \mathbb{N}}{\longleftarrow}}\ 
		\Big((A\otimes_{\Lambda}B)/\big({\mathfrak{M}}^m+ \mathfrak{a}(A\otimes_{\Lambda}B)+\mathfrak{b}(A\otimes_{\Lambda}B)\big)\Big)
		\rightarrow
		\lim\limits_{\underset{m\in \mathbb{N}}{\longleftarrow}}\ \Big((A\widehat{\otimes}_{\Lambda}B)/\big({\mathfrak{M}^{e}}^m+ \mathfrak{a}(A\widehat{\otimes}_{\Lambda}B)+\mab(A\widehat{\otimes}_\Lambda B)\big)\Big)$$
		that fits in the  commutative diagram 
		\begin{center}
			$\begin{CD}
			A @>a\mapsto \Big((a\otimes 1)+\big(\mathfrak{M}^m+(\maa)+(\mab)\big)\Big)_{m\in\mn}>> \lim\limits_{\underset{m\in \mathbb{N}}{\longleftarrow}}\ 
			\Big((A\otimes_{\Lambda}B)/\big({\mathfrak{M}}^m+ \mathfrak{a}(A\otimes_{\Lambda}B)+\mathfrak{b}(A\otimes_{\Lambda}B)\big)\Big)\\
			@V V \text{(nat. epi.)\ }\circ\  \eta_{A,A\widehat{\otimes}_\Lambda B} V @V\cong V\psi \text{\ (indu. by }\ \eta_{A\otimes_{\Lambda}B}) V\\
			\frac{A\widehat{\otimes}_{\Lambda}B}{\mathfrak{a}( A\widehat{\otimes}_{\Lambda}B)+\mathfrak{b}( A\widehat{\otimes}_{\Lambda}B)}
			@>\cong >\text{can. map to the completion}> \lim\limits_{\underset{m\in \mathbb{N}}{\longleftarrow}}\ \Big((A\widehat{\otimes}_{\Lambda}B)/\big({\mathfrak{M}^{e}}^m+ \mathfrak{a}(A\widehat{\otimes}_{\Lambda}B)+\mab(A\widehat{\otimes}_\Lambda B)\big)\Big)
			\end{CD} $.
		\end{center}
		Note that the bottom horizontal map is an isomorphism, because  $(A\widehat{\otimes}_{\Lambda}B)/\big(\mathfrak{a}( A\widehat{\otimes}_{\Lambda}B)+\mathfrak{b}( A\widehat{\otimes}_{\Lambda}B)\big)$ is    a quotient of the complete local ring $A\widehat{\otimes}_\Lambda B$ and so it is also complete. It follows that,
		\begin{equation}
		\label{EhsanWorkHarder}
		A\widehat{\otimes}_\Lambda B/\big(\maa(A\widehat{\otimes}_\Lambda B)+\mab(A\widehat{\otimes}_\Lambda B)\big)
		\cong
		\lim\limits_{\underset{m\in \mathbb{N}}{\longleftarrow}}\ \Big((A\otimes_{\Lambda}B)/\big({\mathfrak{M}}^m+ \mathfrak{a}(A\otimes_{\Lambda}B)+\mab(A\otimes_\Lambda B)\big)\Big)
		\end{equation}
		as $(A/\mathfrak{a})$-algebras where the algebra structure of the left-hand side is induced by  $\eta_{A,A\widehat{\otimes}_{\Lambda}B}$ and the algebra structure of the right-hand side is induced by   the top horizontal map in the diagram.
		
		Also, from  Fact \ref{ElementaryTensorProductFact}(ii) in conjunction with the Third Isomorphism Theorem we obtain the isomorphism $\varphi$ as the bottom horizontal map of the  commutative diagram
		\begin{center}
			$\begin{CD}
			A @>f:=\ a\ \mapsto\  (a+\maa)\otimes (1+\mab)>> (A/\maa)\otimes_\Lambda (B/\mab)\\
			@V a\mapsto \Big((a\otimes 1)+\big(\mathfrak{M}^m+(\maa)+(\mab)\big)\Big)_{m\in\mn} VV @Vg:=\ \text{can. map to the completion} V V\\
			\lim\limits_{\underset{m\in \mathbb{N}}{\longleftarrow}}\ \Big((A\otimes_{\Lambda}B)/\big({\mathfrak{M}}^m+ \mathfrak{a}(A\otimes_{\Lambda}B)+\mathfrak{b}(A\otimes_{\Lambda}B)\big)\Big) @>\cong >\varphi>    \lim\limits_{\underset{m\in \mathbb{N}}{\longleftarrow}}\  \Big((\frac{A}{\maa}\otimes_{\Lambda}\frac{B}{\mab})/\big({\mathfrak{M}}^m(\frac{A}{\maa}\otimes_{\Lambda}\frac{B}{\mab})\big)\Big)
			\end{CD}$.
		\end{center}
		It follows that $\varphi$ is also an isomorphism of $(A/\maa)$-algebras, such that the algebra structure of the target of $\varphi$ is induced by $g\circ f$ as defined in the commutative diagram. From this, (\ref{EhsanWorkHarder}) and  the fact that $(A/\maa)\widehat{\otimes}_{\Lambda}(B/\mab)$ is $\big((A/\maa)\otimes_{\Lambda}(B/\mab)\big)$-algebra isomorphic to the completion of $(A/\maa)\otimes_{\Lambda}(B/\mab)$ (Proposition \ref{CompleteTensorProductNoetherianness}(i)), we deduce the statement.
	\end{proof} 
\end{subprop}

The next result concerns the associativity of certain complete tensor products and it will be used in the next section. Note that by Proposition \ref{IdealIsMaximal}(i), $A\widehat{\otimes}_{\Lambda}B$ (resp.  $B\widehat{\otimes}_{\Gamma}C$) is a complete local  $B$-algebra and thus $\Gamma$-algebra (resp. $\Lambda$-algebra).

\begin{subprop}\label{Associativity} Suppose that $\Lambda$ and $\Gamma$ are  local rings, $A$ is an $\Lambda$-local algebra, $B$ is an $(\Lambda,\Gamma)$-local algebra and  $C$ is an $\Gamma$-local algebra.  If $\mathfrak{M}_{A,\Lambda,B}\in \Max(A\otimes_{\Lambda}B)$ and $\mathfrak{M}_{B,\Gamma,C}\in \Max(B\otimes_{\Gamma}C)$  then  $$A\widehat{\otimes}_{\Lambda}(B\widehat{\otimes}_{\Gamma}C)\cong (A\widehat{\otimes}_{\Lambda}B)\widehat{\otimes}_{\Gamma}C$$ 
	as $A$-algebras (where the left hand side is an $A$-algebra via $\eta_{A,A\widehat{\otimes}_{\Lambda}(B\widehat{\otimes}_{\Gamma}C)}$, and the right hand side is an $A$-algebra via $\eta_{A\widehat{\otimes}_{\Lambda}B,(A\widehat{\otimes}_{\Lambda}B)\widehat{\otimes}_{\Gamma}C}\circ \eta_{A,A\widehat{\otimes}_{\Lambda}B}$.).
	\begin{proof}
		Note that by our hypothesis $B$ is an $(\Lambda,\Gamma)$-bimodule (Convention \ref{FirstConvention}). So the string of isomorphisms
		\begin{align*}
		A\widehat{\otimes}_{\Lambda}(B\widehat{\otimes}_{\Gamma}C)
		&\overset{}{=} 
		\hspace{0.4mm} 	 \lim\limits_{\underset{n\in\mathbb{N}}{\longleftarrow}}\ \Big((A/\mathfrak{m}_A^n)\otimes_{\Lambda}\big((B\widehat{\otimes}_{\Gamma}C)/(\mathfrak{M}^{e}_{B,\Gamma,C})^n\big)\Big) && (\text{by definition})
		&\\&\overset{}{\cong}
		\hspace{0.2mm} \lim\limits_{\underset{n\in\mathbb{N}}{\longleftarrow}}\ \Big((A/\mathfrak{m}_A^n)\otimes_{\Lambda}\big((B\otimes_{\Gamma}C)/(\mathfrak{M}_{B,\Gamma,C})^n\big)\Big) && (\text{Poposition \ref{IdealIsMaximal}(ii)})
		&\\&\overset{}{\cong}
		\hspace{0.10mm} \lim\limits_{\underset{n\in\mathbb{N}}{\longleftarrow}}\ \bigg((A/\mathfrak{m}_A^n)\otimes_{\Lambda}\Big((B\otimes_{\Gamma}C)/\big(\mathfrak{m}^n_{B}(B\otimes_{\Gamma}C)+\mathfrak{m}^n_{C}(B\otimes_{\Gamma}C)\big)\Big)\bigg) && (\text{Lemma \ref{cofinal}(i)})
		&\\&\overset{}{\cong}
		\hspace{-0.3mm} \lim\limits_{\underset{n\in\mathbb{N}}{\longleftarrow}}\ \Big((A/\mathfrak{m}_A^n)\otimes_{\Lambda}\big((B/\mathfrak{m}_B^n)\otimes_{\Gamma}(C/\mathfrak{m}_C^n)\big)\Big) && (\text{Fact \ref{ElementaryTensorProductFact}(ii)})
		&\\&\overset{}{\cong}
		\hspace{-0.3mm} \lim\limits_{\underset{n\in\mathbb{N}}{\longleftarrow}}\ \Big(\big((A/\mathfrak{m}_A^n)\otimes_{\Lambda}(B/\mathfrak{m}_B^n)\big)\otimes_{\Gamma}(C/\mathfrak{m}_C^n)\Big) && (\text{\cite[Exercise 2.15]{AtiyahIntroduction}})
		&\\&\overset{}{\cong }
		\lim\limits_{\underset{n\in\mathbb{N}}{\longleftarrow}}\ \Big(\big((A\otimes_{\Lambda}B)/(\mathfrak{M}_{A,\Lambda,B})^n\big)\otimes_{\Gamma}(C/\mathfrak{m}_C^n)\Big) && (\text{Fact \ref{ElementaryTensorProductFact}+Lemma \ref{cofinal}(i)})
		&\\&\overset{}{\cong}
		\hspace{-0.4mm} \lim\limits_{\underset{n\in\mathbb{N}}{\longleftarrow}}\ \Big(\big((A\widehat{\otimes}_{\Lambda}B)/(\mathfrak{M}_{A,\Lambda,B}^e)^n\big)\otimes_{\Gamma}(C/\mathfrak{m}_C^n)\Big) && (\text{Proposition \ref{IdealIsMaximal}(ii)})
		&\\&\overset{}{=}
		(A\widehat{\otimes}_{\Lambda}B)\widehat{\otimes}_{\Gamma}C && (\text{by definition})
		\end{align*}
		proves our claim. Note that, by a straightforward verification we can see that   the obtained isomorphism takes   $\Big((a+\fm_A^m)\otimes \big(1_{B\widehat{\otimes}_{\Gamma}C}+(\mathfrak{M}_{B,\Gamma,C}^e)^m\big)\Big)_{m\in \mathbb{N}}$ to  $$\bigg(\Big(\big((a+\fm_A^n)\otimes (1_B+\fm_B^n)\big)_{n\in \mathbb{N}}+(\mathfrak{M}_{A,\Lambda,B}^e)^m\Big)\otimes (1_{C}+\fm_C^m)\bigg)_{m\in \mathbb{N}}$$ for each $a\in A$. Consequently, the obtained isomorphism is an isomorphism of $A$-algebras with the mentioned $A$-algebra structures.
	\end{proof}
\end{subprop}

\begin{subrem}
	\emph{It is perhaps worth pointing out that completion of a quasi-local ring   is not necessarily a flat extension, if we drop the Noetherian condition. Namely, in view of   Nagata's  example \cite[Appendix (2), page 69]{Nagata} there exists a non-Noetherian (non-complete) generalized local ring $R$. In the reference it is written only that the example $R$ is  quasi-local with finitely generated maximal ideal, but $R$ is also $\mathfrak{m}_R$-adically  separated as it is a subring of $k[[x,y]]$ over which $k[[x,y]]$ is a local algebra. Therefore, $\widehat{R}$ is not flat over $R$. Because otherwise, the ring homomorphism  $R\rightarrow \widehat{R}$ which maps  $\mathfrak{m}_R$ into  $\mathfrak{m}_{\widehat{R}}$ (see Definition and Remark \ref{NoetherianCompletion}) is a faithfully flat ring homomorphism by \cite[Theorem 7.2]{Matsumura}. Moreover, $\widehat{R}$ is Noetherian (see Definition and Remark \ref{NoetherianCompletion}). These two facts in conjunction with  \cite[Exercise 7.9, page 53]{Matsumura} violate the non-Noetherian property of $R$ and we get a contradiction.}
\end{subrem}

The next proposition provides a sufficient condition for the flatness of $A\widehat{\otimes}_{\Lambda}B$ over $A\otimes_{\Lambda}B$ or  over $A$. In its proof we use the following theorem due to Ogoma.

\begin{subthm}\label{Ogoma}(\cite[Theorem 1]{Ogoma}) Let $\{(A_\lambda,\mathfrak{m}_\lambda)\}$ be a filtered inductive system of Noetherian local rings such that $\mathfrak{m}_{\lambda}A_{\mu}=\mathfrak{m}_\mu$ for $\mu\ge \lambda$. Then the inductive limit $A$ of the  system is Noetherian.	
\end{subthm}

\begin{subprop}\label{FlatProposition} Let   $\Lambda$ be a local ring and  $A,B$ be  $\Lambda$-local algebras. Assume that $\mathfrak{M}\in \Max(A\otimes_{\Lambda}B)$.
	\begin{enumerate}
		\item[(i)] If    $B$  is an Artinian ring, then   $(A\otimes_{\Lambda}B)_\mathfrak{M}$  is Noetherian. In particular, in this case $A\widehat{\otimes}_{\Lambda}B$ is a flat $(A\otimes_{\Lambda}B)$-algebra (via $\eta_{A\otimes_{\Lambda}B}$).
		
		\item[(ii)]  If  $B$ is a  flat $\Lambda$-algebra then $A\widehat{\otimes}_{\Lambda}B$ is a flat (local)  $A$-algebra (via $\eta_{A,A\widehat{\otimes}_{\Lambda}B}$).
		
	\end{enumerate} 
	\begin{proof} (i)  
		We begin by showing that $(A\otimes_{\Lambda}B)_{\mathfrak{M}}$ is Noetherian. 	Let us  prove our claim first in the case where  $B=K_B$  is a field.  Note that in this case, $\mathfrak{M}=\mathfrak{m}_A(A\otimes_{\Lambda}B)$.
		By our hypothesis   $$K_A\otimes_{K_\Lambda}B\cong K_A\otimes_{K_\Lambda}(K_{\Lambda}\otimes_{\Lambda}B)\cong K_A\otimes_\Lambda B\overset{\text{Fact \ref{ElementaryTensorProductFact}(ii)}}{\cong} (A\otimes_{\Lambda}B)/\mathfrak{M}$$ is a field. Thus, Fact \ref{ApplyingSharpResult} implies that  $\mathfrak{m}_A\big((A/\mathfrak{m}_{\Lambda}A)\otimes_{K_\Lambda}\mathcal{L}\big)\in \Max\big((A/\mathfrak{m}_{\Lambda}A)\otimes_{K_\Lambda}\mathcal{L}\big)$ for any subfield  $\mathcal{L}$ of $B$ containing $K_\Lambda$.

		Let $I$ be the directed set consisting of all subfields $\mathcal{L}$ of $B$ such that $\mathcal{L}$ is a finitely generated field extension of $K_\Lambda$ ($I$ is partially ordered with respect to the inclusion). We have $B=\lim\limits_{\underset{L\in I}{\longrightarrow}}\ \mathcal{L}$. So $(A/\mathfrak{m}_\Lambda A)\otimes_{K_\Lambda}B\cong \lim\limits_{\underset{\mathcal{L}\in I}{\longrightarrow}}\big((A/\mathfrak{m}_\Lambda A)\otimes_{K_\Lambda}\mathcal{L}\big)$, as  tensor product commutes with direct limits (\cite[Exercise 20, page 33]{AtiyahIntroduction}). Let, $$\varphi_\mathcal{L}:(A/\mathfrak{m}_\Lambda A)\otimes_{K_\Lambda}\mathcal{L}\overset{\id_{(A/\mathfrak{m}_\Lambda A)}\otimes (\mathcal{L}\hookrightarrow B)}{\longrightarrow} (A/\mathfrak{m}_\Lambda A)\otimes_{K_\Lambda}B$$  be the natural $K_\Lambda$-algebra homomorphism, that is the map to the direct limit ($\mathcal{L}\in I$). As mentioned in the previous paragraph $\mathfrak{m}_A\big((A/\mathfrak{m}_{\Lambda}A)\otimes_{K_\Lambda}\mathcal{L}\big)\in \Max\big((A/\mathfrak{m}_{\Lambda}A)\otimes_{K_\Lambda}\mathcal{L}\big)$, so we have $$\mathfrak{m}_A\big((A/\mathfrak{m}_{\Lambda}A)\otimes_{K_\Lambda}\mathcal{L}\big)=\varphi_\mathcal{L}^{-1}\Big(\mathfrak{m}_A\big((A/\mathfrak{m}_{\Lambda}A)\otimes_{K_\Lambda}B\big)\Big).$$

		From this as well as Fact \ref{DirectLimit} we obtain  
		\begin{equation}		   
		\label{LocalizationAndDirectLimit}
		\big((A/\mathfrak{m}_\Lambda A)\otimes_{K_\Lambda}B\big)_{\mathfrak{m}_A\big((A/\mathfrak{m}_\Lambda A)\otimes_{K_\Lambda}B\big)}
		\cong 
		\lim\limits_{\underset{\mathcal{L}\in I}{\longrightarrow}}\ \Big(\big((A/\mathfrak{m}_\Lambda A)\otimes_{K_\Lambda}\mathcal{L}\big)_{\mathfrak{m}_A\big((A/\mathfrak{m}_\Lambda A)\otimes_{K_\Lambda}\mathcal{L}\big)}\Big).
		\end{equation} 
		Since each $\mathcal{L}\in I$ is a finitely generated field extension of $K_{\Lambda}$, so  each $\big((A/\mathfrak{m}_\Lambda A)\otimes_{K_\Lambda}\mathcal{L}\big)_{\mathfrak{m}_A\big((A/\mathfrak{m}_\Lambda A)\otimes_{K_\Lambda}\mathcal{L}\big)}$ is essentially of finite type  over the Noetherian ring $(A/\mathfrak{m}_\Lambda A)$, and  is therefore Noetherian. From this fact, (\ref{LocalizationAndDirectLimit}) and Theorem \ref{Ogoma}   we  conclude that 
		$\big((A/\mathfrak{m}_\Lambda A)\otimes_{K_\Lambda}B\big)_{\mathfrak{m}_A\big((A/\mathfrak{m}_\Lambda A)\otimes_{K_\Lambda}B\big)}$ is Noetherian. Thus, $$(A\otimes_{\Lambda}B)_\mathfrak{M}\cong \big((A/\mathfrak{m}_\Lambda A)\otimes_{K_\Lambda}B\big)_{\mathfrak{m}_A\big((A/\mathfrak{m}_\Lambda A)\otimes_{K_\Lambda}B\big)}$$ is Noetherian, as claimed ($(A/\mam_\Lambda A)\otimes_{K_\Lambda}B\cong (A\otimes_{\Lambda}K_{\Lambda})\otimes_{K_\Lambda}B\cong A\otimes_{\Lambda}B$).
		
		More generally, suppose that $B$ is an Artinian ring, i.e. $\mathfrak{m}_B$ is a nilpotent ideal. By what we have proved for the case  where $B$ is a field  
		$$(A\otimes_\Lambda B)_\mathfrak{M}/\big(\mam_B(A\otimes_\Lambda B)_\mathfrak{M}\big)
		\cong 
		\Big((A\otimes_\Lambda B)/\big(\mam_B(A\otimes_\Lambda B)\big)\Big)_{(\mathfrak{m}_A)}
		\overset{\text{ Fact \ref{ElementaryTensorProductFact}(ii)}}{\cong} 
		(A\otimes_\Lambda K_B)_{\mam_A(A\otimes_\Lambda K_B)}$$ 
		is Noetherian. From this and the fact that each prime ideal of $(A\otimes_\Lambda B)_\mathfrak{M}$  contains the finitely generated nilpotent ideal $\mam_B(A\otimes_\Lambda B)_\mathfrak{M}$, we conclude that prime ideals of $(A\otimes_\Lambda B)_\mathfrak{M}$ are all finitely generated. Consequently,  $(A\otimes_\Lambda B)_\mathfrak{M}$ is Noetherian by   \cite[Theorem 3.4]{Matsumura}).
		
		By Proposition \ref{CompleteTensorProductNoetherianness}(iii), $A\widehat{\otimes}_{\Lambda}B$  is $(A\otimes_{\Lambda}B)$-algebra isomorphic to the    	  $\big(\mathfrak{M}(A\otimes_{\Lambda}B)_{\mathfrak{M}}\big)$-adic completion, $\widehat{(A\otimes_{\Lambda}B)_{\mathfrak{M}}},$ of the Noetherian ring $(A\otimes_{\Lambda}B)_{\mathfrak{M}}$. But $\widehat{(A\otimes_{\Lambda}B)_{\mathfrak{M}}}$ is flat over  $(A\otimes_{\Lambda}B)_{\mathfrak{M}}$  by \cite[Theorem 8.14.]{Matsumura}, so it is  a flat $(A\otimes_{\Lambda}B)$-algebra by the transitivity of  flatness. Hence $A\widehat{\otimes}_{\Lambda}B$ is a flat $(A\otimes_{\Lambda}B)$-algebra, as was to be proved.

		(ii) Let $n\in \mn$. By  the previous part,
		$(A/\mam_A^n)\widehat{\otimes}_{\Lambda}B$ is  flat over $(A/\mam_A^n)\otimes_{\Lambda}B$  while $(A/\mam_A^n)\otimes_{\Lambda}B$ is flat over $A/\mam_A^n$ in view of the flatness of $B$ over $\Lambda$. Hence, $\eta_{A/\mam_A^n,A/\mam_A^n\widehat{\otimes}_{\Lambda}B}$ is a flat homomorphism i.e. $(A/\mam_A^n)\widehat{\otimes}_{\Lambda}B$ is a  flat $(A/\mam_A^n)$-algebra.
		
		Note that $A\widehat{\otimes}_{\Lambda}B$  is a local $A$-algebra via $\eta_{A,A\widehat{\otimes}_{\Lambda}B}$ (Proposition \ref{IdealIsMaximal}(i)). Thus, it is easily seen that $A\widehat{\otimes}_{\Lambda}B$  is an $\mam_A$-adically ideal-separated $A$-module in the sense of \cite[Definition, page 174]{Matsumura}. Hence, in the light of		 	 \cite[Theorem 22.3(5)]{Matsumura} it suffices to show that $(A\widehat{\otimes}_{\Lambda}B)/\big(\mathfrak{m}^n_{A}(A\widehat{\otimes}_{\Lambda}B)\big)$  is a flat $(A/\mam_A^n)$-algebra via the ring homomorphism induced by $\eta_{A,A\widehat{\otimes}_{\Lambda}B}$ (for each $n\in \mn$).
		But $$A\widehat{\otimes}_{\Lambda}B/\big(\mathfrak{m}^n_{A}(A\widehat{\otimes}_{\Lambda}B)\big)
		\overset{}{\cong} 
		(A/\mam_A^n)\widehat{\otimes}_{\Lambda}B$$ as $(A/\mam_A^n)$-algebras by Proposition \ref{IdealIsMaximal}(iii).   So the statement follows from the conclusion of the previous paragraph. 	\end{proof}
\end{subprop}

\begin{subrem}\label{RemarkSomeResutlsAreNew}
	\emph{To my knowledge,  some of our results in this section are  new. For instance, part (ii) of Proposition \ref{FlatProposition} is stated in \cite[Lemma 19.7.1.2]{GrothendieckIV} but with another  condition that the residue field   of one of the $\Lambda$-algebras  is a finite $\Lambda$-module (cf. with the maximality condition in our statement). Also part (i) of Proposition \ref{FlatProposition} as well as Proposition \ref{Associativity} are new.}
\end{subrem}

\section{Two theorems concerning test modules}

   The main results of this section are Corollary \ref{TestComplexAndFlatHomomorphism} and Corollary \ref{MainResult}.   
Corollary \ref{TestComplexAndFlatHomomorphism} answers positively Question \ref{CWQuestion},  when the residue field $K_R$ is uncountable. Corollary \ref{MainResult} answers affirmatively Question \ref{CDTQuestion}  for all local rings. 

\subsection{Preliminaries}
We devote the first subsection of this section to   definitions and some (rather) elementary facts that will be used in the subsequent subsections.

Let $R$ be a Noetherian ring. The category of finitely generated $R$-modules is denoted by $\text{mod}(R)$. The derived category of $R$ is denoted by $\mathcal{D}(R)$.  A complex $X$ of $R$-modules is said to be a bounded above (resp. bounded below) complex  if $X_i=0$ for all $i\gg 0$ (resp. $i\ll 0$).  It is called a bounded complex provided it is  bounded below and bounded above.  It  is   called a  homologically bounded  above (resp. homologically bounded below) complex if  $H(X)$ is a bounded above (resp. a bounded below) complex. Similarly, $X$ is said to be a homologically bounded complex provided $H(X)$ is a bounded complex. Finally, a homologically bounded complex $X$    is said to be   homologically finite  provided  $H_i(X)$ is a finitely generated $R$-module for all $i$.  In this paper, by a perfect complex we mean a bounded complex consisting of finitely generated free modules (we stress that, in contrast to our definition in this paper, in the literature perfect complexes are considered as bounded complexes consisting of finitely generated projective modules).

 The full subcategory of $\mathcal{D}(R)$ consisting of homologically bounded complexes is denoted by $\mathcal{D}_b(R)$, and the full subcategory of $\mathcal{D}_b(R)$ consisting of homologically finite complexes is denoted by  $\mathcal{D}_b^f(R)$. Quasi-isomorphism of complexes are denoted by $\simeq$.  

For the definition of a projective resolution  and  a flat resolution of a homologically bounded below complex  see parts (F) and (P) of \cite[(A.3.1) Definitions]{ChrisensenGorenstein}, and for the definition of an injective resolution of a homologically bounded above complex see \cite[(A.3.1) Definitions(I)]{ChrisensenGorenstein}. Moreover, see \cite[(A.3.2) Theorem (Existence of Resolutions)]{ChrisensenGorenstein} for the existence of projective, flat and injective resolutions. Using projective and  injective resolutions, we can define the projective dimension of homologically bounded below complexes and the injective dimension of homologically bounded above complexes (see \cite[(A.3.9) Definition]{ChrisensenGorenstein} and \cite[(A.3.8) Definition]{ChrisensenGorenstein}).

For the definition of the derived tensor product $X\otimes_R^{\mathbf{L}}Y$ where  $X$ or $Y$ is a homologically bounded below complex (which is defined using  flat resolutions of complexes)  we refer to \cite[(A.4.11) Definition]{ChrisensenGorenstein}, and in this case we define $\text{Tor}^R_i(X,Y):=H_i(X\otimes^{\mathbf{L}}_RY)$.  Also, for the definition of the derived  Hom-complex $\text{\textbf{R}Hom}_R(X,Y)$  of two $R$-complexes $X$ and $Y$  where $X$ is homologically bounded below or $Y$ is homologically bounded above (which is defined using a projective resolution of $X$ or an injective resolution of $Y$) we refer to \cite[(A.4.2) Definition]{ChrisensenGorenstein}. If $X$ is a homologically bounded below complex or $Y$ is a homologically bounded above complex we define $$\Ext^i_R(X,Y):=H_{-i}\big(\text{\textbf{R}Hom}_R(X,Y)\big).$$

A bounded homologically finite  complex $D$ consisting  of injective $R$-modules is said to be a dualizing complex for $R$ provided  the natural homothety map $R\rightarrow \text{\textbf{R}Hom}_R(D,D)$  is a quasi-isomorphism (see \cite[(7.1.2) Definition]{ChristensenFoxbyHyperHomological} for the definition of the homothety map). It follows from the definition that any shift of a dualizing complex  is again a dualizing complex (see \cite[(A.1.3) Shift]{ChrisensenGorenstein} for the definition of a shift of a complex). For any complex of $R$-modules we define $X^{\vee}:=\text{Hom}_R(X,D)(\simeq\text{\textbf{R}Hom}_R(X,D))$.

\begin{subdefi}\label{quasiDeformationDefinition}
	\emph{A \textit{deformation} of a Noetherian ring $R$ is a ring epimorphism $\pi:A\twoheadrightarrow R$ of Noetherian rings whose kernel is generated by a regular sequence of $A$. A \textit{quasi-deformation} of  a local ring $R$ is a diagram $R\overset{g}{\hookrightarrow} S\overset{\pi}{\twoheadleftarrow} A$ of local homomorphisms of local rings, where $g$ is a flat local homomorphism and $\pi$ is a deformation of $S$.}
\end{subdefi}

Here we borrow the definition of complete intersection dimension of a homologically finite complex  from \cite{WagstaffComplete}, but only over local rings  which is sufficient for our purpose in this paper (it is defined for arbitrary Noetherian rings).  The notion of the complete intersection dimension for homologically finite complexes extends  the concept of complete intersection dimension for finitely generated  modules introduced in \cite{Avramov} (by considering a module as a complex concentrated at  degree $0$).

\begin{subdefi} \emph{Let $R$ be a local ring and $X$ be a homologically finite complex of $R$-modules, i.e. a complex $X\in \mathcal{D}_b^f(R)$.
			The \textit{complete intersection dimension} of $X$ is defined as,
			$$\text{CI-dim}_R(X):=\inf\{\pd_A(X\otimes_RS)-\pd_A(S):\ R\hookrightarrow S \twoheadleftarrow A\text{\ is a quasi-deformation of \ }R\}.$$}
\end{subdefi}

We stress that the zero complex  (so also the zero module) has complete intersection dimension $-\infty$ by convention. 


\begin{subdefi} \emph{Let $\varphi:R\rightarrow S$ be a local homomorphism of local rings.
  \begin{enumerate}
  	\item [(i)] $\varphi$ is said to be  weakly unramified provided $\mathfrak{m}_RS=\mathfrak{m}_S$, and in this case we say that $S$ is a weakly unramified (local) $R$-algebra.
  	\item[(ii)] $\varphi$ is said to be weakly regular provided it is flat and $S/\fm_RS$ is a regular local ring (compare with the concept of a  regular homomorphism, Definition \ref{DefinitionRegularHomomorphism}).
  \end{enumerate}}
\end{subdefi}

\begin{subdefi-rem}\emph{\label{DefinitionRemarkCohenFactorization}(\cite{AvramovFoxbyHerzogCohen})
 Let $\varphi:R\rightarrow S$ be a local homomorphism of local rings. Let $\varphi^{\backprime}:R\overset{\varphi}{\rightarrow} S\overset{s\mapsto (s+\fm_S^n)_{n\in \mathbb{N}}}{\longrightarrow} \widehat{S}$	be the  induced map to $\widehat{S}$. Then a Cohen factorization of $\varphi^{\backprime}:R\rightarrow \widehat{S}$ is a commutative diagram of local homomorphisms} 
 	\begin{center}\begin{tikzcd}
 		  & R' \arrow[dr,"\varphi'"] \\
 		  R\arrow[rr,"\varphi^{\backprime}" swap]\arrow[ur,"\dot{\varphi}"] && \widehat{S}
 	\end{tikzcd}\end{center}
 \emph{such that $\dot{\varphi}$ is a weakly regular homomorphism, $R'$ is a complete local ring and $\varphi'$ is surjective. By \cite[(1.1) Theorem (Existence of Cohen Factorizations)]{AvramovFoxbyHerzogCohen}, we can assign some  Cohen factorization of  $\varphi^{\backprime}$  to every local homomorphism $\varphi$.}
\end{subdefi-rem}

The following nice criterion for finiteness of the flat dimension generalizes the well-known local flatness criterion of finitely generated modules. 

\begin{subfact}(see \cite[II, Lemma 57]{AndreHomologie})\label{FactAndreCharmingResult} Let $\varphi:A\rightarrow B$ be a homomorphism of Noetherian rings and $M$ be a finitely generated $B$-module. If $\tor^A_{n}\big(A/\varphi^{-1}(\fm),M\big)=0$ for    each  $\fm\in \Max(B)$ and some (fixed) $n\in \mathbb{N}$, then $M$ has  flat dimension $\le n-1$ as an $A$-module. In particular, if $\varphi$ is a local homomorphism of local rings then $\tor^A_{n}(K_A,M)=0$ if and only if $M$ has flat dimension $\le n-1$ as an $A$-module.
\end{subfact}

\begin{subnota}\label{NotationRestrictionOfScalars} 	  \emph{ The notation $M\rceil_{h}$ assigned to a ring homomorphism $h:A\rightarrow B$   and an $B$-module $M$   means the Abelian group $M$ considered as an $A$-module by restricting the scalars via $h$.}
\end{subnota}

\begin{subnota}\emph{The notation $\tor^R_{\gg}(M,N)=0$ assigned to some $R$-modules $M$ and $N$ means $\tor^R_i(M,N)=0$ for $i\gg 0$.}
\end{subnota}

\begin{subdefi}
	\emph{Let $R$ be a Noetherian ring and $T$ be a finitely generated $R$-module. We say that $T$ is a \textit{test module} for $R$, or equivalently a \textit{pd-test module}, if the following condition holds for all $M\in \text{mod}(R)$: If $\tor^R_{\gg}(T,M)=0$ then $\pd_R(M)<\infty$.}
\end{subdefi}

\begin{subdefi}
	\emph{Let $R$ be a Noetherian  ring. A complex $T\in \mathcal{D}^f_b(R)$ is said to be a \textit{test complex} for $R$, or equivalently  a \textit{pd-test complex} for $R$,  if the following condition holds for all $X\in \mathcal{D}^f_b(R)$: If $\text{Tor}^R_i(T,X)=0$ for all $i\gg 0$, i.e. if $T\otimes_R^{\mathbf{L}}X\in \mathcal{D}_b(R)$, then $\pd_R(X)<\infty$.}		
\end{subdefi}

\begin{subdefi-rems}
	\emph{Let $R$ be a quasi-local ring. The integer $\text{Char}(R/\mam_R)$  is said to be the  \textit{residual characteristic}  of   $R$   (which is either zero or a prime number). When the residual characteristic of $R$  agrees with $\text{Char}(R)$ we say that $R$  is an  \textit{equicharacteristic} ring, otherwise we say that $R$ has \textit{mixed characteristic}. 
		It is easily seen that a quasi-local ring is equicharacteristic precisely when it contains a field. We say that a quasi-local ring $R$ has mixed characteristic $(0,p)$ provided $\text{Char}(R)=0$ while $\text{Char}(K_R)=p>0$. A quasi-local  ring $(R,\mam_R)$ with residual characteristic $p>0$ is said to be \textit{unramified}    provided $p\notin \mam_R^2$.   A complete unramified discrete valuation ring of mixed characteristic $(0,p)$ is said to be a \textit{$p$-ring} (thus a $p$-ring is a $1$-dimensional complete regular local ring whose maximal ideal is generated by its residual characteristic, $p>0$).}
\end{subdefi-rems}

\begin{subrem}\label{RemarkExtensionOfUnramifiedDVRs}
	\emph{Suppose that $\varphi:V\rightarrow V'$ is a ring homomorphism where $V$ and $V'$ are both unramified discrete valuation rings of mixed characteristic $(0,p)$. Thus $V$ and $V'$ are regular local rings with $\fm_V=pV$ and $\fm_{V'}=pV'$ (see \cite[Theorem 11.2]{Matsumura}). Then $\varphi$ is a weakly unramified local homomorphism, because indeed $\varphi(p)=p$.  Moreover, $\varphi$ is a flat ring homomorphism because $V'$ is a domain thus it is torsion free (see \cite[Exercise 11.9, page 86]{Matsumura}).} 
\end{subrem}

The following definition of coefficient rings is similar to \cite[Tag 0326]{Stacks}, but in contrast to it our coefficient ring  is  always a domain and it does not need  to be a subring of the ambient ring.

\begin{subdefi}\label{coefficientringdefinition}
	\emph{Let $R$ be a complete local ring. In our paper, a \textit{coefficient ring} of $R$ is a local homomorphism of local rings $\lambda_R:C\rightarrow R$ such that the induced map $\overline{\lambda_R}:C/\mathfrak{m}_{C}\rightarrow R/\mathfrak{m}_R$ is an isomorphism and $\begin{cases}C\text{\ is a field}, &\text{\ if } R\text{\  contains a field}\\ C\text{\ is a $p$-ring,}& \text{\ if }R\text{\ has mixed characteristic and }\text{Char}(K_R)=p>0\end{cases}.$ \emph}
\end{subdefi}

\begin{subrem} \emph{In view of  the Cohen Structure Theorem (see \cite[Tag 032A]{Stacks}) any complete local ring admits a coefficient ring (the term ``Cohen ring" in the statement of  \cite[Tag 032A]{Stacks} means the same $p$-ring as defined in our paper, see \cite[Tag 0327]{Stacks}).}
\end{subrem}



The following    facts and remarks will be used in the sequel. Since we do not know any reference for many of the  facts below, we provide their proof here for the sake of completeness.

\begin{subrem}\label{QuasiDeformationWithCompleteDeformation}
	\emph{Suppose that $R$ is a local ring and that $L$ is an $R$-module of finite complete intersection dimension. When we are considering a quasi-deformation $R\hookrightarrow V\twoheadleftarrow W$   with $\pd_W(L\otimes_RV)<\infty$, without loss of generality, we can assume that $V$ and $W$ are complete local rings. Namely, assume that $V=W/\mathbf{x}W$ for some regular sequence $\mathbf{x}$ for $W$. It is well-known that $\widehat{W}$ (resp. $\widehat{V}$) is a  faithfully flat extension of $W$ (resp. $V$) (\cite[Theorem 8.14]{Matsumura}), that the image of  $\mathbf{x}$ in $\widehat{W}$ is also a regular sequence of $\widehat{W}$ (\cite[Theorem 7.4]{Matsumura}) and that $\widehat{W}/\mathbf{x}\widehat{W}\cong \widehat{V}$. Consequently, $\widehat{V}$ is also a flat extension of $R$ in view of the transitivity of flatness  and $\widehat{W}$ is a deformation of $\widehat{V}$, in other words	 $R\hookrightarrow \widehat{V}\twoheadleftarrow \widehat{W}$ is also a quasi-deformation of $R$. Moreover,  we have the  string  $$(L\otimes_RV)\otimes_W\widehat{W}\cong L\otimes_R(V\otimes_W\widehat{W})\cong L\otimes_R(\widehat{W}/\mathbf{x}\widehat{W})\cong L\otimes_R\widehat{V}$$ of $\widehat{W}$-isomorphisms. Therefore, $\pd_{\widehat{W}}(L\otimes_R\widehat{V})=\pd_{\widehat{W}}\big((L\otimes_RV)\otimes_W\widehat{W}\big)=\pd_W(L\otimes_RV)<\infty$, where the last equality holds because $\widehat{W}$ is faithfully flat over $W$. }
\end{subrem}


\begin{subrem}\emph{\label{RemarkTestModulePropertyDesentFromFaithfullyFlatExtension}
	Suppose that   $R$ is a local ring and $S$ is a flat local $R$-algebra. If $T$  is a finitely generated $R$-module such that $T\otimes_RS$ is a test module for $S$, then it is easily seen that $T$ is a test module for $R$. Namely, if $N$ is a finitely generated $R$-module with $\text{Tor}^R_{\gg 0}(T,N)=0$ then by tensoring the flat $R$-algebra $S$  we get $\text{Tor}^S_{\gg 0}(T\otimes_RS,N\otimes_RS)=0$ implying that $\pd_R(N)=\pd_{S}(N\otimes_RS)<\infty$, as $T\otimes_RS$ is an $S$-test module.}
\end{subrem}

\begin{subfact} \label{ResolutionModuleRegularSequence} (see \cite[Proposition 1.1.5]{BrunsHerzogCohenMacaulay}) Let  $R$ be a ring and  $$N_\bullet:=\cdots\rightarrow N_m\overset{\varphi_m}{\rightarrow }N_{m-1}\rightarrow \cdots \rightarrow  N_0\overset{\varphi_0}{\rightarrow }N_{-1}\rightarrow 0$$ an exact complex of $R$-modules. If $\mathbf{x}$ is weakly $N_i$-regular for all $i$, then $N_\bullet \otimes_RR/(\mathbf{x})$ is exact again.
\end{subfact}

\begin{subfact}\label{SyzygyAndGrade} Assume that $R$ is a ring, $N$ is an $R$-module and $\mathbf{x}:=x_1,\ldots,x_n$ is a sequence  of elements of  $R$  with  $\mathbf{x}N\neq N$.
	\begin{enumerate}
		\item [(i)] Presume that there is an exact sequence $0\rightarrow N\rightarrow M_{n-1}\rightarrow M_{n-2}\rightarrow \cdots \rightarrow M_0\rightarrow L\rightarrow 0$ of $R$-modules such that $\mathbf{x}$ is  weakly $M_i$-regular  for each $0\le i\le n-1$. Then $\mathbf{x}$ is a regular sequence on $N$ (note that $\mathbf{x}$ has length $n$).
		\item[(ii)] If   $N$ is an $n$-syzygy  and $\mathbf{x}$ is a  regular sequence on $R$,  then $\mathbf{x}$ is  a regular sequence on $N$ (note that $\mathbf{x}$ has length $n$).
	\end{enumerate} 
	\begin{proof} 
		(i)  If $n=1$‌ then $N$ ‌is a submodule of $M_0$ while $x_1$ is a weakly $M_0$-regular element, thus $x_1$ is regular on its submodule $N$. Hence the statement holds in this case.
		Suppose that $n>1$ and the statement has been proved for smaller values than $n$. We can take into account the exact sequence, $0\rightarrow N\rightarrow M_{n-1}\rightarrow \cdots\rightarrow M_1\rightarrow L'\rightarrow0$ where $L':=\ker (M_0\rightarrow L)$. Then $x_1$ is weakly $L'$-regular (as it is weakly $M_0$-regular),   and it is weakly $M_i$-regular    for each $1\le i\le n-1$ by our hypothesis. Moreover, $x_1$ is $N$-regular by the base of our induction. So in view of Fact \ref{ResolutionModuleRegularSequence}, the complex $$0\rightarrow N/x_1N\rightarrow M_{n-1}/x_1 M_{n-1}\rightarrow \cdots\rightarrow M_{1}/x_1M_{1}\rightarrow L'/x_1L'\rightarrow 0$$ is exact while $x_2,\ldots,x_n$ is a weakly $(M_{i}/x_1M_{i})$-regular sequence for each $1\le i\le n-1$. So by our inductive hypothesis $x_2,\ldots,x_n$ is a regular sequence on $N/x_1N$, i.e. $x_1,\ldots,x_n$ is a regular sequence on $N$.
		
		(ii) This part is a special case of the previous part, because $N$ being an $n$-syzygy fits into an exact sequence $0\rightarrow N\rightarrow F_{n-1}\rightarrow \cdots\rightarrow F_0\rightarrow L\rightarrow 0$ where $F_0,\ldots,F_{n-1}$ are free $R$-modules.
	\end{proof}
\end{subfact}

\begin{subfact}\label{FactTorOverBisTorOverAWhenBflatOverA}
	Suppose that $\varphi:A\rightarrow B$ is a flat local homomorphism, $M$ is an $A$-module and $N$ is an $B$-module. Then $$\tor^B_i(M\otimes_AB,N)\cong \tor^A_i(M,N\rceil_{\varphi}),\ \ \  \forall\ i\in \mathbb{N}.$$ 
	  \begin{proof}
	  	 Let $F_\bullet$ be a free resolution of $M$ (over $A$). Since $B$ is a flat $A$-module by our assumption, so $F_\bullet\otimes_AB$ is a free resolution of $M\otimes_AB$ (over $B$).  Hence, $$\tor^B_i(M\otimes_AB,N)=H_i\big((F_\bullet\otimes_AB)\otimes_B N\big)\cong H_i\big(F_\bullet\otimes_A(B\otimes_BN)\big)\cong H_i(F_\bullet\otimes_AN\rceil_{\varphi})\cong \tor^A_i(M,N\rceil_{\varphi}).\qedhere$$
	  \end{proof}
\end{subfact}

\begin{subfact}\label{FactTorOverRTorOverRxR}
	Suppose that $R$ is a ring and $\mathbf{x}$ is a regular sequence on $R$ and on some $R$-module $M$. Let $N$ be an $R/xR$-module. Then $\tor^R_i(M,N)\cong \tor^{R/(\mathbf{x})}_i(M/\mathbf{x}M,N)$ for each $i$.
	  \begin{proof}
	  	 Let $F_\bullet$ be a free resolution $M$. Since $\mathbf{x}$ is a regular sequence on both of $M$ and $R$ so $F_\bullet/\mathbf{x}F_\bullet$ is an $(R/\mathbf{x}R)$-free resolution of $M/\mathbf{x}M$ (by Fact \ref{ResolutionModuleRegularSequence}). So 
	  	   \begin{align*}
	  	     \tor^R_i(M,N)\cong H_i(F_\bullet\otimes_RN)\cong H_i\big(F_\bullet\otimes_R(R/\mathbf{x}R\otimes_{R/\mathbf{x}R}N)\big)
	  	     &\cong
	  	     H_i(F_\bullet/\mathbf{x}F_\bullet\otimes_{R/\mathbf{x}R}N)
	  	     &\\&
	  	     \cong\tor^{R/\mathbf{x}R}_i(M/\mathbf{x}M,N). \qedhere
	  	   \end{align*}
	  \end{proof}
\end{subfact}

\begin{subfact}\label{FactGradeOverNiceFaithfullyFlatBaseChange}
	  Suppose that $\varphi:A\rightarrow B$ 	is a faithfully flat ring homomorphism of Noetherian rings and that $\mathbf{x}:=x_1,\ldots,x_n$ is a regular sequence on $B$. If $B/\mathbf{x}B$ is also flat over $A$, then $\text{grade}_{(\mathbf{x})}(B\otimes_AM)=n$ for any non-zero finitely generated $A$-module $M$.
	    \begin{proof}
	    	Consider a finitely generated $A$-module $M$. Then $B\otimes_AM$ is a finitely generated $B$-module. As $\mathbf{x}$ is a regular sequence of $B$ so the Koszul complex $K_\bullet(\mathbf{x};B)$ is acyclic (\cite[Corollary 1.6.14]{BrunsHerzogCohenMacaulay}). Even more, $K_\bullet(\mathbf{x};B)$ is an $A$-flat resolution of $B/\mathbf{x}B(=H_0(\mathbf{x};B))$, because $B$ is assumed to be a flat $A$-module. In particular, $\tor^A(B/\mathbf{x}B,-)$ can be computed using $K_\bullet(\mathbf{x};B)\otimes_A-$ (see \cite[Theorem 7.5]{Rotman}). Therefore,
	    	$$ 
	    	  H_i(\mathbf{x};B\otimes_AM)=H_i\big(K_\bullet(\mathbf{x};B\otimes_AM)\big)\cong H_i\big(K_\bullet(\mathbf{x};B)\otimes_AM\big)\cong \tor^A_i(B/\mathbf{x}B,M)=0
	    	$$
	    	for each $i\gneq 0$, because $B/\mathbf{x}B$ is assumed to be flat over $A$. From this as well as \cite[Corollary 1.6.17]{BrunsHerzogCohenMacaulay} we conclude that $\text{grade}_{(\mathbf{x})}(B\otimes_AM)=n$ ($B\otimes_AM\neq 0$, because $\varphi$ is faithfully flat).
	    \end{proof}
\end{subfact}

\begin{subfact}\label{FactTorOfMAndTorOfSpecializationM}
	Suppose that $A$ is a Noetherian ring and $M,N$ are two $A$-modules such that $N$ is finitely generated. Let $\mathbf{x}:=x_1,\ldots,x_n$ be a sequence of elements of $A$ such that $\text{grade}_{(\mathbf{x})}(N)=n$. Assume that for some $a\in A$ and some $i\ge n$ 
	$$
	\tor^A_j(M,N)_a=0,\ \ \forall\ i-n\le j\le i.
	$$
	Then $\tor^A_i(M,N/\mathbf{x}N)_a=0$.
	
	\begin{proof}
		 By \cite[Exercise 1.2.21]{BrunsHerzogCohenMacaulay} without loss of generality we may, and we do, assume that $\mathbf{x}$ is a regular sequence on $N$. Suppose that $n=1$. Then applying $\tor^A(M,-)$ to the exact sequence 
		 $$
		 0\rightarrow N\rightarrow N\rightarrow N/x_1N\rightarrow 0
		 $$
		 and afterwards applying the exact localization functor $(-)_a$ we obtain the exact sequence 
		 $$
		 \tor^A_i(M,N)_a\overset{x_1}{\rightarrow}\tor^A_i(M,N)_a\rightarrow \tor^A_i(M,N/x_1N)_a\rightarrow \tor^A_{i-1}(M,N)_a
		 $$
		 implying that $\tor^A_i(M,N/x_1N)_a=0$.
		 
		 Now suppose that $n>1$ and the statement has been proved for smaller values than $n$. Repeated use of the base of our induction yields $\tor^A_j(M,N/x_1N)_a=0$ for $i-n+1\le j\le i$. Hence from the induction hypothesis we get our desired vanishing $\tor^A_i(M,N/\mathbf{x}N)_a=0$.
	\end{proof}
\end{subfact}

\begin{subfact}\label{EhsanStuck} Let $A$ be an $\Lambda$-algebra. Let  $S$ be a multiplicatively closed subset of $A$. Suppose that $(F_\bullet,\partial^{F_\bullet}_\bullet)$  is a perfect complex over $S^{-1}(A)$.  Then, there is a perfect complex $(G_\bullet,\partial^{G_\bullet}_\bullet)$ of (finite free) $A$-modules and $A$-homomorphisms such that it fulfills the following properties.
	\begin{enumerate} 
		\item [(i)]  $G_\bullet\otimes_A S^{-1}(A)\cong F_\bullet$.
		\item[(ii)] If $\Lambda$ is local and the entries of the matrix of $\partial^{F_\bullet}_i$ belong to $\mam_{\Lambda}S^{-1}(A)$ for every $i$ then the entries of the matrix of $\partial^{G_\bullet}_i$ belong to $\mam_{\Lambda}A$ for every $i$.
	\end{enumerate}	
	\begin{proof}
	  	Suppose that $$F_\bullet:=0\rightarrow S^{-1}(A)^{n_u}\overset{H_u}{\longrightarrow} S^{-1}(A)^{n_{u-1}}\overset{H_{u-1}}{\longrightarrow}\cdots\overset{H_2}{\longrightarrow}  S^{-1}(A)^{n_1}\overset{H_1}{\longrightarrow} S^{-1}(A)^{n_0}\rightarrow 0$$ with differentials given by the matrices $H_i$ $(1\le i\le u$). 
	  
	  Let $H_i=[\ {_ia_{l,k}}/{_it_{l,k}}\ ]$ for each $1\le i\le u$. If the entries of the matrices $H_i$ belong to $\mam_\Lambda S^{-1}(A)$, then it is easily seen that we can present $H_i$ as $H_i=[\ {_ia_{l,k}}/{_it_{l,k}}\ ]$ with $_ia_{l,k}\in \mam_{\Lambda}A$ for each $l,k$. Let $s:=\prod\limits_{\text{all possible}\ i,l,k} {_it_{l,k}}$ and set ${_ia_{l,k}'}:={_ia_{l,k}}(s/{_it_{l,k}})$. Then we get  $H_i=[\ {_ia_{l,k}}/{_it_{l,k}}\ ]=[\ {_ia'_{l,k}}/s\ ]$ for each $1\le i\le u$. Note that if ${_ia_{l,k}}\in \mam_{\Lambda}A$ then evidently ${_ia'_{l,k}}\in \mam_{\Lambda}A$. It follows that, we can find a fixed element $s\in S$ such that 
	  \begin{equation} 
	  	\label{EquationMatricesWithUnifiedDenominator}
	  	H_i=[\ _ia'_{l,k}/s \ ]
	  \end{equation}   so   the entries of all  $H_i$'s all have the same denominator $s$. Moreover, we have ${_ia'_{l,k}}\in \mam_{\Lambda}A$ for each $l,k$ provided the entries of $H_i$ belong to $\fm_{\Lambda}S^{-1}(A)$. 
	  
	  Set, $$ G_\bullet:=0\rightarrow A^{n_u}\overset{B_u}{\longrightarrow} A^{n_{u-1}}\overset{B_{u-1}}{\longrightarrow} \cdots \overset{B_2}{\longrightarrow} A^{n_1}\overset{B_1}{\longrightarrow} A^{n_0}\rightarrow 0$$ with differentials given by the matrices $B_i=[\ _ia'_{l,k}\ ]$ (thus $sH_i=B_i$). 
	  
	  At this form $G_\bullet$ is not a complex necessarily because we do not have $\partial^{G_\bullet}_{i+1}\circ \partial^{G_\bullet}_i=0$ necessarily (unless $S$ is consisting of regular elements of $A$). Therefore we consider a modification $G'_\bullet$ of $G_\bullet$ consisting of the same finite free modules as in $G_\bullet$ but with the differentials given by the matrices $B'_i=[\ t.\ {_ia'_{l,k}}\ ]$  for some fixed $t\in S$. Choosing this fixed $t$ appropriately makes $G'_\bullet$ a complex, as $G_\bullet$ becomes a complex after localizing at $S$. More precisely, $\im (B_{i+1}\circ B_{i})=0$ for each $i$  after localization at $S$, thus $$\im\big((t.B_{i+1})\circ (t.B_{i})\big)=t^2.\im (B_{i+1}\circ B_{i})\subset t.\im (B_{i+1}\circ B_{i})=0$$
	  for some $t\in S$.

	  For each $0\le i\le u$, let $$\psi_i:A^{n_i}\otimes_A S^{-1}(A)\rightarrow S^{-1}(A)^{n_i},\ \ (x_j)_{j=1}^{n_i}\otimes (a/t')\mapsto \big((x_ja)/t'\big)_{j=1}^{n_i}\ \ \ (a,x_1,\ldots,x_{n_i}\in A,\ t'\in S)$$  be the natural isomorphism. Set, $$\varphi_i:=(ts)^i.\psi_i,\ \ (0\le i\le u)$$  ($s$ is appeared as the unified denominator in (\ref{EquationMatricesWithUnifiedDenominator}), and $t$ is the fixed element used to define $B'_i$). Then, it is straightforward to check that $$(\varphi_i)_{i=0}^u:G'_\bullet\otimes_A S^{-1}(A)\rightarrow F_\bullet$$ is a chain map and an isomorphism of complexes.  This proves part (i), and part (ii) also holds in view of the previous arguments.
	\end{proof}
\end{subfact}

\begin{subfact}\label{FaithfullyFlat}
	Suppose that $\varphi:A\rightarrow B$ is a ring homomorphism and $\mathbf{X}$ is a finite sequence of indeterminates. 
	\begin{enumerate}
		\item [(i)]  $\varphi$‌
		is a faithfully flat ring homomorphism (i.e. $B$‌ is  faithfully flat as $A$-module) if and only if it is injective and $B/A$ is a flat $A$-module.
		\item[(ii)] Suppose that $\varphi$ is faithfully flat (e.g. $\varphi$ is a flat local homomorphism of local rings). Then $A[\mathbf{X}]\overset{}{\longrightarrow} B[\mathbf{X}]$ (given by $a\mapsto \varphi(a),\ X_i\mapsto X_i$) is also a faithfully flat ring homomorphism.
	\end{enumerate}	
	\begin{proof}
		(i)	 Suppose that $B$ is a faithfully flat $A$-module. For the injectivity of $\varphi$ see \cite[Theorem 7.5(i)]{Matsumura}. For the flatness of $B/A$,  we (always) have the split exact sequence $0\rightarrow B\overset{\eta}{\rightarrow} B\otimes_AB\rightarrow B\otimes_A (B/A)\rightarrow 0$, wherein $\eta$ is the composition  $(\text{id}_B \otimes \varphi)\circ (B\overset{\cong}{\rightarrow }B\otimes_AA)$. The splitting property of $\eta$ holds in  view of  the multiplicative map $b\otimes b'\rightarrow bb'$.  If follows that $B\otimes_A (B/A)$‌ is a direct summand of $B\otimes_A B$ and thus it is flat. Now, since $B$‌ is faithfully flat as $A$-module so $B/A$ is also a flat $A$-module.

		For the reverse implication,  from the injectivity of $\varphi$ we have the exact sequence $$\mathcal{E}:0\rightarrow A\rightarrow B\rightarrow B/A\rightarrow 0$$ of $A$-modules and then the flatness of $A$ and $B/A$  implies that $B$‌ is also a flat $A$-module. Let $M$‌ be an $A$-module. In view of the flatness of $B/A$,  $\mathcal{E}$‌ is a  pure exact sequence and  therefore it remains exact after tensoring with $M$ implying that $M\rightarrow M\otimes_AB$ is injective. Therefore, $M\otimes_AB\neq 0$ for any non-zero $A$-module $M$. So $B$ is faithfully flat over $A$ (\cite[Theorem 7.2]{Matsumura}). 
		
		(ii) By part (i) it suffices to show that $B[\mathbf{X}]/A[\mathbf{X}]\cong (B/A)[\mathbf{X}]$ is a flat $A[\mathbf{X}]$-module. Again, by part (i) $B/A$ is a flat $A$-module so it is a direct limit of free modules in light of \cite[Theorem 5.40(Lazard)]{Rotman}. Then it is easily seen that $(B/A)[\mathbf{X}]$ is also a direct limit of free $A[\mathbf{X}]$-modules, i.e. a flat $A[\mathbf{X}]$-module.
	\end{proof}
\end{subfact}

\begin{subfact}\label{FactCompleteTensorProductFactCoefficientRingExtension} 
	Suppose that $R$ is a complete local ring and $\lambda_R:C\rightarrow R$ is a coefficient ring of $R$. Let  $\varphi:C\rightarrow C'$ be a local homomorphism that is  either an extension of   fields  or an extension of unramified discrete valuation rings of mixed characteristic $(0,p)$ ($C$ is indeed a $p$-ring in the mixed characteristic case, but we let $C'$ to be possibly non-complete).
	
	\begin{enumerate}
		\item [(i)] We have $\mathfrak{M}\overset{\text{by definition}}{=}\fm_R(R\otimes_CC')+\fm_{C'}(R\otimes_CC')=\fm_R(R\otimes_CC')$ and $\mathfrak{M}\in \text{Max}(R\otimes_CC')$.
		\item[(ii)] 
		Consider the complete tensor product $R\widehat{\otimes}_{C}C'$ assigned to local $C$-algebras $\varphi$ and $\lambda_R$. Then the ring homomorphism $$\eta_{R,R\widehat{\otimes}_{C}C'}:(R,\fm_R)\rightarrow (R\widehat{\otimes}_{C}C',\mathfrak{M}^e)$$ is a weakly unramified flat local extension of complete local rings.
		\item[(iii)] 
		   $K_{R\widehat{\otimes}_{C}C'}$ is  isomorphic to $K_{C'}$ as $K_R$-algebras where the former is considered as an $K_R$-algebra via $\overline{\eta_{R,R\widehat{\otimes}_CC'}}$ and the latter is considered as an $K_R$-algebra via $\overline{\varphi}\circ \overline{\lambda_{R}}^{\ -1}$. In particular,  $\overline{\varphi}$ is finite algebraic (resp. algebraic) field extension if and only if $\overline{\eta_{R,R\widehat{\otimes}_CC'}}$ is also a finite algebraic (resp. algebraic) field extension.
		 \item[(iv)]  The natural map $\lambda_{R\widehat{\otimes}_CC'}(:=\eta_{C',R\widehat{\otimes}_CC'}):C'\rightarrow R\widehat{\otimes}_{C}C'$ is a coefficient ring for $R\widehat{\otimes}_CC'$.
	\end{enumerate}
	\begin{proof}
		(i) and (iii) Note that, if $C$ (equivalently $C'$) is a field  then $K_C=C$, $K_{C'}=C'$ and $\fm_C=0=\fm_{C'}$. Since  $\lambda_R:C\rightarrow R$    is a coefficient ring so the  induced homomorphism $\overline{\lambda_R}:K_{C}\rightarrow K_R$ on  residue fields is an isomorphism. 	Therefore, 
		\begin{align*}
		\numberthis
		\label{Isomoorphism}
		K_R\otimes_{C}K_{C'}
		\cong \Big(K_R\otimes_{C/\mathfrak{m}_{C}}\big(C/\mathfrak{m}_C\big)\Big)\otimes_{C}K_{C'}
		\cong K_R\otimes_{C/\mathfrak{m}_{C}}\Big(\big(C/\mathfrak{m}_{C}\big)\otimes_{C}K_{C'}\Big)
		&  \cong K_R\otimes_{K_C} K_{C'}
		&\\& \overset{\text{}}{\cong} K_{C'}
		\end{align*}
		where the isomorphism is an isomorphism of $K_R$-algebras (with the $K_R$-algebra structure of $K_{C'}$ mentioned in the statement of part (iii) and the natural $K_R$-algebra structure of $K_R\otimes_CK_{C'}$ given by $\kappa\in K_R\mapsto \kappa\otimes 1$). Thus $\mathfrak{M}=\mathfrak{m}_R(R\otimes_{C}C')+\mathfrak{m}_{C'}(R\otimes_{C}C')\in \Max(R\otimes_{C}C')$  by Fact \ref{ElementaryTensorProductFact}(ii).   
		
		Moreover, $K_R\otimes_CK_{C'}\cong K_{R\widehat{\otimes}_CC'}$ as $K_R$-algebras by Proposition \ref{IdealIsMaximal}(ii). From this fact as well as (\ref{Isomoorphism}), we conclude the first statement of part (iii). The second statement of part (iii)  is  evident, as 
		$\overline{\lambda_R}^{\ -1}$ is an isomorphism.

 If $C'$ (equivalently $C$) is a field then $\fm_{C'}=0$ hence  $\mathfrak{M}=\fm_R(R\otimes_{C}C')$. If otherwise $C'$ is an unramified discrete valuation ring of mixed characteristic $(0,p)$ (equivalently $C$  is a $p$-ring)  then $R$ has mixed characteristic  and $\fm_{C'}=pC'$, so $\fm_{C'}(R\otimes_CC')=p(R\otimes_CC')\subseteq \fm_R(R\otimes_CC')$ implying that again $\mathfrak{M}=\fm_R(R\otimes_{C}C')$. This completes the proof of part (i).
		
	(ii)	The fact that $\mathfrak{M}\in \text{Max}(R\otimes_{C}C')$ together with  $\mathfrak{M}=\fm_R(R\otimes_CC')$  shows that $\eta_{R,R\widehat{\otimes}_{C}C'}$ is a  weakly unramified local homomorphism of complete local rings by Proposition \ref{IdealIsMaximal}(i). Moreover, from the maximality of $\mathfrak{M}$, Remark \ref{RemarkExtensionOfUnramifiedDVRs} and Proposition \ref{FlatProposition}(ii) we conclude that $\eta_{R,R\widehat{\otimes}_CC'}$ is a flat  homomorphism (in case $C$ and $C'$ are fields, any field extension is evidently flat).
	
	(iv)  Note that, by part (ii) $(R\widehat{\otimes}_CC',\mathfrak{M}^e)$ is a complete local ring (so  a coefficient ring for $R\widehat{\otimes}_CC'$ makes sense). We have the commutative diagram 
	\begin{center}
			$\begin{CD}
		C'
		@>\lambda_{R\widehat{\otimes}_CC'}=\eta_{C',R\widehat{\otimes}_CC'},\ \ c\mapsto \Big((1+\mam_R^n)\otimes (c+\mam_{C'}^n)\Big)_{n\in \mathbb{N}}>>
		R\widehat{\otimes}_{C}C' \\
		@V\text{surjective}VV @V\text{natural epimorphism}VV\\
		K_R\otimes_{C}K_{C'}\cong K_R\widehat{\otimes}_{C}K_{C'}  
		@< \cong  < \text{Proposition \ref{IdealIsMaximal}(iii)} < 
		R\widehat{\otimes}_CC'/\mathfrak{M}^e
		\end{CD}$
	\end{center}
 whose left vertical map, given by $C'\overset{\text{natural epimorphism}}{\twoheadrightarrow}K_{C'}\overset{\kappa\rightarrow 1\otimes\kappa}{\rightarrow}K_R\otimes_{C}K_{C'}$, is surjective in view of  (\ref{Isomoorphism}). Note that the commutativity of the diagram follows from the fact that the map $R\widehat{\otimes}_CC'/\mathfrak{M}^e\rightarrow K_R\widehat{\otimes}_{C}K_{C'}$ is an $(C'/\fm_{C'})$-algebra isomorphism, in particular so an $C'$-algebra isomorphism (see Proposition \ref{IdealIsMaximal}(iii)). Then an easy diagram chasing shows that $C'\overset{\eta_{C',R\widehat{\otimes}_{C}C'}}{\rightarrow} R\widehat{\otimes}_CC'\twoheadrightarrow R\widehat{\otimes}_CC'/\mathfrak{M}^e$ is also surjective. Consequently the induced map $\overline{\eta_{C',R\widehat{\otimes}_CC'}}:K_{C'}\rightarrow R\widehat{\otimes}_CC'/\mathfrak{M}^e$ is an isomorphism.
	\end{proof}
\end{subfact}

\begin{subfact}\label{RemarkMixedCharacteristicPRingExtension}
	Suppose that  $\lambda_R:C\rightarrow R$ is a coefficient ring of a complete local ring and $h:K_R\rightarrow \mathcal{K}$ is a field extension. Assume that $R$ has mixed characteristic  so $C$ is a $p$-ring ($p=\text{Char}(K_R)$).	
	Then there is an extension of  $p$-rings $\varphi:C\rightarrow C'$  which its induced residue field extension $\overline{\varphi}$   fits into the following commutative diagram of field extension (in particular $\overline{\varphi}\circ \overline{\lambda_R}^{\ -1}:K_R\rightarrow K_{C'}$ is $K_R$-algebra isomorphic to $h:K_R\rightarrow \mathcal{K}$): 
		\[
		\xymatrix{
			K_C \ar[d]_{\overline{\lambda_R}}^{\cong} \ar[rrrr]^{\overline{\varphi}} & &&& K_{C'} \ar[d]_{\cong}^{\phi}\\
			K_R \ar[rrrru]^{\overline{\varphi}\circ \overline{\lambda_R}^{\ -1}} \ar[rrrr]_{h} & & &&\mathcal{K}
		}
		\]
  \begin{proof}
  	  Considering the $p$-ring $C$ and the field extension $h\circ \overline{\lambda_R}:K_C\rightarrow \mathcal{K}$  we appeal to \cite[Theorem 29.1]{Matsumura} to find an unramified discrete valuation ring $V$ of mixed characteristic $(0,p)$, containing $C$,  such that $K_V=\mathcal{K}$ and an extension  $\gamma:C\overset{\text{}}{\rightarrow}V$ lifting the field extension $h\circ \overline{\lambda_R}:K_C\rightarrow \mathcal{K}$. 
  	  Let $C'$ be the $pV$-adic completion of $V$, and let $\varphi:C\rightarrow C'$ be the composition of  $\gamma$ with the canonical map to the completion $\mu:V\rightarrow C'=\widehat{V}$. Let $\phi:=\overline{\mu}^{\ -1}:K_{C'}\rightarrow \mathcal{K}=K_{V}$ ($\overline{\mu}$ is an isomorphism, see \cite[(4), page 63]{Matsumura}). Then, $\phi\circ \overline{\varphi}=\overline{\mu}^{\ -1}\circ (\overline{\mu}\circ \overline{\gamma})=\overline{\gamma}=h\circ \overline{\lambda_R}$, as required.
  \end{proof}		
\end{subfact}

\begin{subfact}\label{PolynomialFact}
	Suppose that $R$ is an integral domain which is an uncountable set. Let $n\in \mathbb{N}$. Then  for any at most countable set $F:=\{0\neq f_i\in R[X_1,\ldots,X_n]\}_{i\in \mathbb{N}}$ consisting of non-zero polynomials, there exists some $(r_1,\ldots,r_n)$ with $f_i(r_1,\ldots,r_n)\neq 0$ for all $i\in \mathbb{N}$.
	\begin{proof}
		If $n=1$ then the statement is easily verified (because any non-zero polynomial of degree $m$ with coefficients in an integral domain admits at most $m$ roots, thus the union of roots of $f_i$ ($i\in I$) is an at most countable set).  So suppose that $n>1$ and the statement has been proved for smaller values than $n$. Let   $F'$ be the subset of $F$ consisting of those polynomials $f_i$  having at least one monomial wherein   $X_n$ appears. Then, $F'$, can be considered as an at most countable set of non-zero polynomials in one indeterminate $X_n$ and with coefficients in the uncountable domain $R[X_1,\ldots,X_{n-1}]$, thus by the base of our induction there exists some $g(X_1,\ldots,X_{n-1})\in R[X_1,\ldots,X_{n-1}]$ such that $f_i(X_1,\ldots,X_{n-1},g)\neq 0$ for each $f_i\in F'$.  If $f_i\notin F'$, then we use  conventions $f_i(X_1,\ldots,X_{n-1},g):=f_i$ and $f_i(X_1,\ldots,X_{n-1},g)(r_1,\ldots,r_{n-1}):=f_i(r_1,\ldots,r_{n-1})$. Then, $\{f_i(X_1,\ldots,X_{n-1},g)\}_{i\in \mathbb{N}}$ is an at most countable set of non-zero polynomials with coefficients in the uncountable ring $R$, and supported at  (at most) $n-1$ number of variables $X_1,\ldots,X_{n-1}$. Hence by our inductive hypothesis there exist $r_1,\ldots,r_{n-1}$ with $$f_i(X_1,\ldots,X_{n-1},X_n)\big(r_1,\ldots,r_{n-1},g(r_1,\ldots,r_{n-1})\big)=f_i(X_1,\ldots,X_{n-1},g)(r_1,\ldots,r_{n-1})\neq 0,$$ for each $i\in \mathbb{N}$. Consequently, then $\big(r_1,\ldots,r_{n-1},g(r_1,\ldots,r_{n-1})\big)$ would be our desired point.
	\end{proof}
\end{subfact}

\begin{subfact}\label{PolynomialOverFiniteField}
	Suppose that $K$ is a finite field and $f\in K[X_1,\ldots,X_n]$ is a non-zero polynomial. Then there is a finite extension $L$ of $K$ and some  $(\lambda_1,\ldots,\lambda_n)\in L^n$ such that $f(\lambda_1,\ldots,\lambda_n)\neq 0$. 
	\begin{proof}
		Let $\overline{K}$ be the algebraic closure of $K$. Since	$\overline{K}$, being algebraically closed, is an  infinite field and $f$ is non-zero so it is well-known that there is some $(\lambda_1,\ldots,\lambda_n)\in \overline{K}$ with $f(\lambda_1,\ldots,\lambda_n)\neq 0$.
		Then, set $L:=K[\lambda_1,\ldots,\lambda_n]$ which is a field and is a finite extension of $K$ because the elements $\lambda_1,\ldots,\lambda_n$ are algebraic over $K$.
	\end{proof}
\end{subfact}

\begin{subfact}\label{FactMaximalSubfieldAlgebraicExtension}
	Suppose that $A$ is a complete local ring and $K$ is a subfield of $A$  that is maximal in $A$ with respect to  inclusion. Then $K\overset{\text{inclustion}}{\hookrightarrow} A\overset{\text{nat. epi.}}{\twoheadrightarrow} K_A$ is an algebraic field extension.
	\begin{proof}
		Consider some $a\in A$   whose residue class $a+\fm_A$ does not belong to the image of $K$ in $K_A$. We claim that  there is some $0\neq \sum\limits_{i=0}^nk_ia^i\in K[a]$ which is not invertible in $A$. Namely,   if all non-zero elements $\sum\limits_{i=0}^n k_ia^i$ ($n\in \mathbb{N}_0$) are invertible in $A$ then we first conclude that $K[a]$ is an integral domain, and then we see that  the inclusion $K[a]\subseteq A$ induces an injective ring homomorphism of $K$-algebra $$K(a)=\text{Frac}(K[a])\rightarrow A,\ \ \ \ \  (\sum\limits_{i=0}^nk_ia^i)/(\sum\limits_{i=0}^mk'_ia^i)\mapsto (\sum\limits_{i=0}^nk_ia^i)(\sum\limits_{i=0}^mk'_ia^i)^{-1}.$$
		However, then we conclude that there is a subfield of $A$ containing $K$ and the element  $a$, while this contradicts  the maximality of the subfield $K$ of $A$.
		
		So our claim holds, and there is some    $0\neq \sum\limits_{i=0}^nk_ia^i$ which is not invertible in $A$, equivalently $\sum\limits_{i=0}^nk_ia^i\in \fm_A$. But then $\sum\limits_{i=0}^n(k_i+\fm_A)(a+\fm_A)^i=0$ in $K_A$ implying that $a+\fm_A$ is algebraic over the image of $K$ in $K_A$. Since $a+\fm_A$ was chosen arbitrary, so the statement follows.
	\end{proof}
\end{subfact}

\begin{subfact}\label{FactFlatAnnihilator}
	Suppose that $f:R\rightarrow S$ is a flat homomorphism of commutative rings and $M$ is a finitely generated $R$-module.
	\begin{enumerate}
		\item [(i)]
	 $0:_S(M\otimes_R S)=(0:_RM)S$. 
	 \item[(ii)] If $f$ is faithfully flat then $\big(0:_S(M\otimes_R S)\big)\cap R=0:_RM$.
	 \end{enumerate}
	\begin{proof}
		(i) The statement is clear provided $M$ is a cyclic module. Suppose that $M=\sum\limits_{i=1}^n Rm_i$ for some $m_1,\ldots,m_n\in M$. Then $0:_RM=\cap_{i=1}^n0:_R(Rm_i).$ Consequently,  
		\begin{align*}
		(0:_RM)S
		&=
		\big(\cap_{i=1}^n0:_R(Rm_i)\big)S\overset{\text{\cite[Theorem 7.4]{Matsumura}}}{=}\cap_{i=1}^n\Big(\big(0:_R(Rm_i)\big)S\Big)\overset{Rm_i\text{\ is cyclic}}{=}\cap_{i=1}^n\Big(0:_S\big((Rm_i)\otimes_RS\big)\Big)
		&\\& 
		=  \cap_{i=1}^n \Big(0:_S\big((m_i\otimes 1)S\big)\Big)=0:_S(M\otimes_R S).
		\end{align*}  
		The first equality in the second row holds because  the image of the embedding $\eta:Rm_i\otimes_R S\rightarrow M\otimes_RS$ is $\text{im}(\eta)=(m_i\otimes_R 1)S$, so $Rm_i\otimes_R S\cong (m_i\otimes1)S$ as $S$-modules.  The second inequality in the second row holds because $M\otimes_RS$ is generated by $m_1\otimes 1,\ldots,m_n\otimes 1$.
		
		(ii)  The statement follows from \cite[Theorem 7.5(ii)]{Matsumura} and the previous part.
	\end{proof}
\end{subfact}


\subsection{Test modules and transcendental extensions}


	In this subsection we show that the test module property ascends under flat local homomorphism $R\rightarrow R[X_1,\ldots,X_n]_{(\fm_R)}$ provided $K_R$ is uncountable or the test module has finite complete intersection dimension (see Corollary \ref{CorollaryTranscendentalExtension}). This is the first  step towards answering Question \ref{CDTQuestion} and Question \ref{CWQuestion}.
	
\begin{subrem} \label{HomIsBounded}	
	\emph{Consider a bounded complex  $$C_\bullet:=0\rightarrow C_l\rightarrow C_{l-1}\rightarrow \cdots\rightarrow C_{m+1}\rightarrow C_m\rightarrow 0$$   over some  ring $R$.}
	\begin{enumerate}
		\item [\emph{(i)}]  \emph{Suppose that  $X_\bullet:=\cdots\rightarrow X_{t+2}\rightarrow X_{t+1}\rightarrow  X_{t}\rightarrow 0$  is a bounded below complex of $R$-modules. 
			  Then an easy calculation shows that $$\text{Hom}_R(X_\bullet,C_\bullet)_{i}=\prod\limits_{k\in \mathbb{Z}}\text{Hom}_R(X_k,C_{k+i})=0,\ \ \forall\ i\gneq l-t.$$}
		\item[\emph{(ii)}] \emph{Now suppose that the complex $X_\bullet:=0\rightarrow X_s\rightarrow X_{s-1}\rightarrow \cdots \rightarrow X_{t+1}\rightarrow X_t\rightarrow 0$  in part (i) is bounded. Then an easy verification shows that  $$\text{Hom}_R(X_\bullet,C_\bullet)_{i}=\prod\limits_{k\in \mathbb{Z}}\text{Hom}_R(X_k,C_{k+i})=0,\ \ \forall\ i\lneq  m-s.$$  
			This together with part (i) implies that then  $\text{Hom}_R(X_\bullet,C_\bullet)$ is a bounded complex.}
	\end{enumerate} 
\end{subrem}

	The next two lemmas will be used in the proof of Corollary \ref{CorollaryUniformTorAnnihilator}.

\begin{sublem}\label{LemmaRingOfOperators}
	Suppose that $S$ is a Noetherian ring and  $Q\twoheadrightarrow S$ is a deformation. Let $M\in \text{mod}(S)$ such that $\pd_Q(M)<\infty$. Assume that $$C_\bullet:=0\rightarrow C_0\rightarrow C_{-1}\rightarrow C_{-2}\rightarrow \cdots\rightarrow C_{-u}\rightarrow 0$$ is a bounded homologically finite complex of $S$-modules and  $w\in \mathbb{N}_0$. Then  there is some $v\in \mathbb{N}_0$ such that $$\big(\bigcap\limits_{w\le i\le w+v+1}0:_S\text{\emph{Ext}}^i_S(M,C_\bullet)\big)\subseteq 0:_S\text{\emph{Ext}}^j_S(M,C_\bullet), \ \ \forall\  j\ge w.$$
	\begin{proof}
		Let $F_\bullet$ be the  projective resolution of $M$ over $Q$ consisting of finitely generated free $Q$-modules which is a bounded complex ($\pd_Q(M)<\infty$). Since $F_\bullet$ and $C_\bullet$ are both bounded, so  $\text{Hom}_Q(F_\bullet,C_\bullet)$ is also a bounded complex implying that $\text{Ext}^{i}_Q(M,C_\bullet)=0$ for all but finitely many $i\in \mathbb{Z}$ (see Remark \ref{HomIsBounded}(ii)). Also,  $\text{Hom}_Q(F_\bullet,C_\bullet)=\text{\textbf{R}Hom}_Q(M,C_\bullet)$ is  homologically finite in view of \cite[(A.4.4) Lemma]{ChrisensenGorenstein}. 
		It follows that, $\text{Ext}^{*}_Q(M,C_\bullet)=\bigoplus\limits_{n\in \mathbb{N}_0}\text{Ext}^{n}_Q(M,C_\bullet)$ is a finite direct sum of Noetherian $S$-modules ($\text{Ext}^i_Q(M,C_\bullet)=H_{-i}\big(\text{Hom}_Q(F_\bullet,C_\bullet)\big)=0$ for all $i\lneq 0$ by Remark \ref{HomIsBounded}(i)). Hence $\text{Ext}^*_Q(M,C_\bullet)$ is a Noetherian $S$-module (equivalently $Q$-module). Consequently, in the light of and beauty of \cite[Theorem, page 708]{AvramovSunCohomology}, over the graded ring of operators $$\mathscr{S}:=S[\mathcal{X}_1,\ldots,\mathcal{X}_c;\ \text{deg}(\mathcal{X}_i)=2] $$ $\bigoplus\limits_{i\in \mathbb{N}_0}\text{Ext}^i_S(M,C_\bullet)$ (which admits a  graded $\mathscr{S}$-module  structure)  is a  finitely generated $\mathscr{S}$-module (the subscript $c$ appeared in the ring of operators is the height of the defining ideal of $S$ as a quotient of $Q$).
		Hence there is some $v\in \mathbb{N}_0$ such that $$\text{Ext}^0_{S}(M,C_\bullet),\ldots,\text{Ext}^v_{S}(M,C_\bullet)$$ generates $\text{Ext}^*_{S}(M,C_\bullet)$ over $\mathscr{S}$ (by Remark \ref{HomIsBounded}(i)   all  negative Ext modules $\text{Ext}^i_S(M,C_\bullet)$ vanish). 
		
		Set $\maa:=\bigcap\limits_{w\le i\le w+v+1}0:_S\text{Ext}^i_S(M,C_\bullet)$. It suffices to show that $\maa \big(\text{Ext}^i_S(M,C_\bullet)\big)=0$ for each $i\ge v+w+2$.  Pick some $e_i\in \text{Ext}^i_{S}(M,C_\bullet)$ for some  $i\ge v+w+2$. Then,  $e_i=\sum\limits_{j}h_j(\underline{\mathcal{X}})e_{j}$ for some homogeneous elements $h_j(\underline{\mathcal{X}})\in \mathscr{S}$ and $e_j\in\text{Ext}^{k_j}_{S}(M,C_\bullet)$ such that   $0\le k_j\le v$ and $\deg\big(h_j(\underline{\mathcal{X}})\big)\ge w+1$. It is easily seen that $$h_j(\underline{\mathcal{X}})=\sum_{l_{j,1}+\cdots+l_{j,c}=\lceil w/2\rceil}\limits\mathcal{X}_1^{l_{j,1}}\cdots\mathcal{X}_c^{l_{j,c}}h'_{j,\underline{l}}(\underline{\mathcal{X}})$$ for some elements $h'_{j,\underline{l}}(\underline{\mathcal{X}})\in \mathscr{S}$. Then, since each $\mathcal{X}_1^{l_{j,1}}\cdots\mathcal{X}_c^{l_{j,c}}e_j\in \dsum\limits_{h=w}^{w+v+1}\text{Ext}^h_{S}(M,C_\bullet)$ so it is killed by $\maa$. It follows that $e_i$ is killed by $\maa$.
	\end{proof}
\end{sublem}

\begin{sublem}\label{SpectralSequenceFact}	Suppose that $R$ is a  ring and
	$C_\bullet$ is an arbitrary complex. Consider another complex  $$D_\bullet:=0\rightarrow D_0\rightarrow D_{-1} \rightarrow \cdots \rightarrow D_{-t+1}\rightarrow D_{-t}\rightarrow 0$$    consisting of injective $R$-modules. Set  $\mathfrak{a}:=\bigcap\limits_{i=-n-t}^{-n}0:_RH_i(C_\bullet)$ for some $n\in\mathbb{Z}$. Then we have $$\mathfrak{a}^{t+1}\subseteq 0:_RH_{n}\big(\Hom_R(C_\bullet,D_\bullet)\big)$$ (note that $t$ is the length of the complex $D_\bullet$).
	\begin{proof} 
		Since $D_\bullet$ is a bounded complex, so we could give a  proof in one step by using one spectral sequence concerning the unbounded complex $C_\bullet$ and the bounded complex $D_\bullet$. However to make the proof and our spectral sequence compatible with the third quadrant spectral sequences in \cite{Rotman}, we break the proof in three steps.
		
		\textbf{Step 1:}  In this step, we consider the case where  $C_\bullet:=\cdots \rightarrow C_m\rightarrow \cdots C_1\rightarrow C_0\rightarrow 0$ is bounded below and concentrated in non-negative degrees. 
		
		We bear in  mind the spectral sequence arising from the  third quadrant bicomplex $$M_{p,q}:=\text{Hom}_R(C_{-p},D_q).$$  
		We recall that $\text{Tot}(M)=\text{Hom}_R(C_\bullet,D_\bullet)$, and the differentials of the bicomplex are pointing left and down. We should consider the second filtration that is  obtained by   restricting  the considered rows. 
		Namely, setting $^{\text{II}}F^{p}\big(\text{Tot}(M)\big)$ to be the subcomplex of $\text{Tot}(M)$ with $^{\text{II}}F^{p}\big(\text{Tot}(M)\big)_{n}:=\oplus_{j\le p} M_{n-j,j}$, we have the filtration 
		$$ \cdots\subseteq\  ^{\text{II}}F^{-2}\big(\text{Tot}(M)\big)\subseteq \ ^{\text{II}}F^{-1}\big(\text{Tot}(M)\big)\subseteq \ ^{\text{II}}F^{0}\big(\text{Tot}(M)\big)=\text{Tot}(M).$$
		
		Hence the $^{\text{II}}E^1$ page  of the spectral sequence is given by taking the homologies of the rows of the  bicomplex. Thus 
		\begin{equation}
		\label{EquationE1Page}
		^{\text{II}}E^1_{p,q}=H_{p}\big(\text{Hom}_R(C_\bullet,D_{q})\big)\overset{D_q\text{\ is  injective}}{\cong}\text{Hom}_R\big(H_{-p}(C_\bullet),D_q\big).  
		\end{equation}		
		$H_{n}\big(\text{Tot}(M)\big)$  has the filtration, 
		\begin{align*}
		\label{EquationHomologyFiltration}
		\numberthis
		0= 
		\phi^{n-1}\Big(H_{n}\big(\Tot(M)\big)\Big)
		\subseteq \phi^{n}\Big(H_{n}\big(\Tot(M)\big)\Big)
		\subseteq \cdots 
		\subseteq\phi^{-1}\Big(H_{n}\big(\Tot(M)\big)\Big)
		&\subseteq \phi^0\Big(H_{n}\big(\Tot(M)\big)\Big)
		&\\&
		=H_{n}\big(\Tot(M)\big)
		\end{align*}
		such that  
		\begin{equation}
		\label{EquationEInfty}
		\phi^q\Big(H_{n}\big(\Tot(M)\big)\Big)/\phi^{q-1}\Big(H_{n}\big(\Tot(M)\big)\Big)\cong {^{\text{II}}E^\infty_{n-q,q}}
		\end{equation}
		is a subquotient of ${^{\text{II}}E^1_{n-q,q}}\overset{\text{(\ref{EquationE1Page})}}{\cong} \text{Hom}_R\big(H_{q-n}(C_\bullet),D_q\big).$		   
		
		If $q\in (-\infty,-t-1)\cup \mathbb{N}$ then, $D_q=0$, so $\text{Hom}_R\big(H_{q-n}(C_\bullet),D_q\big)$ and ${^{\text{II}}E^\infty_{n-q,q}}$ vanish. If, otherwise $q\in [-t,0]$ then $-t-n\le q-n\le -n$, so by our hypothesis $\mathfrak{a}\subseteq 0:_R\text{Hom}_R\big(H_{q-n}(C_\bullet),D_q\big)\subseteq 0:_R{^{\text{II}}E^\infty_{n-q,q}}$. From this, (\ref{EquationHomologyFiltration}) and (\ref{EquationEInfty}) it is easily seen that $\mathfrak{a}^{t+1}\Big(H_{n}\big(\text{Tot}(M)\big)\Big)=0$.
		
		\textbf{Step 2:} Let $C_\bullet:=\cdots\rightarrow C_m\rightarrow \cdots \rightarrow C_{u+1}\rightarrow C_u\rightarrow 0$ be an arbitrary bounded below complex. Note that $\sum\limits^{-u}C_\bullet$ satisfies the assumption of the previous step. By our hypothesis, $\mathfrak{a}\big(H_{-u-n-i}(\sum\limits^{-u}C_\bullet)\big)=\mathfrak{a}\big(H_{-n-i}(C_\bullet)\big)=0$ for each $0\le i\le t$. So,  $\mathfrak{a}^{t+1}\subseteq 0:_R\Big(H_{u+n}\big(\text{Hom}_R(\sum\limits^{-u}C_\bullet, D_\bullet)\big)\Big)$ by step 1. But, $H_{u+n}\big(\text{Hom}_R(\sum\limits^{-u}C_\bullet, D_\bullet)\big)=H_{u+n}\big(\sum\limits^u\text{Hom}_R(C_\bullet,D_\bullet)\big)=H_{n}\big(\text{Hom}_R(C_\bullet,D_\bullet)\big)$. Hence, $$\mathfrak{a}^{t+1}\subseteq 0:_RH_{n}\big(\text{Hom}_R(C_\bullet,D_\bullet)\big)$$ as claimed.
		
		\textbf{Step 3:} Let $(C_\bullet,\partial_\bullet^C)$ be an arbitrary complex. Set, $$C'_\bullet:=0\longrightarrow C_{-n-t-2}\overset{\partial^C_{-n-t-2}}{\longrightarrow} C_{-n-t-1}\overset{\partial^C_{-n-t-1}}{\longrightarrow} \cdots\rightarrow C_{-n+1}\overset{\partial^C_{-n+1}}{\longrightarrow}C_{-n+2}\longrightarrow 0.$$
		
		Note that $\text{Hom}_R(C_\bullet,D_\bullet)_{n\pm 1}=\prod\limits_{i\in \mathbb{Z}}\text{Hom}_R(C_{i-n\mp 1},D_i)=\bigoplus\limits_{i=-t}^0\text{Hom}_R(C_{i-n\mp 1},D_i)=\text{Hom}_R(C'_\bullet,D_\bullet)_{n\pm 1}$ implying that $H_{n}\big(\text{Hom}_R(C_\bullet,D_\bullet)\big)=H_{n}\big(\text{Hom}_R(C'_\bullet,D_\bullet)\big)$. By our hypothesis  $$\mathfrak{a}\big(H_{-n-i}(C'_\bullet)\big)=\mathfrak{a}\big(H_{-n-i}(C_\bullet)\big)=0,\ \ \forall -t\le i\le 0.$$ So, Step 2 yields $\mathfrak{a}^{t+1}\subseteq 0:_RH_{n}\big(\text{Hom}_R(C'_\bullet,D_\bullet)\big)=0:_RH_{n}\big(\text{Hom}_R(C_\bullet,D_\bullet)\big)$, as was to be proved.   	\end{proof}	 
\end{sublem}

The next corollary will be used in the proof of Lemma \ref{LemmaDesiredPerfectComplexExists}.

\begin{subcor}\label{CorollaryUniformTorAnnihilator}
	Suppose that $S$ is a Noetherian ring and $Q\twoheadrightarrow S$ is a deformation. Let $M$ be a finitely generated $S$-module with $\pd_Q(M)<\infty$. Assume that $S$ admits a dualizing complex $$D_\bullet:=0\rightarrow D_0\rightarrow D_{-1}\rightarrow \cdots \rightarrow D_{-t}\rightarrow 0.$$ Then for any  $N\in \text{mod}(S)$ and  $w\ge t$ ($t$ is the subscript appeared in the dualizing complex), there exists $v\in \mathbb{N}_0$ such that 
	$$ \Big(\bigcap\limits_{i=w-t}^{v+w+1}0:_S\emph{Tor}^S_i(M,N)\Big)^{2t+2} \subseteq 0:_S\emph{Tor}^S_j(M,N),\ \ \ \forall\  j\ge w.$$
	\begin{proof}
		Fix some $w\ge t$. Let $F^M_\bullet$ be a free resolution for $M$ over $S$ consisting of finitely generated free $S$-modules. Set $X_\bullet:=F^M_\bullet\otimes_SN$.  Since $$X_\bullet^{\vee}:=\text{Hom}_S(X_\bullet,D_\bullet)
		=\text{Hom}_S(F^M_\bullet\otimes_SN,D_\bullet)
		\overset{\text{\cite[(A.2.8) Adjointness]{ChrisensenGorenstein}}}{\cong}
		\text{Hom}_S\big(F^M_\bullet,\text{Hom}_S(N,D_\bullet)\big),$$ so 
		\begin{equation}
		\label{EquationHomologyToExt}
		H_{-i}(X_\bullet^\vee)\cong \text{Ext}^i_S\big(M,\text{Hom}_S(N,D_\bullet)\big),\ \ \forall \ i\in \mathbb{Z}.
		\end{equation}
		Note that  by \cite[A.4.4]{ChrisensenGorenstein}, the bounded complex $\text{Hom}_S(N,D_\bullet)\simeq \text{\textbf{R}Hom}_S(N,D_\bullet)$ is homologically finite. So, as $\pd_Q(M)<\infty$,  from Lemma \ref{LemmaRingOfOperators}     we can conclude that there exists $v\in \mathbb{N}_0$ with 
		$$\Big(\bigcap\limits_{w\le i\le w+v+1}0:_S\text{Ext}^i_S\big(M,\text{Hom}_S(N,D_\bullet)\big)\Big) \subseteq 0:_S\text{Ext}_S^{j}\big(M,\text{Hom}_S(N,D_\bullet)\big),\ \ \forall\ j\ge w.$$ 
		From this and (\ref{EquationHomologyToExt})  we get    the second containment in the following display:
		$$
		\big(\bigcap\limits_{i=w-t}^{v+w+1}0:_SH_i(X_\bullet)\big)^{t+1}\overset{\text{Lemma \ref{SpectralSequenceFact}}}{\subseteq}\bigcap\limits_{-w-v-1\le i\le -w}0:_SH_{i}(X_\bullet^\vee)\subseteq 0:_SH_{j}(X_\bullet^\vee),\ \ \forall\ j\le  -w.
		$$ 
		
		Setting $\mathfrak{a}:=\bigcap\limits_{i=w-t}^{v+w+1}0:_SH_i(X_\bullet)$, by the above display  we have $\mathfrak{a}^{t+1}\subseteq \bigcap\limits_{j\le -w}0:_SH_{j}(X_\bullet^\vee)$. Thus, again from Lemma \ref{SpectralSequenceFact} we get 
		\begin{equation}
		\label{EquationExtAnnihilator}
		\mathfrak{a}^{2t+2}\subseteq 0:_SH_j(X_\bullet^{\vee\vee}),\ \ \forall j\ge w.
		\end{equation}
		Since $X_\bullet$ is a homologically finite bounded below complex of $S$-modules and $D_\bullet$ is a dualizing complex for $S$, so $X_\bullet\simeq X_\bullet^{\vee\vee}:=\text{Hom}_S\big(\text{Hom}_S(X_\bullet,D_\bullet),D_\bullet\big)$ 
		by \cite[(7.2.2) Theorem]{ChristensenFoxbyHyperHomological} (see \cite[Lists of Symbols]{ChristensenFoxbyHyperHomological} for the notation used in the statement of \cite[(7.2.2) Theorem]{ChristensenFoxbyHyperHomological}). From this and  (\ref{EquationExtAnnihilator}) we deduce that
		$$
		\maa^{2t+2} \subseteq 0:_SH_j(X_\bullet),\ \ \ \forall\  j\ge w.
		$$
		Hence, our proof is complete as $H_i(X_\bullet)=\text{Tor}^S_i(M,N)$ for all $i$.
	\end{proof}
\end{subcor}

\begin{subfact}\label{FactKernelBeingRegularSequence}
	Let $\mathbf{X}:=X_1,\ldots,X_n$ be a sequence of indeterminates over some ring $R$.  Consider the ring homomorphism $\varphi_n:R[\mathbf{X}]\rightarrow R$ given by the rule $X_i\mapsto c_i$ and $r\in R\mapsto r$. Then $\text{Ker}(\varphi_n)$ is generated by $X_1-c_1,\ldots,X_n-c_n$  which is  a regular sequence on $R[\mathbf{X}]$.
	\begin{proof}
		We use induction on $n$.  Suppose that $n=1$, so that only one indeterminate is involved. Since $X_1-c_1$ is a monic polynomial so it is a regular element (see \cite[Exercise 2(iii), page 11]{AtiyahIntroduction}). In order  to verify that $X_1-c_1$ generates $\Ker(\varphi_1)$, let $f=\sum\limits_{j=0}^{m}r_jX_1^j\in \Ker(\varphi_1)$. Then, 
		\begin{align*}
		\sum\limits_{j=0}^{m}r_jX_1^j=\sum\limits_{j=0}^{m}r_j\big((X_1-c_1)+c_1\big)^j
		&
		=\sum\limits_{j=0}^{m}r_j\big(\sum\limits_{k=0}^j\binom{j}{k}(X_1-c_1)^kc_1^{j-k}\big)
		&\\&=
			(X_1-c_1)g+\sum\limits_{j=0}^mr_jc_1^j\ \ (\text{for some\ }g\in R[X_1])
		&\\&=
		(X_1-c_1)g+\varphi_1(f)
		&\\&=
		(X_1-c_1)g.
		\end{align*}	   
		Consequently $f\in (X_1-c_1)R[X_1]$ as required.
		
		Suppose that $n>1$ and the statement has been proved when the number of indeterminates is less than $n$. By our inductive hypothesis $X_1-c_1,\ldots,X_n-c_n$ is a regular sequence of $R[X_1,\ldots,X_{n-1}]$ and it is a generating set of the kernel of the ring epimorphism $\varphi_{n-1}:=R[X_1,\ldots,X_{n-1}]\rightarrow R$ given by the rule $X_i\mapsto c_i$ ($1\le i\le n-1$). Let $\psi$ be the ring epimorphism $R[X_1,\ldots,X_n]\rightarrow R[X_1,\ldots,X_{n-1}]$ that fixes $R$ and  $X_i$ for $1\le i\le n-1$ but $\psi(X_n)= c_n$. So, $\varphi_{n-1}\circ \psi=\varphi_n$. Hence, if $f\in \Ker(\varphi_{n})$ then  $\psi(f)\in \Ker\ \varphi_{n-1}=(X_i-c_i)_{i=1}^{n-1}$ by our inductive hypothesis.  We consider the commutative diagram 
		
		\begin{center}
			\[	\xymatrix{
				R[X_1,\ldots,X_n]  \ar[rr]^{\psi} \ar[d]_{\zeta}^{\text{natural epimorphism}} &&R[X_1,\ldots,X_{n-1}]\ar[d]^{\pi}_{\text{natural epimorphis}}\\
				R[X_1,\ldots,X_n]/(X_i-c_i)_{i=1}^{n-1} \ar[r]_{\Theta\ \ \ \ \ } &\big(R[X_1,\ldots,X_{n-1}]/(X_i-c_i)_{i=1}^{n-1}\big)[X_n] \ar[r]_{\ \ \ \overline{\psi}} & R[X_1,\ldots,X_{n-1}]/(X_i-c_i)_{i=1}^{n-1}
			} \]
		\end{center}
		where $\Theta$ is the isomorphism  induced by $$R[X_1,\ldots,X_n]\rightarrow \big(R[X_1,\ldots,X_{n-1}]/(X_i-c_i)_{i=1}^{n-1}\big)[X_n]$$ defined by $$r\mapsto r+(X_i-c_i)_{i=1}^{n-1},\  X_j\mapsto X_j+(X_i-c_i)_{i=1}^{n-1}\ (r\in R,1\le j\le n-1),\ X_n\mapsto X_n.$$ Moreover, $\overline{\psi}$ is the ring homomorphism that is identity on $R[X_1,\ldots,X_{n-1}]/(X_i-c_i)_{i=1}^{n-1}$ while $\overline{\psi}(X_n)=c_n+(X_i-c_i)_{i=1}^{n-1}$. Then from the  commutativity of the diagram and the above arguments we conclude that $\Theta\circ \zeta(f)\in \Ker\ \overline{\psi}$. Now, the base of our inductions implies that $\Theta\circ \zeta(f)$ belongs to the ideal generated  by $X_n-\big(c_n+(X_i-c_i)_{i=1}^{n-1}\big)$. Consequently, by applying $\Theta^{-1}$ we get $\zeta(f)$ belongs to the ideal generated by $(X_n-c_n)+(X_i-c_i)_{i=1}^{n-1},$ i.e. $f\in (X_i-c_i)_{i=1}^{n}$ as was to be proved.
		
		Finally, $X_n-\big(c_n+(X_i-c_i)_{i=1}^{n-1}\big)$ is a regular element of  $(R[X_1,\ldots,X_{n-1}]/(X_i-c_i)_{i=1}^{n-1}\big)[X_n]$  by the base of our induction, hence applying the isomorphism $\Theta^{-1}$ mentioned in the above diagram we conclude that  $(X_i-c_i)_{i=1}^{n}$ is a regular sequence of $R[\mathbf{X}]$.
	\end{proof}
\end{subfact}

The following lemma will be used in the proof of Corollary \ref{CorollaryTranscendentalExtension} which is the main result of this subsection. The lemma shows that the existence of a perfect complex over $R[\mathbf{X}]$ satisfying parts (i)-(iv) of the lemma, with mild conditions, implies that the module $T$ appeared in part (iv) is not a test $R$-module.

\begin{sublem}\label{EhsanWhenWillThisEnd}  Suppose that  $R$ is a  local ring of depth $s$ and $T\in \text{mod}(R)$. Let $\mathbf{X}:=X_1,\ldots,X_n$ be a sequence of indeterminates over $R$. Assume that there is a perfect complex $(G_\bullet,\partial^{G_\bullet}_\bullet)\in \mathcal{D}_b(R[\mathbf{X}])$ with the following properties:
	\vspace{1mm}
	\begin{enumerate}
		\item[(i)]   $G_{i}= 0$  if $i<0$ and  $G_i\neq 0$   if $0\le i\le n+s+3$ where $n$ is the number of indeterminates in $\mathbf{X}$ ($s=\depth(R)$).
		\item[(ii)] All entries of the matrices of   differentials of $G_\bullet$ belong to $\mam_RR[\mathbf{X}]$.
		\item [(iii)] There is  $g\in R[\mathbf{X}]\backslash \mam_RR[\mathbf{X}]$  such that   $$g\in \bigcap\limits_{i=1}^{n+s+2} 0:_{R[\mathbf{X}]}H_i(G_{\bullet}).$$ 
		\item[(iv)] For some $u\ge 1$ and for each $i\ge u$ there is $\zeta_i\in R[\mathbf{X}]\backslash \mam_RR[\mathbf{X}]$ such that  $$\zeta_i\in 0:_{R[\mathbf{X}]}\Big(\text{\emph{Tor}}^{R[\mathbf{X}]}_i\big(H_0(G_\bullet),T\otimes_RR[\mathbf{X}]\big)\Big).$$
		
	\end{enumerate}
	Finally, assume that there is a sequence  $c_1,\ldots,c_n$ of elements of $R$ such that  $g(c_1,\ldots,c_n)\in R\backslash\mam_R$ and $\zeta_i(c_1,\ldots,c_n)\in R\backslash \mam_R$  for each $i\ge u$. Then $T$ is not a test module for $R$.
	\begin{proof}
		We suppose to the contrary that $T$ is a test $R$-module, and we get a contradiction. 
		
		We define the  ring epimorphism $$\varphi:R[\mathbf{X}]/(X_i-c_i)_{i=1}^n\rightarrow R$$  that is induced by the rule  $r\mapsto r$ and $X_i\mapsto c_i$ for each $r\in R$ and $1\le i\le n$ (the elements $c_1,\ldots,c_n$ are mentioned in the statement). By Fact \ref{FactKernelBeingRegularSequence} $\varphi$ is an isomorphism and $(X_i-c_i)_{i=1}^n$ is a regular sequence of $R[\mathbf{X}]$ which we denote this sequence by the notation $\mathbf{y}$ throughout the rest of the proof.
		
		Let $$\lambda:R\rightarrow R[\mathbf{X}]/(\mathbf{y}),\ \ r\mapsto r+(\mathbf{y}).$$ From $\varphi\circ \lambda=\id_R$, we deduce that $\lambda$ is an isomorphism. In particular, $T\otimes_R(R[\mathbf{X}]/\big(\mathbf{y})\big)$ is a test module for $R[\mathbf{X}]/(\mathbf{y})$.  Consequently since $(T\otimes_RR[\mathbf{X}])/\mathbf{y}(T\otimes_RR[\mathbf{X}])\cong T\otimes_R(R[\mathbf{X}]/\big(\mathbf{y})\big)$ by Fact \ref{ElementaryTensorProductFact}(i), so $$(T\otimes_RR[\mathbf{X}])/\mathbf{y}(T\otimes_RR[\mathbf{X}])$$ is a test module for $R[\mathbf{X}]/(\mathbf{y})$. Let $$G'_\bullet:=0\rightarrow G_{n+s+3}\overset{\partial^{G_\bullet}_{n+s+3}}{\rightarrow} \cdots\rightarrow G_{n+1}\overset{\partial^{G_\bullet}_{n+1}}{\rightarrow} G_{n}\rightarrow 0$$ be a hard  truncation of $G_\bullet$. Note that $H_0(G_\bullet')/\mathbf{y}H_0(G_\bullet')\cong G_n/\big(\text{im}(\partial^{G_\bullet}_{n+1})+(\mathbf{y})G_n\big)\neq 0$. Namely, if otherwise, by projecting $G_n=R[\mathbf{X}]^{m}$ (for some $m$) to any free summand $R[\mathbf{X}]$ of it we would get $1\in \mathfrak{m}_RR[\mathbf{X}]+(\mathbf{y})$ (here we are applying the condition (ii) in the statement). Then we would have $1+(\mathbf{y})\in \mathfrak{m}_R\big(R[\mathbf{X}]/(\mathbf{y})\big)$ which contradicts with the fact that $\lambda:R\rightarrow R[\mathbf{X}]/(\mathbf{y})$ defined above is an isomorphism.
		
		In order to get our desired contradiction we aim to  show that
		\begin{equation}
		\label{EquationGetContradiction}
		\tor^{R[\mathbf{X}]/(\mathbf{y})}_{\gg}\Big(H_0(G_\bullet')/\mathbf{y}H_0(G_\bullet'), T\otimes_RR[\mathbf{X}]/\mathbf{y}(T\otimes_RR[\mathbf{X}])\Big)=0.
		\end{equation}
		
		Let $k\ge u$ and $\zeta:=\prod\limits_{j=k}^{k+n}\zeta_j$ so that 
		\begin{equation}
		\label{EquationZetaInAnnihilatorOfTors}
		\zeta\in 0:_{R[\mathbf{X}]}\Big(\tor^{R[\mathbf{X}]}_{j}\big(H_0(G_\bullet),T\otimes_RR[\mathbf{X}]\big)\Big),\ \ \ \forall\  k\le j\le k+n.
		\end{equation}
		By our hypothesis, $\varphi(\zeta_j+(\mathbf{y}))=\zeta_j(c_1,\ldots,c_n)\in R\backslash \fm_R$ for each $u\le j$ and $\varphi(g+(\mathbf{y}))=g(c_1,\ldots,c_n)\in R\backslash \fm_R$, consequently $g\zeta+(\mathbf{y})$ is an invertible element of $R[\mathbf{X}]/(\mathbf{y})$ ($\varphi$ is an isomorphism). Therefore, 
		\begin{center}
			\begin{equation}
			\label{IsomorphismOfTruncation}
			H_0(G_\bullet')/\mathbf{y}H_0(G_\bullet')\cong \big(H_0(G_\bullet')/\mathbf{y}H_0(G_\bullet')\big)_{g\zeta+(\mathbf{y})}\cong H_0({G_\bullet'}_{g\zeta})/\mathbf{y}\big(H_0({G_\bullet'}_{g\zeta})\big).
			\end{equation}
		\end{center}  
		
		In view of the condition (iii) mentioned in  the statement, the complex 
		\begin{equation} 
		\label{EquationLongExactSequence}
		0\rightarrow  H_0({G_\bullet'}_{g\zeta})\overset{\overline{{\partial^{{G_\bullet}}_{n}}_{g\zeta}}}{\rightarrow}{G_{n-1}}_{g\zeta}\overset{{\partial^{G_\bullet}_{n-1}}_{g\zeta}}{\rightarrow} {G_{n-2}}_{g\zeta}\overset{{\partial^{G_\bullet}_{n-2}}_{g\zeta}}{\rightarrow} \cdots \overset{{\partial^{G_\bullet}_{1}}_{g\zeta}}{\rightarrow} {G_0}_{g\zeta}\rightarrow H_0({G_\bullet}_{g\zeta})\rightarrow 0
		\end{equation}
		is exact and from Fact \ref{SyzygyAndGrade}(i) we can deduce that $\mathbf{y}$ is a regular sequence on $H_0({G'_\bullet}_{g\zeta})$ (we remind that $H_0({G'_\bullet}_{g\zeta})/\mathbf{y}\big(H_0({G'_\bullet}_{g\zeta})\big)\neq 0$ in view of (\ref{IsomorphismOfTruncation}) and the paragraph before (\ref{EquationGetContradiction})). So        \begin{align*}
		\numberthis
		\label{EquationTorIsomorphismss}
		\tor^{R[\mathbf{x}]/(\mathbf{y})}_k\big(H_0({G_\bullet'}_{g\zeta})/\mathbf{y}H_0({G_\bullet'}_{g\zeta}),T\otimes_RR[\mathbf{X}]/\mathbf{y}(T\otimes_RR[\mathbf{X}])\big)
		&\cong \tor^{R[\mathbf{x}]}_k\big(H_0({G_\bullet'}_{g\zeta}),T\otimes_RR[\mathbf{X}]/\mathbf{y}(T\otimes_RR[\mathbf{X}])\big)
		&\\&
		\overset{\text{(\ref{EquationLongExactSequence})}}{\cong} \tor^{R[\mathbf{x}]}_{k+n}\big(H_0({G_\bullet}_{g\zeta}),T\otimes_RR[\mathbf{X}]/\mathbf{y}(T\otimes_RR[\mathbf{X}])\big)
		&\\& 
		\cong \tor^{R[\mathbf{x}]}_{k+n}\big(H_0({G_\bullet}),T\otimes_RR[\mathbf{X}]/\mathbf{y}(T\otimes_RR[\mathbf{X}])\big)_{g\zeta}.
		\end{align*}
        where the first isomorphism holds in view of Fact \ref{FactTorOverRTorOverRxR}.

		Setting $i:=k+n$ from (\ref{EquationZetaInAnnihilatorOfTors})  we get $\tor^{R[\mathbf{x}]}_{j}\big(H_0({G_\bullet}),T\otimes_RR[\mathbf{X}]\big)_{g\zeta}=0$ for $i-n\le j\le i$. On the other hand, from Fact \ref{FactGradeOverNiceFaithfullyFlatBaseChange}   we get $\text{grade}_{(\mathbf{y})}(T\otimes_RR[\mathbf{X}])=n$. Hence Fact \ref{FactTorOfMAndTorOfSpecializationM} implies that $\tor^{R[\mathbf{x}]}_{k+n}\big(H_0({G_\bullet}),T\otimes_RR[\mathbf{X}]/\mathbf{y}(T\otimes_RR[\mathbf{X}])\big)_{g\zeta}=0$, i.e. $$\tor^{R[\mathbf{x}]/(\mathbf{y})}_k\big(H_0({G_\bullet'})/\mathbf{y}H_0({G_\bullet'}),T\otimes_RR[\mathbf{X}]/\mathbf{y}(T\otimes_RR[\mathbf{X}])\big)=0$$ by (\ref{EquationTorIsomorphismss}) in conjunction with (\ref{IsomorphismOfTruncation}).
		
		As $k\ge u$ was arbitrary, so our desired vanishing in (\ref{EquationGetContradiction}) holds. But then  
		\begin{align*}
		\label{EquationProjectiveDimensionContradiction}
		\numberthis
		 \pd_{R[\mathbf{x}]/\mathbf{y}R[\mathbf{X}]}\big(H_0(G'_\bullet)/\mathbf{y}H_0(G'_\bullet)\big)\overset{(\ref{IsomorphismOfTruncation})}{\cong} \pd_{R[\mathbf{x}]/\mathbf{y}R[\mathbf{X}]}\big(H_0({G'_\bullet}_{g\zeta})/\mathbf{y}H_0({G'_\bullet}_{g\zeta})\big)\le s
		 \end{align*}
		 		 in view of the Auslander-Buchsbaum formula in conjunction with the fact that $T\otimes_RR[\mathbf{X}]/\mathbf{y}(T\otimes_RR[\mathbf{X}])$ is a test module for the depth $s$ local ring $R[\mathbf{X}]/(\mathbf{y})$. 
		 		 
		 		 But ${G'_\bullet}_{g\zeta}$ is a part of a flat resolution of $H_0({G'_\bullet}_{g\zeta})$ (by condition (iii)) and ${G'_\bullet}_{g\zeta}/\mathbf{y}({G'_\bullet}_{g\zeta})$ is a part of a free resolution of $H_0({G'_\bullet}_{g\zeta})/\mathbf{y}H_0({G'_\bullet}_{g\zeta})$ by Fact \ref{ResolutionModuleRegularSequence}. However, $${G'_\bullet}_{g\zeta}/\mathbf{y}({G'_\bullet}_{g\zeta})\cong (G'_\bullet/\mathbf{y}G'_\bullet)_{g\zeta+(\mathbf{y})}\cong G'_\bullet/\mathbf{y}G'_\bullet$$  is a part of the minimal free resolution because of the condition (ii) while it has length strictly more than $s$ in view of condition (i). This   contradicts with (\ref{EquationProjectiveDimensionContradiction}).
	\end{proof} 
\end{sublem}

The following lemma which will be used in the proof of Corollary \ref{CorollaryTranscendentalExtension}, shows that under certain vanishing of Tor modules over $R[\mathbf{X}]_{(\fm_R)}$,	a perfect complex $G_\bullet$ as in the statement of the previous lemma exists.

\begin{sublem}\label{LemmaDesiredPerfectComplexExists}
	Let $R$ be a local ring, $T\in \text{mod}(R)$  and $N\in \text{mod}(R[\mathbf{X}]_{(\mam_R)})$, where $\mathbf{X}:=X_1,\ldots,X_n$ is a sequence of indeterminates over $R$. Suppose that $$\text{\emph{Tor}}^{R[\mathbf{X}]_{(\mam_R)}}_{\gg}\big(N,T\otimes_R R[\mathbf{X}]_{(\mam_R)}\big)=0$$  and $\emph{pd}_{R[\mathbf{X}]_{(\mam_R)}}(N)=\infty$. Then  there is  a prefect complex $G_\bullet$ in $\mathcal{D}_b(R[\mathbf{X}])$ with the following properties:
	\begin{enumerate}
		\item [(i)]
		$G_\bullet$ satisfies the conditions $(i)-(iv)$ of the statement of Lemma \ref{EhsanWhenWillThisEnd}. 
		\item[(ii)] If, moreover, $T$ has finite complete intersection dimension then part (iv) of the statement of Lemma \ref{EhsanWhenWillThisEnd} can be improved and restated as follows: 
		
		There is some   $\varsigma\in R[\mathbf{X}]\backslash \mam_R[R\mathbf{X}]$ such that  $$\varsigma\in 0:_{R[\mathbf{X}]}\Big(\text{\emph{Tor}}^{R[\mathbf{X}]}_{i}\big(H_0(G_\bullet),T\otimes_RR[\mathbf{X}]\big)\Big), \ \ \ \forall\ i \gg 0.$$
	\end{enumerate} 
	\begin{proof}
		
		Suppose that $s:=\text{depth}(R)$ and let  $$F^N_{\bullet}:=0\rightarrow F^N_{n+s+3}\overset{}{\rightarrow} F^N_{n+s+2}\overset{}{\rightarrow}\cdots\overset{}{\rightarrow} F^N_{0}\rightarrow0,$$
		be the length $(3+n+s)$-hard truncation of the $R[\mathbf{X}]_{(\mam_R)}$-minimal
		resolution of $N$   ($n$ is the number of indeterminates in the sequence $\mathbf{X}$).
		Note that the entries of the matrices of the differentials of $F^N_\bullet$  belong to  $\mam_RR[\mathbf{X}]_{(\mam_R)}$,  and $F_i\neq 0$ for $1\le i\le n+s+3$ as $\pd_{R_{[\mathbf{X}]_{(\mam_R)}}}(N)=\infty$.
		Applying Fact \ref{EhsanStuck} to the complex $F^N_\bullet$, we obtain the complex $G_\bullet$ of finitely generated free $R[\mathbf{X}]$-modules with $G_\bullet\otimes_{R[\mathbf{X}]}(R[\mathbf{X}]_{(\mam_R)})\cong F^N_\bullet$ (Fact \ref{EhsanStuck}(i)) and such that  $G_\bullet$ satisfies condition (ii) of Lemma \ref{EhsanWhenWillThisEnd} (by Fact \ref{EhsanStuck}(ii)). Furthermore, from $G_\bullet\otimes_{R[\mathbf{X}]}(R[\mathbf{X}]_{(\mam_R)})\cong F^N_\bullet$  we conclude that 
		the complex $G_\bullet$ satisfies condition (i) of Lemma \ref{EhsanWhenWillThisEnd}. Moreover, since for each $1\le i\le n+s+2$ we have 
		$$\big(H_i(G_\bullet)\big)_{(\mam_R)}\cong H_i(G_\bullet)\otimes_{R[\mathbf{X}]}(R[\mathbf{X}]_{(\mam_R)})\cong H_i\big(G_\bullet\otimes_{R[\mathbf{X}]}(R[\mathbf{X}]_{(\mam_R)})\big)\cong H_i(F^N_\bullet)=0$$    so we may, and we do, find (and fix) some polynomial
		\begin{equation}
		\label{Elementg}
		g\in\Big(\bigcap\limits _{i=1}^{n+s+2}0:_{R[\mathbf{X}]}H_{i}(G_{\bullet})\Big)\backslash\mam_RR[\mathbf{X}].
		\end{equation}

		Hence $G_\bullet$ satisfies condition (iii) of Lemma \ref{EhsanWhenWillThisEnd}. By our hypothesis there is some $u\in \mn$ such that  for each $i\ge u$ 
		\begin{align*}
		\Big(\tor^{R[\mathbf{X}]}_{i}\big(H_0(G_\bullet),T\otimes_{R}R[\mathbf{X}]\big)\Big)_{(\mam_R)}
		&\cong
		\Big(\tor^{R[\mathbf{X}]}_{i}\big(H_0(G_\bullet),T\otimes_{R}R[\mathbf{X}]\big)\Big)\otimes_{R[\mathbf{X}]}R[\mathbf{X}]_{(\mam_R)}
		&\\&\cong
		\tor^{R[\mathbf{X}]_{(\mam_R)}}_{i}\big(H_0(G_\bullet\otimes_{R[\mathbf{X}]}R[\mathbf{X}]_{(\mam_R)}),T\otimes_{R}R[\mathbf{X}]_{(\mam_R)}\big)
		&\\&\cong
		\tor^{R[\mathbf{X}]_{(\mam_R)}}_{i}\big(N,T\otimes_RR[\mathbf{X}]_{(\mam_R)}\big)=0,
		\end{align*}
		and whence there is some    
		$\zeta_i\in R[\mathbf{X}]\backslash \mam_RR[\mathbf{X}]$   	   such that 
		\begin{equation}
		\label{EquationAnnihilatorElements}
		\zeta_i\in 0:_{R[\mathbf{X}]}\Big(\tor^{R[\mathbf{X}]}_{i}\big(H_0(G_\bullet),T\otimes_{R}R[\mathbf{X}]\big)\Big).
		\end{equation}
		Consequently, $G_\bullet$ also satisfies condition (iv) of Lemma \ref{EhsanWhenWillThisEnd}. 
		
		It remains to prove that when $T$ has finite complete intersection dimension, we can find a uniform annihilator  $\varsigma\in R[\mathbf{X}]\backslash \mam_RR[\mathbf{X}]$ for all  of Tor modules  $\tor^{R[\mathbf{X}]}_{i}\big(H_0(G_\bullet),T\otimes_{R}R[\mathbf{X}]\big)$ with $i\gg 0$.
		
		Let $R\overset{g}{\hookrightarrow} R'\twoheadleftarrow Q$ be a quasi-deformation of $R$ with $\pd_{Q}(R'\otimes_RT)<\infty$, in which we can assume without loss of generality that $R'$ is a complete local ring (Remark \ref{QuasiDeformationWithCompleteDeformation}). Then, $g_{\mathbf{X}}:R[\mathbf{X}]\overset{r\mapsto g(r),\ X_i\mapsto X_i}{\longrightarrow} R'[\mathbf{X}]$ is a faithfully flat extension (Fact \ref{FaithfullyFlat}(ii))  and 
		\begin{align*}
		\numberthis
		\label{EquationFinitenessProjectiveDimension}
		\pd_{Q[\mathbf{X}]}(T\otimes_RR'[\mathbf{X}])=\pd_{Q[\mathbf{X}]}\big(T\otimes_R(R'\otimes_QQ[\mathbf{X}])\big)&=\pd_{Q[\mathbf{X}]}\big((T\otimes_RR')\otimes_QQ[\mathbf{X}]\big)&\\&\overset{Q[\mathbf{X}]\text{\ is\ }Q\text{-flat}}{\le} \pd_Q(T\otimes_RR')<\infty.
		\end{align*}
		
		Moreover, from (\ref{EquationAnnihilatorElements})  and flatness of $g_{\mathbf{X}}$ we get 
		\begin{align*}
		\label{EquationIndexedAnnihilator}
		\numberthis
		g_{\mathbf{X}}(\zeta_i)&\in 0:_{R'[\mathbf{X}]}\Big(\tor^{R'[\mathbf{X}]}_i\big(H_0(G_\bullet),T\otimes_RR[\mathbf{X}]\big)\otimes_{R[\mathbf{X}]}R'[\mathbf{X}]\Big)&\\&=0:_{R'[\mathbf{X}]}\Big(\tor^{R'[\mathbf{X}]}_i\big(H_0(G_\bullet)\otimes_{R[\mathbf{X}]}R'[\mathbf{X}],T\otimes_RR'[\mathbf{X}]\big)\Big)
		\end{align*}
		for all $i\ge u$. Since $R'$ is a complete local ring so $R'[\mathbf{X}]$ is a homomorphic image of Gorenstein ring of finite Krull dimension and thus it admits a dualizing complex (\cite[(7.1.6) Remark]{ChristensenFoxbyHyperHomological}). So by Corollary \ref{CorollaryUniformTorAnnihilator} and the display (\ref{EquationFinitenessProjectiveDimension}) there is $t,v\in \mathbb{N}_0$ such that setting $w:=u+t$ we have 
		\begin{align*}
		\label{EquationCommonAnnihilator}
		\numberthis 
		\Big(\bigcap\limits_{i=u}^{w+v+1}0:_{R'[\mathbf{X}]}\tor^{R'[\mathbf{X}]}_i\big(H_0(G_\bullet)\otimes_{R[\mathbf{X}]}R'&[\mathbf{X}],T\otimes_RR'[\mathbf{X}]\big)\Big)^{2t+2}
		\subseteq &\\&0:_{R'[\mathbf{X}]}\tor^{R'[\mathbf{X}]}_j\big(H_0(G_\bullet)\otimes_{R[\mathbf{X}]}R'[\mathbf{X}],T\otimes_RR'[\mathbf{X}]\big),\ \ \forall\ j\ge w.
		\end{align*}
		Setting $\varsigma:=(\prod\limits_{i=u}^{v+w+1}\zeta_i)^{2t+2}$, we have $\varsigma\in R[\mathbf{X}]\backslash \mathfrak{m}_RR[\mathbf{X}]$ while $$g_{\mathbf{X}}(\varsigma)\in 0:_{R'[\mathbf{X}]}\tor^{R'[\mathbf{X}]}_j\big(H_0(G_\bullet)\otimes_{R[\mathbf{X}]}R'[\mathbf{X}],T\otimes_RR'[\mathbf{X}]\big),\ \ \forall\ j\ge w$$ by (\ref{EquationIndexedAnnihilator}) and (\ref{EquationCommonAnnihilator}). Consequently, $\varsigma\in 0:_{R[\mathbf{X}]}\text{Tor}^{R[\mathbf{X}]}_{j}\big(H_0(G_\bullet),T\otimes_RR[\mathbf{X}]\big)$ for all $j\ge w$ by Fact \ref{FactFlatAnnihilator}(ii) as $g_\mathbf{X}$ is faithfully flat. 
	\end{proof}
\end{sublem}

Now we can deduce the main result of this subsection.

\begin{subcor}\label{CorollaryTranscendentalExtension}
	Let $R$ be a  complete local ring, $T$ be a test module for $R$ and $\mathbf{X}:=X_1,\ldots,X_n$ be a  sequence of indeterminates over $R$. Then $T\otimes_RR[\mathbf{X}]_{(\fm_R)}$ is a test module for $R[\mathbf{X}]_{(\fm_R)}$ provided one of the following conditions holds:
	\begin{enumerate}
		\item [(i)] There exists a flat  local $R$-algebra $R'$ such that $K_{R'}$ is uncountable and $T\otimes_RR'$ is an $R'$-test module (evidently, one can set $R':=R$ provided $K_R$ is uncountable).
		\item[(ii)] $\text{CI-dim}_R(T)<\infty$.
	\end{enumerate}
	\begin{proof}
		(i) Using our hypothesis we can, without loss of generality, assume that $K_R$ is uncountable. Namely, $R[\mathbf{X}]_{(\fm_R)}\rightarrow R'[\mathbf{X}]_{(\fm_{R'})}$ is also a faithfully flat  ring homomorphism because by Fact \ref{FaithfullyFlat}(ii) $R[\mathbf{X}]\rightarrow R'[\mathbf{X}]$ is (faithfully)  flat while   $(\fm_{R'}R'[\mathbf{X}])\cap R[\mathbf{X}]=\fm_RR[\mathbf{X}]$ (as $\fm_{R'}\cap R=\fm_R$). So setting $T':=T\otimes_R R',$ if we can prove that $$T'\otimes_{R'}R'[\mathbf{X}]_{(\fm_{R'})}\cong T\otimes_RR'[\mathbf{X}]_{(\fm_{R'})}\cong  (T\otimes_RR[\mathbf{X}]_{(\fm_R)})\otimes_{R[\mathbf{X}]_{(\fm_{R})}}R'[\mathbf{X}]_{(\fm_{R'})}$$ is a test module for $R'[\mathbf{X}]_{(\fm_{R'})}$ then we can conclude that $T\otimes_R R[\mathbf{X}]_{(\fm_R)}$ is also a test module for $R[\mathbf{X}]_{(\fm_R)}$ (see Remark \ref{RemarkTestModulePropertyDesentFromFaithfullyFlatExtension}).

		Suppose to the contrary that $T\otimes_R(R[\mathbf{X}]_{(\fm_R)})$ is not a test module hence there is $N\in \text{mod}(R[\mathbf{X}]_{(\fm_R)})$  such that $\pd_{R[\mathbf{X}]_{(\fm_R)}}(N)=\infty$ and $\tor^{R[\mathbf{X}]_{(\fm_R)}}_{\gg}(N,T\otimes_RR[\mathbf{X}]_{(\fm_R)})=0$. Then by Lemma \ref{LemmaDesiredPerfectComplexExists}  there is a perfect complex $G_\bullet\in \mathcal{D}_b(R[\mathbf{X}])$  satisfies the condition (i)-(iv) of 
		Lemma \ref{EhsanWhenWillThisEnd}. In view of Lemma \ref{EhsanWhenWillThisEnd}, in order to get our desired contradiction it suffices to find a sequence $c_1,\ldots,c_n$ of elements of $R$ such that $\zeta_j(c_1,\ldots,c_n)\in R\backslash \fm_R$ for any $j\ge u$ and $g(c_1,\ldots,c_n)\in R\backslash \fm_R$. This is equivalent to find some elements $k_1,\ldots,k_n\in K_R$ such that $\overline{\zeta_j}(k_1,\ldots,k_n)\neq 0$ and $\overline{g}(k_1,\ldots,k_n)\neq 0$ where $\overline{\zeta_j}$ (resp. $\overline{g}$) is the image of $\zeta_j$ (resp. $g$) under the natural ring epimorphism $R[\mathbf{X}]\rightarrow K_R[\mathbf{X}]$ induced by the natural epimorphism $R\twoheadrightarrow K_R$. However, such a sequence $k_1,\ldots,k_n\in K_R$ exists in view of Fact \ref{PolynomialFact} because $\overline{g}\neq 0$  and $\overline{\zeta_j}\neq 0$ for each $j\ge u$ (recall that $g\notin \fm_RR[\mathbf{X}]$ and $\zeta_j\notin \fm_RR[\mathbf{X}]$ ($j\ge u$)). 
		
		(ii)   If $\text{CI-dim}_R(T)=-\infty$, i.e. $T=0$, then   $R$ (and so $R[\mathbf{X}]_{(\fm_R)}$) is a   regular local ring and there is nothing to prove. So, we suppose that $\text{CI-dim}_R(T)$ is finite.
		
		 Arguing as in the proof of the previous part we assume, to the contrary, that $T\otimes_RR[\mathbf{X}]_{(\fm_R)}$ is not a test module  and we arrive at a perfect complex $G_\bullet\in \mathcal{D}_b(R[\mathbf{X}])$  satisfying parts (i)-(iv) in Lemma \ref{EhsanWhenWillThisEnd}. Even more,  since it is assumed that $\text{CI-dim}_R(T)<\infty$ (and $T\neq 0$) so by Lemma \ref{LemmaDesiredPerfectComplexExists}(ii) we can assume that $G_\bullet$ satisfies an  improved version of part (iv) in Lemma \ref{EhsanWhenWillThisEnd}. Namely, without loss of generality we can assume, about the elements $\zeta_j$ in part (iv) in Lemma \ref{EhsanWhenWillThisEnd}, that $\zeta_j=\varsigma$ for each $i\ge u$ and for some  (fixed) $\varsigma\in R[\mathbf{X}]\backslash \fm_RR[\mathbf{X}]$. Hence, for getting our desired contradiction, following the statement of Lemma \ref{EhsanWhenWillThisEnd} it suffices  to find a sequence of elements $c_1,\ldots,c_n$ of $R$ such that $\varsigma(c_1,\ldots,c_n)\neq 0$ and $g(c_1,\ldots,c_n)\neq 0$. If  $K_R$ is an infinite field then there is $k_1,\ldots,k_n\in K_R$ such that $\overline{g}\overline{\varsigma}(k_1,\ldots,k_n)\neq 0$, so we set the desired sequence $c_1,\ldots,c_n\in R$ to be any lift of $k_1,\ldots,k_n$.  So the statement holds if $K_R$ is an infinite field.
		
		Suppose that $K_R$ is a finite field. Then  there is a finite algebraic field extension $h:K_R\rightarrow \mathcal{K}$  of $K_R$  and a sequence $\kappa_1,\ldots,\kappa_n\in \mathcal{K}$ such that $\overline{g}\overline{\varsigma}(\kappa_1,\ldots,\kappa_n)\neq 0$ (see Fact \ref{PolynomialOverFiniteField}). Since $R$ is a complete local ring so $R$ admits a coefficient ring $\lambda_R:C\rightarrow R$ which induces the isomorphism $\overline{\lambda_R}:K_C\rightarrow K_R$ by definition. 
		
		If $R$ is equicharacteristic then $C$ is a field and we set $C':=\mathcal{K}$. Afterwards, we consider local $C$-algebra $\varphi:C\rightarrow C'$ as the composition of ring homomorphisms $C\overset{\lambda_R}{\rightarrow}R\overset{\pi}{\twoheadrightarrow}K_R\overset{h}{\longrightarrow}\mathcal{K}$ where $\pi$ is the natural epimorphism (i.e. $\varphi=h\circ \pi \circ  \lambda_R$, and it is  a finite algebraic field extension as $\pi\circ \lambda_R$ is an isomorphism). 
		
		If $R$ has mixed characteristic then $C$ is a $p$-ring where $p$ is the characteristic of $K_R$. Then by Fact \ref{RemarkMixedCharacteristicPRingExtension} there is a $p$-ring $C'$ and an   extension  of $p$-rings $\varphi:C\rightarrow C'$ such that $\overline{\varphi}\circ \overline{\lambda_R}^{\ -1}:K_R\rightarrow K_{C'}$ and $h:K_R\rightarrow \mathcal{K}$ are isomorphic $K_R$-algebras. 		
		 In particular, $K_{C'}$ is a finite field extension of $K_R$ and  there is a sequence of elements $\kappa'_1,\ldots,\kappa'_n\in K_{C'}$ such that $\overline{g}\overline{\varsigma}(\kappa'_1,\ldots,\kappa'_n)\neq 0$. 
		
		In any case (mixed characteristic or equicharacteristic), using the local $C$-algebras $\lambda_R:C\rightarrow R$ and $\varphi:C\rightarrow C'$  and Fact \ref{FactCompleteTensorProductFactCoefficientRingExtension}  we can take into account  the complete tensor product $R\widehat{\otimes}_CC'$ and  the weakly unramified  flat local homomorphism  $$\eta_{R,R\widehat{\otimes}_CC'}:R\rightarrow R\widehat{\otimes}_CC'$$ which is a residually finite  local homomorphism as $K_{R\widehat{\otimes}_CC'}$ is $K_R$-algebra isomorphic to the finite algebraic field extension  $\overline{\varphi}\circ \overline{\lambda_R}^{\ -1}:K_R\rightarrow K_{C'}$. 
		
		Let us set $R_{\mathcal{K}}:=R\widehat{\otimes}_{C}C'$ and $\eta:=\eta_{R,R\widehat{\otimes}_{C}C'}$. Then, $\eta_{\mathbf{X}}:R[\mathbf{X}]\overset{r\mapsto \eta(r),\ X_i\mapsto X_i}{\longrightarrow} R_{\mathcal{K}}[\mathbf{X}]$ is also a faithfully flat ring homomorphism (Fact \ref{FaithfullyFlat}(ii)).  Since $\eta_{\mathbf{X}}$ and $\eta$ are faithfully flat ring homomorphisms so $P_\bullet:=G_\bullet\otimes_{R[\mathbf{X}]}(R_{\mathcal{K}}[\mathbf{X}])$ is a perfect complex in $\mathcal{D}_b(R_{\mathcal{K}}[\mathbf{X}])$ such that it fulfills te following conditions:
		\begin{enumerate}
			\item[(i)]   $P_{i}= 0$  if $i<0$ and  $P_i\neq 0$   if $0\le i\le n+s+3$ where $n$ is the number of indeterminates in $\mathbf{X}$ ($s=\depth(R)=\depth(R_{\mathcal{K}})$).
			\item[(ii)] All entries of the matrices of the  differentials of $P_\bullet$ belong to $\mam_RR_{\mathcal{K}}[\mathbf{X}]=\mam_{R_{\mathcal{K}}}R_{\mathcal{K}}[\mathbf{X}]$.
			\item [(iii)]   $\eta_{\mathbf{X}}(g)\in R_{\mathcal{K}}[\mathbf{X}]\backslash \mam_{R_{\mathcal{K}}}R_{\mathcal{K}}[\mathbf{X}]$  and  $$\eta_{\mathbf{X}}(g)\in \bigcap\limits_{i=1}^{n+s+2} 0:_{R_{\mathcal{K}}[\mathbf{X}]}H_i(P_{\bullet}).$$ 
			\item[(iv)] We have $\eta_{\mathbf{X}}(\varsigma)\in R_{\mathcal{K}}[\mathbf{X}]\backslash \fm_{R_\mathcal{K}}R_{\mathcal{K}}[\mathbf{X}]$ and $$\eta_{\mathbf{X}}(\varsigma)\in 0:_{R_{\mathcal{K}}[\mathbf{X}]}\Big(\text{Tor}^{R_{\mathcal{K}}[\mathbf{X}]}_i\big(H_0(P_\bullet),(T\otimes_RR_{\mathcal{K}})\otimes_{R_\mathcal{K}}R_{\mathcal{K}}[\mathbf{X}]\big)\Big),\ \ \forall\ i\gg 0.$$
		\end{enumerate}
		Hence Lemma \ref{EhsanWhenWillThisEnd} implies that $T\otimes_RR_{\mathcal{K}}$ is not a test module for $R_{\mathcal{K}}$ because $\overline{\eta_{\mathbf{X}}(g)}\overline{\eta_{\mathbf{X}}(\varsigma)}(\kappa''_1,\ldots,\kappa''_n)\neq 0$ for some $\kappa''_1,\ldots,\kappa''_n\in K_{R_{\mathcal{K}}}$ (again  because $\overline{\eta_{R,R\widehat{\otimes}_{C}C'}}:K_R\rightarrow K_{R_{\mathcal{K}}},\ \overline{\varphi}\circ \overline{\lambda_R}^{\ -1}:K_R\rightarrow K_{C'}$ and $h:K_R\rightarrow \mathcal{K}$ are all isomorphic $K_R$-algebras).
		But, this is a contradiction because  $T\otimes_{R}R_{\mathcal{K}}$ is a test module for $R_{\mathcal{K}}$ by \cite[Theorem 3.5]{CelikbasWagstaff}.
	\end{proof}
\end{subcor}

\subsection{Test modules and flat coefficient ring base changes}

In this subsection  we generalize the main result of the previous subsection concerning  purely transcendental  extensions  to the context of flat coefficient ring base changes, namely to the extensions of the form $R\widehat{\otimes}_CC'$ where $C$ is a coefficient ring for the complete local ring $R$ and $C\rightarrow C'$ is either an arbitrary field extension or a (residually purely transcendental)  extension of unramified discrete valuation rings (see Corollary \ref{BasedFieldExtensionLemma}). We shall begin by considering residually finitely generated case, but we need first the following  lemma which will be used in the sequel a few times.

\begin{sublem}\label{LemmaFactorizationPRingExtensionToResiduallyTranscendentalAndResiduallyAlgebraic}
	Suppose that $C$ is a $p$-ring and $\varphi:C\rightarrow C'$ is an extension of $p$-rings. Let $\{\kappa_i\}_{i\in I}$ be a set of elements of $C'$ such that the set of residue classes $\{\overline{\kappa}_i\}_{i\in I}\subseteq K_{C'}$ is transcendental over $\overline{\varphi}(K_C)$. Consider the set of indeterminates $\mathbf{X}:=\{X_i\}_{i\in I}$ and the ring homomorphism $$\varphi_{\mathbf{X}}:C[\mathbf{X}]\rightarrow C',\ \ c\mapsto \varphi(c),\ X_i\mapsto \kappa_i.$$
	\begin{enumerate}
		\item [(i)] Let $C[\mathbf{X}]_{(p)}^{\widehat{\ }_{(p)}}$ be the $(p)$-adic completion  of $C[\mathbf{X}]_{pC[\mathbf{X}]}$. Then  $C[\mathbf{X}]_{(p)}^{\widehat{\ }_{(p)}}$ is a $p$-ring.
		\item[(ii)] $\varphi_{\mathbf{X}}$ induces the ring homomorphism $\phi:C[\mathbf{X}]_{(p)}^{\widehat{\ }_{(p)}}\rightarrow C'$.
		\item[(iii)] Let $\psi:C\rightarrow C[\mathbf{X}]_{(p)}^{\widehat{\ }_{(p)}}$ be the natural embedding $c\mapsto \big(c/1+(p^n)\big)_{n\in\mathbb{N}}$. Then $\phi\circ \psi=\varphi$. 
		\item[(iv)] Let $\mathcal{L}$ be the transcendental extension of $\overline{\varphi}(K_C)$ generated by $\{\overline{\kappa}_i\}_{i\in I}$. Then there is an isomorphism $\gamma:K_{C[\mathbf{X}]_{(p)}^{\widehat{\ }_{(p)}}}\rightarrow \mathcal{L}$ that fits into the following commutative diagram:
		 \[
		   \xymatrix{
		   	  K_C \ar[rrrr]^{\overline{\psi}} \ar[d]_{\overline{\varphi}}  &&&& K_{C[\mathbf{X}]_{(p)}^{\widehat{\ }_{(p)}}} \ar[d]^{\gamma}_{\cong} \ar[lllld]_{\overline{\phi}}\\
		   	  K_{C'} &&&& \mathcal{L} \ar[llll]^{\text{inclusion}} 
		   }
		 \]
	\end{enumerate}
  \begin{proof}
  	  (i) This is obvious ($C[\mathbf{X}]_{pC[\mathbf{X}]}$ is  a $1$-dimensional regular local domain  whose maximal ideal is generated by $p$, its residual characteristic).
  	  
  	  (ii) and (iii) If we can show that $\varphi_{\mathbf{X}}(g)$ is invertible in $C'$ for each $g\in C[\mathbf{X}]\backslash pC[\mathbf{X}]$, then in view of the  universal property of localization   $\varphi_{\mathbf{X}}$ induces the ring homomorphism $C[\mathbf{X}]_{pC[\mathbf{X}]}\rightarrow C'$ such that it takes $f/g$ to $\varphi_{\mathbf{X}}(f)\varphi_{\mathbf{X}}(g)^{-1}$. Afterwards, we get the desired ring homomorphism $\phi$ as the composition $$C[\mathbf{X}]_{(p)}^{\widehat{\ }_{(p)}}=\lim\limits_{\underset{n\in \mathbb{N}}{\longleftarrow}}\big(C[\mathbf{X}]_{pC[\mathbf{X}]}/(p^n)\big)\rightarrow \lim\limits_{\underset{n\in \mathbb{N}}{\longleftarrow}} C'/p^nC'\overset{\cong}{\rightarrow}C'$$
  	  because $C'$ is complete ($\lim\limits_{\underset{n\in \mathbb{N}}{\longleftarrow}} C'/p^nC'\rightarrow C'$ is the inverse of the canonical map $C'\rightarrow \lim\limits_{\underset{n\in \mathbb{N}}{\longleftarrow}} C'/p^nC'$, and the first map is such that $\big(g_n/1+(p^n)\big)_{n\in\mathbb{N}}\mapsto \big(\varphi_{\mathbf{X}}(g_n)+(p^n)\big)_{n\in\mathbb{N}}$). Then part (iii) also follows. Hence to complete the proof it suffices to consider some $g\in C[\mathbf{X}]\backslash pC[\mathbf{X}]$ and then to show that $\varphi_{\mathbf{X}}(g)$ is invertible. Note that the residue class of $\varphi_{\mathbf{X}}(g)$ in $K_{C'}$ is a polynomial in the transcendental basis $\{\overline{\kappa}_i\}_{i\in I}$ with coefficients in $\overline{\varphi}(K_{C})$, and with some non-zero coefficient  by our assumption that $g\in C[\mathbf{X}]\backslash pC[\mathbf{X}]$. Therefore, the residue class of $\varphi_{\mathbf{X}}(g)$ in $K_{C'}$ is non-zero in view of the transcendental property of $\{\overline{\kappa}_i\}_{i\in I}$ over $\overline{\varphi}(K_C)$, i.e. $\varphi_{\mathbf{X}}(g)\notin pC'$ as was to be proved.
  	  
  	  (iv) The commutativity of the upper triangle in the diagram is obvious in view of part (iii). For the isomorphism $\gamma$ and the commutativity of the lower triangle, it suffices to notice that  $\mathcal{L}$ coincides with the image of  the  induced map $\overline{\phi}:K_{C[\mathbf{X}]_{(p)}^{\widehat{\ }_{(p)}}}\rightarrow K_{C'}$, by $\phi$ (note that $K_{C[\mathbf{X}]_{(p)}}\rightarrow K_{C[\mathbf{X}]_{(p)}^{\widehat{\ }_{(p)}}}$ is an isomorphism, so to compute the image of $\overline{\phi}$, it suffices to compute the image of $K_{C[\mathbf{X}]_{(p)}}\rightarrow K_{C'}$).
  \end{proof}
\end{sublem}

\begin{sublem}\label{FinitelyGeneratedExension}  Suppose that $R$ is a complete local ring with a coefficient ring $\lambda_R:C\rightarrow R$. Assume that $\varphi:C\rightarrow C'$ is either an extension of  fields or an extension of $p$-rings such that $\overline{\varphi}:K_C\rightarrow K_{C'}$ is a finitely generated field extension.  Consider the flat local weakly unramified homomorphism $\eta_{R,R\widehat{\otimes}_CC'}:R\rightarrow R\widehat{\otimes}_CC'$ where the complete tensor product is assigned to $\lambda_R$ and $\varphi$. If $T$ is an  test module for $R$ and either of the two conditions in the statement of Corollary \ref{CorollaryTranscendentalExtension} holds, then $T\otimes_R(R\widehat{\otimes}_CC')$ is an $R\widehat{\otimes}_CC'$-test module.
	\begin{proof}
		Since $K_{C'}$ is a finitely generated field extension of $K_C$ by our hypothesis, so there is a subfield $\mathcal{L}=\overline{\varphi}(K_C)(\overline{\kappa}_1,\ldots,\overline{\kappa}_n)$ of $K_{C'}$ such that $\{\overline{\kappa}_i\}_{i=1}^n\subseteq K_{C'}$ is transcendental over $\overline{\varphi}(K_C)$ while $K_{C'}$ is a finite algebraic extension of $\mathcal{L}$ ($\overline{\kappa}_i$ is the residue class of $\kappa_i\in C'$). 
		
		 To make the proof more readable we separate the equicharacteristic case of the proof from the mixed characteristic case. However, both cases follow the same procedure and we could present a structural proof working for both cases.
		 
		 \underline{Case 1, $R$ has equicharacteristic:} So $C$ and $C'$ are both fields, $K_C=C$, $K_{C'}=C'$ and $\overline{\kappa}_i=\kappa_i$ $(1\le i\le n)$. Considering the sequence of indeterminates $\mathbf{X}:=X_1,\ldots,X_n$ over $C$, due to our hypothesis we have $C(\mathbf{X})\cong\mathcal{L}$ by the map $\phi(X_i)=\kappa_i$ and $\phi(c)=\varphi(c)$. Moreover, it is easily seen that there is an isomorphism $\theta:R\otimes_CC(\mathbf{X})\rightarrow S^{-1}(R[\mathbf{X}])$ of $R$-algebras where $S$ is the image of the multiplicatively closed set $C[\mathbf{X}]\backslash \{0\}$ under the ring homomorphism $C[\mathbf{X}]\overset{X_i\mapsto X_i,\ c\mapsto \lambda_R(c)}{\rightarrow}R[\mathbf{X}]$. In particular, being isomorphism of $R$-algebras we have $\theta(\mathfrak{M})=\theta\Big(\mathfrak{m}_R\big(R\otimes_CC(\mathbf{X})\big)\Big)=\fm_RS^{-1}(R[\mathbf{X}])$ while $\mathfrak{M}\in \text{Max}\big(R\otimes_CC(\mathbf{X})\big)$ by Fact \ref{FactCompleteTensorProductFactCoefficientRingExtension}(i). It follows that $\big(R\otimes_CC(\mathbf{X})\big)_{\mathfrak{M}}\cong R[\mathbf{X}]_{(\fm_R)}$ as $R$-algebras. Consequently, $T\otimes_R(R\otimes_CC\big(\mathbf{X})\big)_{\mathfrak{M}}$ is a test module for over  $(R\otimes_CC\big(\mathbf{X})\big)_{\mathfrak{M}}$ by Corollary \ref{CorollaryTranscendentalExtension}. So $T\otimes_R\big(R\widehat{\otimes}_CC(\mathbf{X})\big)$ is a test module over $R\widehat{\otimes}_CC(\mathbf{X})$ by \cite[Theorem 3.5]{CelikbasWagstaff} and Proposition \ref{CompleteTensorProductNoetherianness}(iii).
		 
		 Since $C'$ is finite and algebraic over its subfield $\mathcal{L}$, so $\Phi:C(\mathbf{X})\overset{\phi(\cong)}{\rightarrow } \mathcal{L}	\overset{\text{inclusion}}{\rightarrow}C'$ is also a finite algebraic extension of fields. Therefore, considering local $C(\mathbf{X})$-algebras $\Phi$ and $\eta_{C(\mathbf{X}),R\widehat{\otimes}_CC(\mathbf{X})}$ from parts (ii) and (iii) of Fact \ref{FactCompleteTensorProductFactCoefficientRingExtension} we deduce that
		 $$ 
		 \eta_{R\widehat{\otimes}_CC(\mathbf{X}),\big(R\widehat{\otimes}_CC(\mathbf{X})\big)\widehat{\otimes}_{C(\mathbf{X})}C'}:R\widehat{\otimes}_{C}C(\mathbf{X})\rightarrow \big(R\widehat{\otimes}_CC(\mathbf{X})\big)\widehat{\otimes}_{C(\mathbf{X})}C'
		 $$
		 is a residually finite flat 	weakly unramified extension (note that by Fact \ref{FactCompleteTensorProductFactCoefficientRingExtension}(iv), $\eta_{C(\mathbf{X}),R\widehat{\otimes}_CC(\mathbf{X})}$ is a coefficient ring for $R\widehat{\otimes}_CC(\mathbf{X})$). Thus, applying \cite[Theorem 3.5]{CelikbasWagstaff} to the residually finite weakly unramified homomorphism $\eta_{R\widehat{\otimes}_CC(\mathbf{X}),\big(R\widehat{\otimes}_CC(\mathbf{X})\big)\widehat{\otimes}_{C(\mathbf{X})}C'}$, from the concluding assertion of the previous paragraph we deduce that 
		 $$
		 \Big(T\otimes_R\big(R\widehat{\otimes}_CC(\mathbf{X})\big)\Big)\otimes_{R\widehat{\otimes}_CC(\mathbf{X})}\Big(\big(R\widehat{\otimes}_CC(\mathbf{X})\big)\widehat{\otimes}_{C(\mathbf{X})}C'\Big) \cong T\otimes_R\Big(\big(R\widehat{\otimes}_CC(\mathbf{X})\big)\widehat{\otimes}_{C(\mathbf{X})}C'\Big)$$
		 is a test module over $\big(R\widehat{\otimes}_CC(\mathbf{X})\big)\widehat{\otimes}_{C(\mathbf{X})}C'$. 
		 However, there is an isomorphism of $R$-algebras 
		 $$
		   \Big(R\widehat{\otimes}_CC(\mathbf{X})\Big)\widehat{\otimes}_{C(\mathbf{X})}C' \overset{\text{Proposition \ref{Associativity}}}{\cong} R\widehat{\otimes}_C\big(C(\mathbf{X})\widehat{\otimes}_{C(\mathbf{X})}C'\big)\overset{\text{Remark \ref{RemarkNaturalIsomorphism}}}{\cong}R\widehat{\otimes}_CC'.
		 $$
		 Consequently, $T\otimes_R(R\widehat{\otimes}_CC')$ is a test module over $R\widehat{\otimes}_CC'$ as was to be proved.
		 
		 \underline{Case 2, $R$ has mixed characteristic:} We consider the sequence of indeterminates $\mathbf{X}:=X_1,\ldots,X_n$ over $C$ and we denote the $\big(p(C[\mathbf{X}])_{pC[\mathbf{X}]}\big)$-adic completion of $C[\mathbf{X}]_{pC[\mathbf{X}]}$ by $C_\mathcal{L}$. In what follows we consider $C[\mathbf{X}]_{(p)}$ as an $C$-algebra in the natural way, and we consider $C_\mathcal{L}$ also as an $C$-algebra in the natural way via the ring homomorphism $\psi$ in Lemma \ref{LemmaFactorizationPRingExtensionToResiduallyTranscendentalAndResiduallyAlgebraic}(iii).
		 
		 Note that $C[\mathbf{X}]_{(p)}$ is an unramified discrete valuation ring of mixed characteristic $(0,p)$. We have the isomorphism of $R$-algebras $\theta:R\otimes_CC[\mathbf{X}]_{(p)}\rightarrow S^{-1}(R[\mathbf{X}])$  where $S$ is the image of $C[\mathbf{X}]\backslash (p)$ under $C[\mathbf{X}]\overset{X_i\mapsto X_i,\ c\mapsto \lambda_R(c)}{\longrightarrow}R[\mathbf{X}]$. Consequently, we have the $R$-algebra isomorphism $(R\otimes_CC[\mathbf{X}])_{(p)})_{\mathfrak{M}}\cong R[\mathbf{X}]_{(\fm_R)}$ as in the equicharacteristic case, implying that $T\otimes_R(R\otimes_CC[\mathbf{X}]_{(p)})_{\mathfrak{M}}$ is a test module over $(R\otimes_CC[\mathbf{X}]_{(p)})_{\mathfrak{M}}$ in view of Corollary \ref{CorollaryTranscendentalExtension}. Therefore $T\otimes_R(R\widehat{\otimes}_{C}C[\mathbf{X}]_{(p)})$ is a test module over  $$R\widehat{\otimes}_{C}C[\mathbf{X}]_{(p)}$$ by  \cite[Theorem 3.5]{CelikbasWagstaff} and Proposition \ref{CompleteTensorProductNoetherianness}(iii). Equivalently, $T\otimes_R(R\widehat{\otimes}_CC_{\mathcal{L}})$ is a test module over $R\widehat{\otimes}_CC_{\mathcal{L}}$ (see Remark \ref{RemarkDefinitionCompleteTensorProductDoesntAffectCompleteness}).
		 
		 Now we can consider the map $\phi:C_\mathcal{L}\rightarrow C'$ as in Lemma \ref{LemmaFactorizationPRingExtensionToResiduallyTranscendentalAndResiduallyAlgebraic}(ii). Note that by Lemma \ref{LemmaFactorizationPRingExtensionToResiduallyTranscendentalAndResiduallyAlgebraic}(iv) and the commutative diagram therein, $\overline{\phi}$ is a finite algebraic extension of fields (as so is $\mathcal{L}\rightarrow K_{C'}$).  Moreover, $\eta_{C_{\mathcal{L}},R\widehat{\otimes}_CC_{\mathcal{L}}}$ is a coefficient ring for $R\widehat{\otimes}_CC_{\mathcal{L}}$ by Fact \ref{FactCompleteTensorProductFactCoefficientRingExtension}(iv). Thus considering the complete tensor product assigned to to the local homomorphism $\eta_{C_{\mathcal{L}},R\widehat{\otimes}_CC_{\mathcal{L}}}$ and $\phi$, from parts (ii) and (iii) of Fact \ref{FactCompleteTensorProductFactCoefficientRingExtension} we conclude that
		 $$
		 \eta_{R\widehat{\otimes}_CC_{\mathcal{L}},(R\widehat{\otimes}_CC_{\mathcal{L}})\widehat{\otimes}_{C_\mathcal{L}}C'}:
		 R\widehat{\otimes}_CC_{\mathcal{L}}\rightarrow (R\widehat{\otimes}_CC_{\mathcal{L}})\widehat{\otimes}_{C_\mathcal{L}}C'
		 $$
		 is a residually finite weakly unramified flat local homomorphism. Thus \cite[Theorem 3.5]{CelikbasWagstaff} as well as the concluding assertion of the previous paragraph implies that $T\otimes_R\big((R\widehat{\otimes}_CC_{\mathcal{L}})\widehat{\otimes}_{C_{\mathcal{L}}}C'\big)$ is a test module over  $(R\widehat{\otimes}_CC_{\mathcal{L}})\widehat{\otimes}_{C_{\mathcal{L}}}C'$. Finally, from the $R$-algebra isomorphisms
		 $$
		 (R\widehat{\otimes}_{C}C_\mathcal{L})\widehat{\otimes}_{C_\mathcal{L}}C'\overset{\text{Proposition \ref{Associativity}}}{\cong}R\widehat{\otimes}_C(C_{\mathcal{L}}\widehat{\otimes}_{C_\mathcal{L}}C') \overset{\text{Remark \ref{RemarkDefinitionCompleteTensorProductDoesntAffectCompleteness}}}{\cong}R\widehat{\otimes}_CC'
		 $$
		 the statement follows (the isomorphism $C_{\mathcal{L}}\widehat{\otimes}_{C_{\mathcal{L}}}C'\cong C'$ is an isomorphism of $C_\mathcal{L}$-algebras, hence to get our desired final result here we are also using the fact that $C\overset{\psi}{\rightarrow}C_{\mathcal{L}}\overset{\phi}{\rightarrow}C'$ coincides with $\varphi:C\rightarrow C'$, see Lemma \ref{LemmaFactorizationPRingExtensionToResiduallyTranscendentalAndResiduallyAlgebraic}(iii)).
	\end{proof}
\end{sublem}

In the next corollary, we aim to relax the ``residually finitely generated condition"  in the statement of Lemma \ref{FinitelyGeneratedExension}. 
 But, we first give  the following two preparatory lemmas.

\begin{sublem}\label{LemmaTestModuleDirectLimit}
	Suppose that $R$ is a local ring and $A$ is a local $R$-algebra admitting a direct limit presentation $A=\lim\limits_{\underset{i\in I}{\longrightarrow}}A_i$  assigned to a direct system $(A_i,\ \varphi_{i,j}:A_i\rightarrow A_{j})_{\substack{i,j\in I\\ i\le j}}$ of local $R$-algebras and  flat local $R$-algebra homomorphisms. Assume that $T$ is an $R$-module such that $T\otimes_RA_i$ is an $A_i$-test module for each $i\in I$. Then $T\otimes_RA$ is a test $A$-module. 
	  \begin{proof}
	  	   It is easily seen $\fm_A=\lim\limits_{\underset{i\in I}{\longrightarrow}}\mam_{A_i}$.
	  	   Moreover, it is easily verified that the ring homomorphism to the direct limit $\varphi_i:A_i\rightarrow A$ is also a flat local  homomorphism, because $A=\lim\limits_{\underset{\substack{j\in I\\j\ge i}}{\longrightarrow}}A_j$ and each $A_j$ is a flat local $A_i$-algebra by our hypothesis for each $j\ge i$ (here we are using the fact that direct limit of flat modules is again flat). We consider a   finitely generated $A$-module $N$ such that $\tor^A_{\gg}(N,T\otimes_RA)=0$   and we prove that $\pd_A(N)<\infty$. 
	  	   
	  	   We take into account a presentation $A^n\overset{H}{\rightarrow}A^m\rightarrow N\rightarrow 0$ of $N$ where $H=[a_{l,k}]$ is a matrix representing the homomorphism  $A^n\rightarrow A^m$ of finitely generated free modules. For some $i\in I$, there exist elements $\alpha_{l,k}\in A_i$ such that $a_{l,k}=\varphi_i(\alpha_{l,k})$ for each $1\le l\le m$ and $1\le k\le n$. Let $H_{i}=[\alpha_{l,k}]$ be the matrix corresponding to $H$  (which has entries in $A_i$ by our choice of $\alpha_{l,k}$) and let $N_i$ be the cokernel of the homomorphism of finitely generated free modules $A_i^n\overset{H_i}{\rightarrow} A_i^m$. Then $N_i\otimes_{A_i}A\cong N$. Therefore, since $A$ is a faithfully flat $A_i$-algebra we have $\tor^{A_i}_{\gg}(N_i,T\otimes_RA_i)=0$, as $\tor^{A_i}_j(N_i,T\otimes_RA_i)\otimes_{A_i}A\cong \tor^A_j(N,T\otimes_RA)=0$ for any $j\gg 0$. Consequently, $\pd_{A_i}(N_i)<\infty$ because $T\otimes_RA_i$ is a test module over $A_i$ by our assumption. But then $$\pd_{A}(N)=\pd_{A}(N_i\otimes_{A_i}A)\overset{\varphi_{i}\text{\ is faithfully flat}}{=}\pd_{A_i}(N_i)<\infty.\qedhere$$
	  \end{proof}
\end{sublem}

\begin{sublem}\label{LemmaCoefficientRingBaseChange}
	Suppose that $R$ is a complete local ring and $\lambda_R:C\rightarrow R$ is a coefficient ring.  Let $T$ be a test module over $R$. Assume that either of the two conditions in the statement of Corollary \ref{CorollaryTranscendentalExtension} holds.
	\begin{enumerate} 
		\item[(i)] Suppose  that $\varphi:C\rightarrow C'$ is  an extension of fields (thus $R$ has equicharacteristic). Assume that  $C'=\lim\limits_{\underset{i\in I}{\longrightarrow }}V_i$   (as $C$-algebras)  such that each $V_i$ is a finitely generated field extension of $C$ via $\theta_i:C\rightarrow V_i$ and the transitions maps $\varphi_{i,j}:V_i\rightarrow V_j$ satisfy $\varphi_{i,j}\circ \theta_i=\theta_j$ (i.e. $\varphi_{i,j}$ is a homomorphism of $C$-algebras). Then, $T\otimes_R (R\widehat{\otimes}_CC')$ is a test module over $R\widehat{\otimes}_CC'$.
		
		\item[(ii)]  Suppose that $\varphi:C\rightarrow C'$ is an  extension of mixed characteristic $(0,p)$ unramified discrete valuation rings (thus $R$ has mixed characteristic and $K_R$ has prime characteristic $p>0$). Assume that  $C'=\lim\limits_{\underset{i\in I}{\longrightarrow }}V_i$ (as $C$-algebras) such that each  $V_i$ is an unramified discrete valuation ring of mixed characteristic $(0,p)$ that is an $C$-algebra via $\theta_i:C\rightarrow V_i$. Presume that the transitions maps $\varphi_{i,j}:V_i  \rightarrow V_j$ satisfy $\varphi_{i,j}\circ \theta_i=\theta_j$  and that $V_i$ is essentially of finite type over $C$ (via $\theta_i$)  (for each $i,j\in I,\ i\le j$). Then $T\otimes_R (R\widehat{\otimes}_CC')$ is a test module over $R\widehat{\otimes}_CC'$.
		\end{enumerate}
	\begin{proof}
		 The next  paragraph works for both cases (i) and (ii) of the proof.
		 
		  By Fact \ref{FactCompleteTensorProductFactCoefficientRingExtension}(i),  $\mathfrak{M}=\fm_R(R\otimes_CC')$ and it is a maximal ideal of $R\otimes_CC'$.  In view of our hypothesis that the transitions maps $\varphi_{i,j}$ are $C$-algebra homomorphisms, we are given with the transition maps $\text{id}_R\otimes \varphi_{i,j}:R\otimes_{C}V_i\rightarrow R\otimes_{C}V_j$ that are both $C$-algebra and $R$-algebra homomorphisms. Therefore, since tensor product commutes with direct limit we have $R$-algebra isomorphism 
		 \begin{equation}
		 \label{EquationDirectlimitPresentationResiduallyFinitelyGenratedEquicharacteristic}
		 R\otimes_CC'\cong \lim\limits_{\underset{i\in I}{\longrightarrow }}(R\otimes_CV_i).
		 \end{equation}  
		 
		(i) Since $\mathfrak{M}=\fm_R(R\otimes_CC')$ and it is a maximal ideal (Fact \ref{FactCompleteTensorProductFactCoefficientRingExtension}(i)), so from Fact \ref{ApplyingSharpResult} we deduce  that $\fm_R(R\otimes_CV_i)\in \text{Max}(R\otimes_CV_i)$ for each $i\in I$ by which we get $\big(\fm_R(R\otimes_CC')\big)\bigcap (R\otimes_CV_i)=\fm_R(R\otimes_CV_i)$ for each $i\in I$. This fact in conjunction with (\ref{EquationDirectlimitPresentationResiduallyFinitelyGenratedEquicharacteristic}) and  Fact \ref{DirectLimit} yields the direct limit presentation 
		\begin{equation}
		\label{EquationLocalizedDirectlimitPresentationResiduallyFinitelyGenratedEquicharacteristic}
		(R\otimes_CC')_{(\fm_R)}\cong \lim\limits_{\underset{i\in I}{\longrightarrow }}(R\otimes_CV_i)_{(\fm_R)}
		\end{equation}
		of $(R\otimes_CC')_{(\fm_R)}$ assigned to a direct system consisting of weakly unramified flat local homomorphisms. Namely, it is perhaps necessary to illustrate why $(\text{id}_R\otimes \varphi_{i,j})_{(\fm_R)}:(R\otimes_CV_i)_{(\fm_R)}\rightarrow (R\otimes_CV_j)_{(\fm_R)}$ is a flat homomorphism for each $i,j\in I$ with $i\le j$. But, indeed $R\rightarrow (R\otimes_CV_j)_{(\fm_R)}$ is a flat homomorphism as $C\rightarrow V_j$ is flat and localization is flat. Also, indeed $(R\otimes_RV_i)_{(\fm_R)}/(\fm_R)\rightarrow (R\otimes_RV_j)_{(\fm_R)}/(\fm_R)$ is  flat as it is a field extension. Consequently, our desired flatness follows from \cite[Tag 00MP]{Stacks}.   The rest of the proof of the equicharacteristic case will be given later after the following discussion on the mixed characteristic case of the lemma.

		  (ii)   $\mathfrak{M}=\fm_R(R\otimes_CC')$ and it is a maximal ideal (Fact \ref{FactCompleteTensorProductFactCoefficientRingExtension}(i)), so $\fm_R\big((R/pR)\otimes_{K_C}K_{C'})\in \text{Max}\big((R/pR)\otimes_{K_C}K_{C'})\big)$ (Fact \ref{ElementaryTensorProductFact}(ii)). Hence, from Fact \ref{ApplyingSharpResult} we deduce  that $$\fm_R\big((R/pR)\otimes_{K_C}K_{V_i}\big)\in \text{Max}\big((R/pR)\otimes_{K_C}K_{V_i}\big),\ \ \forall\ i\in I$$ i.e. $\fm_R(R\otimes_CV_i)\in \text{Max}(R\otimes_CV_i)$ for each $i\in I$. Consequently,  $\big(\fm_R(R\otimes_CC')\big)\bigcap (R\otimes_CV_i)=\fm_R(R\otimes_CV_i)$ for each $i\in I$. This fact in conjunction with (\ref{EquationDirectlimitPresentationResiduallyFinitelyGenratedEquicharacteristic}) and  Fact \ref{DirectLimit} yields the direct limit presentation 
		  \begin{equation}
		  \label{EquationLocalizedDirectlimitPresentationResiduallyFinitelyGenratedMixedCharacteristic}
		  (R\otimes_CC')_{(\fm_R)}\cong \lim\limits_{\underset{i\in I}{\longrightarrow }}(R\otimes_CV_i)_{(\fm_R)}
		  \end{equation}
		  of $(R\otimes_CC')_{(\fm_R)}$ assigned to a direct system consisting of weakly unramified flat local homomorphisms. It is perhaps necessary to illustrate why $(\text{id}_R\otimes_RV_i)_{(\fm_R)}:(R\otimes_CV_i)_{(\fm_R)}\rightarrow (R\otimes_CV_j)_{(\fm_R)}$ is flat for each $i,j\in I$ with $i\le j$. But this flatness, follows by the same argument as in the equicharacteristic in the previous paragraph, once we recall that any extension of such discrete valuation rings (in particular $C\rightarrow V_j$) is flat (Remark \ref{RemarkExtensionOfUnramifiedDVRs}).  Note that, for each $i\in I$ the field extension $K_C\rightarrow K_{V_i}$ is  finitely generated. Namely,   we can assume that $V_i\cong (C[\mathbf{X}]/\fa)_{\fq}$ for some $\fq \in \text{Spec}(C[\mathbf{X}]/\fa)$  ($V_i$ is essentially of finite type over $C$ by our hypothesis). Thus, $V_i/\fm_{V_i}\cong \text{Frac}(C[\mathbf{X}]/\fq)$ is a finitely generated field extension of $K_C$.
		  
		  Then the remainder of the proof works for both cases (i) and (ii).
		  
		   By our assumption each $V_i$ is essentially of finite type over $C$, so $(R\otimes_{C} V_i)_{(\fm_R)}$ is  essentially of finite type over $R$  for each $i\in I$. Consequently, $(R\otimes_{C} V_i)_{(\fm_R)}$ is Noetherian for each $i\in I$. Thus, the direct limit presentation in (\ref{EquationLocalizedDirectlimitPresentationResiduallyFinitelyGenratedEquicharacteristic}) and (\ref{EquationLocalizedDirectlimitPresentationResiduallyFinitelyGenratedMixedCharacteristic}) in conjunction with Theorem \ref{Ogoma} implies that $ (R\otimes_CC')_{(\fm_R)}$ is a Noetherian ring.  Therefore, in view of \cite[Theorem 3.5]{CelikbasWagstaff} as well as  Proposition \ref{CompleteTensorProductNoetherianness}(iii),  our desired test module property in the statement follows from the test module property of $T\otimes_R(R\otimes_CC')_{\mathfrak{M}}$ over $(R\otimes_CC')_{\mathfrak{M}}$. 
		   Thus, in view of  Lemma \ref{LemmaTestModuleDirectLimit} as well as (\ref{EquationLocalizedDirectlimitPresentationResiduallyFinitelyGenratedEquicharacteristic}) and (\ref{EquationLocalizedDirectlimitPresentationResiduallyFinitelyGenratedMixedCharacteristic}) it suffices to show that $T\otimes_R(R\otimes_CV_i)_{(\fm_R)}$ is a test module for $(R\otimes_CV_i)_{(\fm_R)}$ for each $i\in I$. 
		    But, from Lemma \ref{FinitelyGeneratedExension}   we deduce that $T\otimes_R(R\widehat{\otimes}_{C}V_i)$ is a test module over $R\widehat{\otimes}_{C}V_i$ (we recall that $R\widehat{\otimes}_{C}\widehat{V_i}\cong R\widehat{\otimes}_CV_i$). Consequently, $T\otimes_R(R\otimes_CV_i)_{(\fm_R)}$ is a test module over $(R\otimes_CV_i)_{(\fm_R)}$ by Remark \ref{RemarkTestModulePropertyDesentFromFaithfullyFlatExtension}, as $R\widehat{\otimes}_{C}V_i$ is $R$-algebra isomorphic to the $\big(\fm_R(R\otimes_{C}V_i)_{(\fm_R)}\big)$-adic completion of $(R\otimes_CV_i)_{(\fm_R)}$ by Proposition \ref{CompleteTensorProductNoetherianness}(iii) and Fact \ref{FactCompleteTensorProductFactCoefficientRingExtension}(i).
	\end{proof}
\end{sublem}

The following corollary is the main result of this subsection. Although  case (i) of the next corollary works for arbitrary field extension $C\rightarrow C'$, but case (ii) is stated only for  residually  purely transcendental   extension $C\rightarrow C'$.  However, in Lemma \ref{LemmaCoefficientPRingBaseChangeIsOKWhenSeparableResidueFieldExtension}, an analogue to case (ii)  is proved for arbitrary  $p$-ring extension $C'$ of $C$ with the only extra assumption that the residue field extension $K_C\rightarrow K_{C'}$ is separable. The proof of Lemma \ref{LemmaCoefficientPRingBaseChangeIsOKWhenSeparableResidueFieldExtension} uses some deep results in  commutative algebra.

\begin{subcor}\label{BasedFieldExtensionLemma} Let   $T$ be a  test module over a complete local ring $R$ with coefficient ring $\lambda_R:C\rightarrow R$.  Assume that either of the following conditions  holds.
	   \begin{enumerate}
	   	\item[(i)] $\varphi:C\rightarrow C'$ is a field extension (thus $R$ has equicharacteristic), or
	   	\item [(ii)] $C$ is a $p$-ring, $\{\mathbf{X}\}$ is an arbitrary set of indeterminates over $C$, $C':=C[\mathbf{X}]_{(p)}$ and finally $\varphi:C\rightarrow C'$ is the natural map given by $\varphi(c)=c$ (thus $R$ has mixed characteristic and $K_R$ has prime characteristic $p>0$).	   	
	   \end{enumerate}
   If, furthermore, either of the two conditions in the statement of    Corollary \ref{CorollaryTranscendentalExtension} holds, then  
   $T\otimes_R(R\widehat{\otimes}_{C}C')$   is a test module for $R\widehat{\otimes}_{C}C'$.
	\begin{proof}
		 In the  case (i), the statement  easily follows from Lemma \ref{LemmaCoefficientRingBaseChange}(i) because $C'=\bigcup\limits_{c_1,\ldots,c_n\in C'}C(c_1,\ldots,c_n)$.
		 
		In the case (ii), we have $C[\mathbf{X}]=\bigcup\limits_{X_1,\ldots,X_n\in  \{\mathbf{X}\}}C[X_1,\ldots,X_n]$ and it is clear that $$(pC[\mathbf{X}])\cap C[X_1,\ldots,X_n]=pC[X_1,\ldots,X_n],\ \ \forall\ X_1,\ldots,X_n\in \{\mathbf{X}\}.$$ So $C[\mathbf{X}]_{(p)}\cong \lim\limits_{\underset{X_1,\ldots,X_n\in \{\mathbf{X}\}}{\longrightarrow}}C[X_1,\ldots,X_n]_{(p)}$ (Fact \ref{DirectLimit}), and the statement  holds (Lemma \ref{LemmaCoefficientRingBaseChange}(ii)).
    \end{proof}
\end{subcor}

\subsection{The Celikbas and Sather-Wagstaff question, and, the Celikbas, Dao and Takahashi question}

As a first corollary to the previous subsection, we show that the module version of Question \ref{CWQuestion} has an affirmative answer in the category of local rings with uncountable residue field.
 The following corollary is the only result of our  paper where   the  nice  Sather-Wagstaff's result \cite[Theorem 4.8]{WagstaffAscent} plays a role.

\begin{subcor}\label{TestModuleAndFlatHomomorphism} Suppose that $\varphi:R\rightarrow S$ is a flat local homomorphism  such that $S/\fm_RS$ is regular, and $T$ is a test $R$-module. If either $K_R$ is uncountable or  $\text{CI-dim}_R(T)<\infty$, then $T\otimes_R S$ is a test module for $S$. 
	\begin{proof}
		Without loss of generality we can assume that $R$ and $S$ are both complete local rings. Namely, $\widehat{S}$ is flat over $S$, thus it is also flat over $R$ by the transitivity of  flatness. Thus, by applying \cite[Tag 00MP]{Stacks} to $R\rightarrow \widehat{R}\rightarrow \widehat{S}$ we deduce  that the natural map $\widehat{\varphi}:\widehat{R}\rightarrow \widehat{S}$ (i.e. given by $(r+\fm_R^n)_{n\in\mathbb{N}}\mapsto (\varphi(r)+\fm_S^n)_{n\in\mathbb{N}}$) is flat, implying that $\widehat{\varphi}$ is  weakly regular because $\widehat{S}/\fm_{\widehat{R}}\widehat{S}=\widehat{S}/\fm_R\widehat{S}\cong \widehat{S/\fm_RS}$ is regular. Moreover,    $T\otimes_R\widehat{R}$ is an $\widehat{R}$-test module  by virtue of \cite[Theorem 3.5]{CelikbasWagstaff}. If the statement of the corollary is proved for complete local rings, then  the test module property of $T\otimes_R\widehat{R}$ is preserved under the weakly regular flat local homomorphism $\widehat{R}\rightarrow \widehat{S}$. So, in view of the isomorphisms $$(T\otimes_{R}\widehat{R})\otimes_{\widehat{R}}\widehat{S}\cong T\otimes_R\widehat{S}\cong (T\otimes_RS)\otimes_S\widehat{S}$$ $(T\otimes_RS)\otimes_S\widehat{S}$ is an $\widehat{S}$-test module, implying that $T\otimes_RS$ is an $S$-test module by Remark \ref{RemarkTestModulePropertyDesentFromFaithfullyFlatExtension} (here we are also using the facts that $K_{\widehat{R}}\cong K_R$ and $\text{CI-dim}_R(T)=\text{CI-dim}_R(T\otimes_R\widehat{R})$, see \cite[(1.13) Proposition(2)]{Avramov}). 
		
		Now, assuming that $R$ and $S$ are complete local rings we divide the proof into the following two parts.
		
		\underline{Case 1, $R$ contains a field:} Let $\lambda_R:C\rightarrow R$ be a coefficient field of $R$. By Zorn's lemma, there exists a maximal subfield $C'$ of $S$ containing the field $\varphi\circ \lambda_R(C)$.  Taking into account the complete tensor product $R\widehat{\otimes}_CC'$ assigned to 
		the ring homomorphisms $\lambda_R:C\rightarrow R$ and $\varphi\circ \lambda_R:C\rightarrow C'$, we have the $R$-algebra homomorphism 
		\begin{equation}
	\label{Factorization}\varphi_{C'}:R\widehat{\otimes}_CC'=\lim\limits_{\underset{n\in \mathbb{N}}{\longleftarrow}}(R/\fm_R^n\otimes_CC')\rightarrow \lim\limits_{\underset{n\in\mathbb{N}}{\longleftarrow}}S/\fm_S^n=\widehat{S}\overset{\cong}{\longrightarrow}S
		\end{equation}
		where the left $R$-algebra homomorphism is induced by $$R/\fm_R^n\otimes_CC'\rightarrow S/\fm_S^n,\ \ \ (r+\fm_R^n)\otimes c'\mapsto \varphi(r)c'+\fm_S^n.$$ It is clear that then, $\varphi$ factors through $\eta_{R,R\widehat{\otimes}_CC'}:R\rightarrow R\widehat{\otimes}_CC'$ and $\varphi_{C'}:R\widehat{\otimes}_CC'\rightarrow S$, i.e. $\varphi=\varphi_{C'}\circ \eta_{R,R\widehat{\otimes}_CC'}$.
		
	We have $\mathfrak{M}^e=\fm_R(R\widehat{\otimes}_CC')\in \Max(R\widehat{\otimes}_CC')$ in view of parts (i) and (ii) of Fact   \ref{FactCompleteTensorProductFactCoefficientRingExtension}, and also by Fact \ref{FactTorOverBisTorOverAWhenBflatOverA} as well as Proposition \ref{FlatProposition}(ii) 
	\begin{align*}
	\numberthis
	\label{EquationTorVanishingImplyingFlatnessOfFactorization}
	\tor^{R\widehat{\otimes}_CC'}_{i}\big((R\widehat{\otimes}_CC')/\mathfrak{M}^e,S\rceil_{\varphi_{C'}}\big)
	&\cong
	\tor^{R\widehat{\otimes}_CC'}_{i}\big((R/\fm_R)\otimes_R(R\widehat{\otimes}_CC'),S\rceil_{\varphi_{C'}}\big)
	&\\&\cong
	\tor^{R}_{i}(R/\mam_R,S\rceil_{\varphi_{C'}\circ \eta_{R,R\widehat{\otimes}_CC'}})
	&\\&=
	\tor^{R}_{i}(R/\mam_R,S\rceil_{\varphi})
	&\\&=0,
	\ \ \forall\ i>0.
	\end{align*}
	 Thus, it follows from
	 Fact \ref{FactAndreCharmingResult} that $\varphi_{C'}:R\widehat{\otimes}_CC'\rightarrow S$ is also a flat local homomorphism ($(R\widehat{\otimes}_CC',\mathfrak{M}^e)$ is a complete local ring by Fact \ref{FactCompleteTensorProductFactCoefficientRingExtension}(ii)).  Moreover,  $S/\mathfrak{M}^eS=S/\fm_RS$ is a regular local ring by our hypothesis, thus $\varphi_{C'}$ is a weakly regular homomorphism.
		
		By Corollary \ref{BasedFieldExtensionLemma},  $T\otimes_R(R\widehat{\otimes}_CC')$ is a test module for  $R\widehat{\otimes}_CC'$.  Thus, if we can show that $\varphi_{C'}$ is a residually algebraic (local) ring homomorphism then  from  \cite[Theorem 4.8]{WagstaffAscent} we can conclude that $\big(T\otimes_R(R\widehat{\otimes}_CC')\big)\otimes_{R\widehat{\otimes}_CC'}S\cong T\otimes_RS$ is a test module for $S$, as was to be proved.
		But,  $C'\overset{\text{inclusion}}{\hookrightarrow}S\overset{\text{nat. epi.}}{\longrightarrow}K_S$ is an algebraic filed extension by Fact \ref{FactMaximalSubfieldAlgebraicExtension} that fits into the  commutative diagram of fields and field extensions:
		
		\begin{center}
			$\begin{CD} 
			C' @>\overline{\eta_{C',R\widehat{\otimes}_CC'}},\ \ \ c'\mapsto \big((1+\fm_R^n)\otimes c'\big)_{n\in \mn}+\mathfrak{M}^e>\cong > K_{R\widehat{\otimes}_CC'}\\
			@V \text{algebraic extension by Fact \ref{FactMaximalSubfieldAlgebraicExtension}}VV @VV\Psi V\\
			K_S @<<\Phi< K_{\widehat{S}}
			\end{CD}$
		\end{center}
		such that $\Phi\circ \Psi=\overline{\varphi_{C'}}$ (see (\ref{Factorization}), where $\varphi_{C'}$ is introduced). Note that the top row of the diagram is isomorphism in view of Fact \ref{FactCompleteTensorProductFactCoefficientRingExtension}(iv).
		 Consequently, $\overline{\varphi_{C'}}=\Phi\circ \Psi$ is am algebraic extension, as  required. 
		
		\underline{Case 2, $R$ has mixed characteristic:} Let $\lambda_R:C\rightarrow R$ be a coefficient ring  of $R$, thus $C$ is a $p$-ring where $p=\text{Char}(K_R)$. Let $\{\kappa_i\}_{i\in I}$ be a transcendental basis of $K_S$ over $\overline{\varphi}(K_R)=\overline{\varphi}\circ \overline{\lambda_R}(K_{C})$ and pick some lift $\{\zeta_{i}\}_{i\in I}\subseteq S$ of $\{\kappa_i\}_{i\in I}$ in $S$ (that is, for  every $i\in I$, $\zeta_i\in S$ has been chosen such that $\zeta_i+\fm_S=\kappa_i$). 
		
		 Let $\{X_i\}_{i\in I}$ be a set of indeterminates over $C$, so $C[\mathbf{X}]_{(p)}$ is an unramified discrete valuation ring of mixed characteristic $(0,p)$.  Let $C_\mathcal{L}$ be the $(p)$-adic completion of $C[\mathbf{X}]_{(p)}$.
		 Then $T\otimes_R(R\widehat{\otimes}_C(C[\mathbf{X}]_{(p)})$ is a test module over $R\widehat{\otimes}_C(C[\mathbf{X}]_{(p)})$ in view of Corollary \ref{BasedFieldExtensionLemma}. 
		  Equivalently, $T\otimes_R(R\widehat{\otimes}_CC_{\mathcal{L}})$ is a test module over $R\widehat{\otimes}_CC_{\mathcal{L}}$ (see Remark \ref{RemarkDefinitionCompleteTensorProductDoesntAffectCompleteness}). 
		  
		  The natural map $\gamma:C[\mathbf{X}]\overset{c\mapsto \varphi\circ \lambda_R(c),\ X_i\mapsto \zeta_i}{\longrightarrow} S$   induces the ring homomorphism $\gamma_{(p)}:C[\mathbf{X}]_{(p)}\rightarrow S$. Namely, it suffices to show that $\gamma(f)\in S\backslash \fm_S$ for each $f\in C[\mathbf{X}]\backslash (p)$ (then the universal property of the localization shows that $\gamma_{(p)}$ is induced by $\gamma$). Let $$f=\sum\limits_{\underline{i}:=i_1,\ldots,i_t\in I,\ \underline{n}:=n_1,\ldots,n_t\in \mathbb{N}_0^t}c_{\underline{i},\underline{n}}X_{i_1}^{n_1}\ldots X_{i_t}^{n_t}\in C[\mathbf{X}]\backslash (p)$$ where  this formal sum is taken over finitely many distinct monomial in $\mathbf{X}$. So $$\gamma(f)=\sum\limits_{\underline{i}:=i_1,\ldots,i_t\in I,\ \underline{n}:=n_1,\ldots,n_t\in \mathbb{N}_0^t}\varphi\circ\lambda_R{(c_{\underline{i},\underline{n}})}\zeta_{i_1}^{n_1}\ldots \zeta_{i_t}^{n_t}.$$ If $\gamma(f)\in \fm_S$ then we get 
		  $$\overline{\gamma(f)}=\sum\limits_{\underline{i}:=i_1,\ldots,i_t\in I,\ \underline{n}:=n_1,\ldots,n_t\in \mathbb{N}_0^t}\overline{\varphi}\circ\overline{\lambda_R}{(c_{\underline{i},\underline{n}}+pC)}\kappa_{i_1}^{n_1}\ldots \kappa_{i_t}^{n_t}=0\in K_S$$ implying that $c_{\underline{i},\underline{n}}\in pC$ for each $\underline{i},\underline{n}$  (because $\{\kappa_i\}_{i\in I}$ is transcendental over $\overline{\varphi}\circ \overline{\lambda_R}(K_C)$). But this contradicts with $f\notin pC[\mathbf{X}]$, consequently we have $\gamma(f)\notin \fm_S$ as  required. Now, from $\gamma_{(p)}:C[\mathbf{X}]_{(p)}\rightarrow S$ we obtain the ring homomorphism $$\widehat{\gamma_{(p)}}:C_{\mathcal{L}}\rightarrow \widehat{S},\ \ \  (f_n/g_n+p^nC[\mathbf{X}]_{(p)})_{n\in \mathbb{N}}\mapsto (\gamma(f_n)\gamma(g_n)^{-1}+\fm_S^n)_{n\in \mathbb{N}}.$$
		  
		  Let $\eta_{S,\widehat{S}}^{-1}$ be the  isomorphism $\widehat{S}\rightarrow S$ that is the inverse to $S\overset{s\mapsto (s+\fm_S^n)_{n\in \mathbb{N}}}{\longrightarrow} \widehat{S}$ and set $$\Gamma:=\eta_{S,\widehat{S}}^{-1}\circ \widehat{\gamma_{(p)}}:C_{\mathcal{L}}\rightarrow S.$$ 
		  
		  Then  $\varphi\circ \lambda_R:C\overset{}{\rightarrow}S$ factors ($\varphi\circ \lambda_R=\Gamma\circ e$)  $$C\overset{e:\ c\mapsto (c/1+p^nC[\mathbf{X}]_{(p)})_{n\in \mathbb{N}}}{\longrightarrow} C_{\mathcal{L}}\overset{\Gamma}{\longrightarrow } S.$$ 
		  
		  So, considering the complete tensor product $R\widehat{\otimes}_CC_{\mathcal{L}}$ assigned to the ring homomorphisms $\lambda_R:C\rightarrow R$ and $e:C\rightarrow C_{\mathcal{L}}$   we have the $R$-algebra local homomorphism $$\varphi_{C_{\mathcal{L}}}:R\widehat{\otimes}_CC_{\mathcal{L}}=\lim\limits_{\underset{n\in\mathbb{N}}{\longleftarrow}}\big((R/\fm_R^n)\otimes_C(C_\mathcal{L}/p^nC_\mathcal{L})\big)\rightarrow \lim\limits_{\underset{n\in\mathbb{N}}{\longleftarrow}}(S/\fm_S^n)=\widehat{S}\overset{\cong}{\rightarrow} S$$ where the first homomorphism is induced by the family of homomorphisms $$(R/\fm_R^n)\otimes_C(C_\mathcal{L}/p^nC_\mathcal{L})\overset{}{\rightarrow} S/\fm_S^n,\ \ \ (r+\fm_R^n)\otimes (c'+(p)^n)\mapsto \varphi(r)\Gamma(c')+\fm_S^n.$$ Thus,
		  we  can factor $R\rightarrow S$ as $$R\overset{\eta_{R,R\widehat{\otimes}_CC_{\mathcal{L}}}}{\longrightarrow} R\widehat{\otimes}_{C}C_{\mathcal{L}}\overset{\varphi_{C_{\mathcal{L}}}}{\longrightarrow} S$$ while $T\otimes_R(R\widehat{\otimes}_CC_{\mathcal{L}})$ is a test module over $R\widehat{\otimes}_CC_{\mathcal{L}}$ as mentioned in the second paragraph of the mixed characteristic case of the proof. Hence if we can show that the local $R$-algebra homomorphism $\varphi_{C_{\mathcal{L}}}$ is weakly regular and residually algebraic, then from  \cite[Theorem 4.8]{WagstaffAscent} we can deduce that $$\big(T\otimes_R(R\widehat{\otimes}_CC_{\mathcal{L}})\big)\otimes_{R\widehat{\otimes}_CC_{\mathcal{L}}}S \cong T\otimes_RS$$ is an $S$-test module as was to be proved.
		  
		  But exactly by the same argument as in the equicharacteristic case of the proof, the analogue to the Tor vanishing in (\ref{EquationTorVanishingImplyingFlatnessOfFactorization}) holds, i.e. $\tor^{R\widehat{\otimes}_CC_{\mathcal{L}}}_i\big((R\widehat{\otimes}_CC_{\mathcal{L}})/\mathfrak{M}^e,S\rceil_{\varphi_{C_\mathcal{L}}}\big)=0$ for $i>0$ and this implies that $\varphi_{C_\mathcal{L}}$ is flat by 
		  Fact \ref{FactAndreCharmingResult}. Also, again $S/\mathfrak{M}^eS=S/\fm_RS$ (parts (i) and (ii) of Fact \ref{FactCompleteTensorProductFactCoefficientRingExtension}) is regular by our hypothesis, thus $\varphi_{C_{\mathcal{L}}}$ is weakly regular.
		  
		   It remains to show that $\varphi_{\mathcal{L}}$ is residually algebraic. To this aim, we first show that $\Gamma$ is  residually algebraic. From $\Gamma\big((c/1+(p)^n)_{n\in \mathbb{N}}\big)=\varphi\circ\lambda_R(c)$ for each $c\in C$, we deduce that $\overline{\varphi}\circ \overline{\lambda_R}(K_C)\subseteq \overline{\Gamma}(K_{C_\mathcal{L}})$. From $\Gamma\big((X_i/1+(p)^n)_{n\in \mathbb{N}}\big)=\zeta_i$ for each $i\in I$, we conclude that $\kappa_i\in \overline{\Gamma}(K_{C_\mathcal{L}})$ for each $i\in I$. It follows that $\overline{\Gamma}(K_{C_\mathcal{L}})$ contains the subfield $\big(\overline{\varphi}\circ \overline{\lambda_R}(K_C)\big)(\kappa_i:i\in I)$ of $K_S$. Consequently, $K_S$ is algebraic over $\overline{\Gamma}(K_{C_\mathcal{L}})$, because $\{\kappa_i\}_{i\in I}$ is a transcendental basis of $K_S$ over $\overline{\varphi}\circ \overline{\lambda}_R(K_C)$. Now, we consider the commutative diagram 
		   	\begin{center}
		   	$\begin{CD} 
		   	K_{C_{\mathcal{L}}}@>\overline{\eta_{C_{\mathcal{L}},R\widehat{\otimes}_CC_{\mathcal{L}}}},\ \ \ c'+pC_{\mathcal{L}}\mapsto \big((1+\fm_R^n)\otimes (c'+p^nC_{\mathcal{L}})\big)_{n\in \mathbb{N}}+\mathfrak{M}^e>\cong > K_{R\widehat{\otimes}_CC_{\mathcal{L}}}\\
		   	@V \overline{\Gamma}V\text{algebraic extension} V @VV\Psi V\\
		   	K_S @<<\Phi< K_{\widehat{S}}
		   	\end{CD}$
		   \end{center}
		    where $\Phi\circ \Psi=\overline{\varphi_{C_\mathcal{L}}}$ (see where $\varphi_{C_{\mathcal{L}}}$ is introduced). Note that the top row of the diagram is isomorphism in view of Fact \ref{FactCompleteTensorProductFactCoefficientRingExtension}(iv). It follows that $\overline{\varphi_{C_\mathcal{L}}}(=\Phi\circ \Psi)$ is an  algebraic field extension and the proof is complete.
		 	\end{proof}
\end{subcor}

Now, to answer Question \ref{CWQuestion} for rings with uncountable residue field we show that the complex case of Corollary \ref{TestModuleAndFlatHomomorphism} holds.

\begin{subcor}\label{TestComplexAndFlatHomomorphism}
 Suppose that $\varphi:R\rightarrow S$ is a flat local homomorphism  with regular closed fiber  $S/\fm_RS$, and $T$ is a test complex for $R$. If  $K_R$ is uncountable or $\text{CI-dim}_R(T)<\infty$, then $T\otimes_R S$ is a test complex for $S$.
	\begin{proof}
	
		First note that since $S$ is flat over $R$, so $T\otimes_R^{\mathbf{L}}S=T\otimes_RS$. Suppose that  $X$  is a homologically finite complex of $S$-modules with $X\otimes_S^{\mathbf{L}}(T\otimes_R S)\in \mathcal{D}_b(S)$. Let $$P^X:=\cdots\rightarrow  P^X_i\rightarrow P^X_{i-1}\rightarrow P^X_{i-2}\rightarrow \cdots\rightarrow P^X_t\rightarrow 0,$$ be a     projective resolution for $X$ consisting of finitely generated free $S$-modules (such a projective resolution exists in view of part (L) of \cite[(A.3.2) Theorem (Existence of Resolutions)]{ChrisensenGorenstein}), and $$P^T:=\cdots\rightarrow P^T_i\rightarrow P^T_{i-1}\rightarrow \cdots\rightarrow P^T_s\rightarrow 0$$ be a projective resolution of $T$ by finite free $R$-modules. As $X$ is homologically finite, so there is a hard truncation  $$P^X_{\ge m}:=\cdots\rightarrow  P^X_{m+1}\rightarrow P^X_{m}\rightarrow0$$ of $P$  which is an acyclic complex, where $m \ge \text{Sup}\ X+1$ is a natural number ($\text{Sup}\ X$ denotes the largest integer $i$ with $\text{H}_{i}(X)\neq 0$). Set, $C^X:=\text{H}_m(P^X_{\ge m})$,  thus some shift of $P^X_{\ge m}$   forms a projective resolution of $C^X$ and $C^X$ is a finitely generated $S$-module (because $P^X$   is consisting of finitely generated free $S$-modules). In the same vein,  there is a finitely generated $R$-module $C^T$ by taking cokernel of the tail of some acyclic hard truncation  $P^T_{\ge m}$ of  $P^T$ (we can pick $m$ so that it works for both of $P^T$ and $P^X$). 
		
		Note that $C^T$ is a test module for $R$, because $\tor^R_\gg(M,C^T)=0$  if and only if $M\otimes_RP^T_{\ge m}\in \mathcal{D}_b(R)$ and this is the case where $M\otimes_RP^T=M\otimes_R^{\mathbf{L}}T\in \mathcal{D}_b(R)$, for each finitely generated $R$-module $M$.
		 We have  $\tor^S_{\gg }(C^T\otimes_RS,C^X)=0$. Namely,  tensoring the exact sequence $$0\rightarrow P^T_{\ge m}\rightarrow P^T\rightarrow P^T_{<m}\rightarrow 0$$ of complexes to $S$ we get the  exact sequence of complexes of $S$-modules
		 	$$ \ 0\rightarrow P^{T}_{\ge m}\otimes_RS\rightarrow P^T\otimes_RS\rightarrow  P^T_{<m}\otimes_RS \rightarrow 0$$
 where $P^T_{<m}:=0\rightarrow P^T_{m-1}\rightarrow \cdots\rightarrow P^T_s\rightarrow 0$.

		 Then tensoring the resulted  exact sequence of  complexes to $X$, we get the   exact sequence
		$$ \ 0\rightarrow (P^{T}_{\ge m}\otimes_RS)\otimes_S X\rightarrow (P^T\otimes_RS)\otimes_S X\rightarrow (P^T_{<m}\otimes_RS)\otimes_S X\rightarrow 0.$$ This obtained exact sequence shows that  
        $(P^{T}_{\ge m}\otimes_RS)\otimes_SX\in\mathcal{D}_b(S)$,
		 because  $(P^T\otimes_R S)\otimes_S X\in \mathcal{D}_b(S)$ by our hypothesis while $$(P^T_{<m}\otimes_RS)\otimes_S X\overset{\text{by definition}}{\simeq}(P^T_{<m}\otimes_RS)\otimes^{\mathbf{L}}_S X\in \mathcal{D}_b(S)$$ by \cite[(A.5.6) FD Corollary]{ChrisensenGorenstein} (note that $(P^T_{<m}\otimes_RS)$ is a bounded complex of projetive $S$-modules). Consequently   		
		 \begin{align*}
		 \numberthis
		 \label{EquationHomologiacllyFinite}
		 (P^{T}_{\ge m}\otimes_RS)\otimes_S X\simeq (P^{T}_{\ge m}\otimes_RS)\otimes_S P^X\in \mathcal{D}_b(S)
		 \end{align*}
		    by \cite[(A.4.1) Presentation of Quasi-Isomorphisms and Equivalences]{ChrisensenGorenstein}. 
		 Setting $$P^X_{<m}:=0\rightarrow P^X_{m-1}\rightarrow\cdots\rightarrow P^X_t\rightarrow 0$$  and discussing similarly as above we can conclude that there is an exact sequence 
		 $$0\rightarrow (P^T_{\ge m}\otimes_R S)\otimes_S P^X_{\ge m} \rightarrow   (P^T_{\ge m}\otimes_R S)\otimes_S P^X\rightarrow  (P^T_{\ge m}\otimes_R S)\otimes_S P^X_{<m}\rightarrow 0.$$
		 This exact sequence and (\ref{EquationHomologiacllyFinite}) imply that $(P^{T}_{\ge m}\otimes_RS)\otimes_S P^{X}_{\ge m}\in \mathcal{D}_b(S)$, i.e.  $$\tor^S_{\gg}(C^T\otimes_RS,C^X)=0.$$
		By Corollary \ref{TestModuleAndFlatHomomorphism} $C^T\otimes_RS$ is a test module for $S$ (note that  by \cite[Proposition 3.7]{WagstaffComplete}, $\text{CI-dim}_R\ C^T<\infty $  provided $\text{CI-dim}_R\ T<\infty$). So the above Tor vanishing yields $\pd_S(C^X)<\infty$, thus $\pd_S(X)<\infty$.
	\end{proof}
\end{subcor}

\begin{subrem} \emph{In the statement of   Corollary \ref{TestComplexAndFlatHomomorphism}, the assumption on the regularity of the closed fiber $S/\fm_RS$  is necessary,  see \cite[Example 3.6]{CelikbasWagstaff}.}
\end{subrem}

As the  main result of this section,  we answer Question \ref{CDTQuestion} affirmatively.

\begin{subcor}\label{MainResult} Suppose that $R$ is local ring  possessing a test module  $T$  such that $\text{CI-dim}_R(T)<\infty$. Then $R$ is complete intersection.
	\begin{proof}
	 If $\text{CI-dim}_R(T)=-\infty$, i.e. $T=0$, then $R$ is regular and we are done. So we suppose that $T$ has finite complete intersection dimension.
	 
	 There is a quasi-deformation $R\overset{g}{\rightarrow} R'\twoheadleftarrow A$ with $\pd_A(T\otimes_R R')<\infty$. 
	  Furthermore, by virtue of \cite[Theorem F]{Wagstaff}  we can,  without loss of generality, assume  that the closed fiber $R'/\mathfrak{m}_RR'$ is a  Gorenstein Artinian ring. Afterwards, applying Remark \ref{QuasiDeformationWithCompleteDeformation} we can assume that $R'$ and $A$ are all complete local rings (while the Artinian Gorenstein property of $R'/\fm_RR'$ is preserved). Since $R'$ is a complete local ring, so we have $g=\mu_{R'}^{-1}\circ \widehat{g}\circ \mu_R$  where $\mu_R:R\rightarrow \widehat{R}$ and $\mu_{R'}:R'\rightarrow \widehat{R'}$ are   canonical map to the completion and $\widehat{g}:\widehat{R}\rightarrow \widehat{R'}$ is given by the rule $(r_n+\fm_R^n)_{n\in\mathbb{N}}\mapsto (g(r_n)+\fm_{R'}^n)_{n\in\mathbb{N}}$. Thus, applying   \cite[Tag 00MP]{Stacks} to the local homomorphisms $R\overset{\mu_R}{\rightarrow} \widehat{R}\overset{\mu_{R'}^{-1}\circ \widehat{g}}{\rightarrow} R'$ we conclude that $\mu_{R'}^{-1}\circ \widehat{g}:\widehat{R}\rightarrow R'$ is a flat local homomorphism. So $\widehat{R}\overset{\mu_{R'}^{-1}\circ \widehat{g}}{\longrightarrow}R'\twoheadleftarrow A$ is a quasi-deformation of $\widehat{R}$ while $\pd_A\big((T\otimes_R\widehat{R})\otimes_{\widehat{R}}R'\big)=\pd_A(T\otimes_RR')<\infty$, implying that $\text{CI-dim}_{\widehat{R}}(T\otimes_R\widehat{R})<\infty$. Moreover, $T\otimes_R\widehat{R}$  is an $\widehat{R}$-test module by virtue of \cite[Theorem 3.5]{CelikbasWagstaff}. If  the statement is true for complete local rings then we can conclude that $\widehat{R}$, equivalently $R$, is a complete intersection and the statement follows in general. It follows that, without loss of generality, we can presume that $R$, $R'$ and $A$ are complete local rings while $R'/\fm_RR'$ is Artinian and Gorenstein.
	  
	  Set $T':=T\otimes_R R'$ and $R':=A/\mathbf{x}A$ where  $\mathbf{x}$ is a regular sequence of $A$.   In the  light of   Definition and Remark \ref{DefinitionRemarkCohenFactorization},
	  there is a Cohen factorization $$R\overset{\delta}{\rightarrow} B\overset{\rho}{\twoheadrightarrow} R'$$ of $g:R\rightarrow R'$, i.e. $g=\rho\circ \delta$, $\delta$ is a weakly regular homomorphism of complete local rings and $\rho$  is a surjective homomorphism.
	  Using Corollary \ref{TestModuleAndFlatHomomorphism}, we know that $T\otimes_R B$ is a test module for $B$. 
	   
	   By \cite[Corollary 3.9]{CelikbasWagstaff}  (cf. \cite[Corollary 3.4]{Celikbas}), we already know that $R$ is Gorenstein. Hence, since the closed fiber $R'/\fm_RR'$ is Gorenstein as discussed in the first paragraph of the proof and since  Gorensteinness deforms  so  $R'$ and $A$ are both Gorenstein rings (see \cite[Corollary 3.3.15]{BrunsHerzogCohenMacaulay} and \cite[Exercise 18.1, page 152]{Matsumura}).
	   
	    Let $\Omega$ be a sufficiently high syzygy of the residue field $K_A$ of $A$ (over $A$), so that  $\mathbf{x}\subseteq A$  is also an $\Omega$-regular sequence in view of Fact \ref{SyzygyAndGrade}(ii).
       Then, $$\tor^R_{\gg}(T,\Omega/\mathbf{x}\Omega\rceil_{g})\overset{\text{Fact \ref{FactTorOverBisTorOverAWhenBflatOverA}}}{\cong}\tor^{R'}_{\gg}(T',\Omega/\mathbf{x}\Omega)\overset{\text{Fact \ref{FactTorOverRTorOverRxR}}}{=}\tor^A_{\gg}(T',\Omega)\overset{\pd_A(T')<\infty}{=}0.$$  

		 Thus  $$ \tor^B_{\gg}(T\otimes_RB ,\Omega/\mathbf{x}\Omega\rceil_{\rho})\overset{\text{Fact \ref{FactTorOverBisTorOverAWhenBflatOverA}}}{=} \tor^R_{\gg}(T,\Omega/\mathbf{x}\Omega\rceil_{\rho\circ \delta})=\tor^R_{\gg}(T,\Omega/\mathbf{x}\Omega\rceil_{g})=0$$ implying that $\pd_B(\Omega/\mathbf{x}\Omega\rceil_{\rho})<\infty$  in view of the test module property of the $B$-module $T\otimes_RB$ (note that $\Omega/\mathbf{x}\Omega\rceil_{\rho}$ is  a finitely generated $B$-module as $\rho$ is surjective). 
		 
		 Consequently, 
		 \begin{align*}
		   \tor^A_{\gg}(K_R\otimes_RR',\Omega)
		   \overset{\text{Fact \ref{FactTorOverRTorOverRxR}}}{\cong}
		   \tor^{R'}_{\gg}(K_R\otimes_RR',\Omega/\mathbf{x}\Omega)
		   \overset{\text{Fact \ref{FactTorOverBisTorOverAWhenBflatOverA}}}&{\cong}
		   \tor^{R}_{\gg}(K_R,\Omega/\mathbf{x}\Omega\rceil_{\underset{=g}{\underbrace{\rho\circ \delta}}})
		   \\\overset{\text{Fact \ref{FactTorOverBisTorOverAWhenBflatOverA}}}&{\cong}
		   \tor^B_{\gg}(K_R\otimes_RB,\Omega/\mathbf{x}\Omega\rceil_{\rho})
		   \\\overset{}&{=}0\ \ (\pd_B(\Omega/\mathbf{x}\Omega\rceil_{\rho})<\infty).
		 \end{align*}  
		    But $\Omega$ is a syzygy of $K_A$, therefore $\tor^A_{\gg}(K_R\otimes_RR',K_A)=0$ in view of the above display. It follows that, $\pd_A(K_R\otimes_RR')<\infty$ and that $\text{CI-dim}_R(K_R)<\infty$. Now the statement follows from \cite[(1.3) Theorem]{Avramov}.
		 	\end{proof}
	 \end{subcor}

\section{Alternative proofs for special cases}
 In this section,   we give  alternative proofs for  Corollary \ref{TestModuleAndFlatHomomorphism}  that is not  using the  recent  Sather-Wagstaff's result \cite[Theorem 4.8]{WagstaffAscent}, provided  either $R$ is Cohen-Macaulay and equicharacteristic (Corollary \ref{CorollaryAlternativeProofCohenMacaulayCase}), or $\overline{\varphi}:K_R\rightarrow K_S$ is  a separable   field extension (Corollary \ref{CorollaryAlternativeProofCWQuestionSeparableCase}).  The results of this section also show that the equicharacteristic case of Corollary \ref{MainResult}, which answers positively Question \ref{CDTQuestion}, can be deduced without applying \cite[Theorem 4.8]{WagstaffAscent}, because Question \ref{CDTQuestion} is a question in the category of Cohen-Macaulay rings (see   Corollary \ref{CorollaryAlternativeProofCDTQuestion}). 
 
 \subsection{Preliminaries}

We first recall the complete version of  Nakayama's lemma.

\begin{subfact}\label{FactCompleteNakayamsLemma} (\cite[Exercise 7.2, page 203]{EisenbudCommutative}) 
	Suppose that $R$ is an $\fa$-adically complete ring for some ideal $\fa$ of $R$ and that $M$ is an $R$-module with $\bigcap\limits_{n\in \mathbb{N}}\fa^nM=0$   (we do not assume that $M$ is finitely generated in advance). If $M/\fa M$ is generated by  the image of $m_1,\ldots,m_n\in M$  then $M$ is generated by $m_1,\ldots,m_n$. 
\end{subfact}

The next elementary fact will be used in the proof of Proposition \ref{ArtinianFactorization}. Note that an Artinian equicharacteristic ring  $R$ is complete, thus it admits a coefficient field $\lambda_R:C\rightarrow R$. Moreover,  in the statement of the next fact   we may, and we do, assume that $\lambda_R$ is inclusion and $C$ is a subfield of $R$ (by replacing $C$ with $\lambda_R(C)$).

\begin{subfact} \label{FactArtinianIsPolynomialQuotient} Suppose that $R$ is an equicharacteristic Artinian local ring and that $C\subseteq R$ is a coefficient field of $R$. Then $(R,\fm_R)\cong \big(C[\mathbf{X}]/\fa,(\mathbf{X})\big)$ for some sequence of indeterminates   $\mathbf{X}:=X_1,\ldots,X_n$  over $C$.
	\begin{proof}
		Assume that $\fm_R$ is generated by $r_1,\ldots,r_n\in R$ and consider $\gamma:C[\mathbf{X}]\overset{c\mapsto c,\ X_i\mapsto r_i}{\longrightarrow} R$ where $\mathbf{X}:=X_1,\ldots,X_n$ is a sequence of indeterminates over $R$. Since $R$ is Artinian, so there is some $h\in \mathbb{N}$ such that $X_i^h\in \text{Ker}(\gamma)$ for each $1\le i\le n$. Consequently, $\gamma$ induces the ring homomorphism $\overline{\gamma}:C[\mathbf{X}]/(X_1^h,\ldots,X_n^h)\rightarrow R$ of Artinian rings. But $C[\mathbf{X}]/(X_1^h,\ldots,X_n^h)$ is an $(\mathbf{X})$-adically complete  (Artinian) local ring while $R/(X_1,\ldots,X_n)R=K_R$ is a cyclic $C[\mathbf{X}]/(X_1^h,\ldots,X_n^h)$-module  generated by $1_{K_R}$ (because $C$ is a coefficient field of $R$ and $C\overset{\text{incl.}}{\rightarrow} R\twoheadrightarrow K_R$ is an isomorphism). So $R$ is generated by $1_R$ over $C[\mathbf{X}]/(X_1^h,\ldots,X_n^h)$ by Fact \ref{FactCompleteNakayamsLemma}. This implies that  $\overline{\gamma}$ is surjective and the statement follows.
	\end{proof}
\end{subfact}

\begin{subfact}\label{FactArtinianFiniteDimensiona}
	Suppose that $K$ is a field and $f:K\rightarrow A$ is a ring homomorphism where $A$ is an Artinian local ring. If $K\rightarrow A\overset{\text{nat. epi.}}{\twoheadrightarrow}K_A$ is a finite field extension, then $A$ is a finite dimensional $K$-vector space.
	\begin{proof}
		Let $\pi_n:A\rightarrow A/\fm_A^n$ be the natural homomorphism for each $n\in \mathbb{N}$. By our hypothesis $A/\fm_A$ is a $K$-finite dimensional vector space via $\pi_1\circ f$. Suppose,  inductively, that $A/\fm_A^n$ is proved to be a finitely generated $K$-module via $\pi_n\circ f$. 
		As $A$ is Noetherian and $\fm_A$ is finitely generated so $\fm_A^n/\fm_{A}^{n+1}$ is a finitely generated $(A/\fm_A^n)$-module, hence it is a finite dimensional $K$-vector space via $\pi_n\circ f$.
		Consequently, the exact sequence $0\rightarrow \fm_A^n/\fm_A^{n+1}\rightarrow A/\fm_A^{n+1}\rightarrow A/\fm_A^n\rightarrow 0$ of $K$-modules (where the $A/\fm_A^{n+1}$ is a $K$-module via $\pi_{n+1}\circ f$) shows that $A/\fm_A^{n+1}$ is finitely generated over $K$. Since $A$ is Artinian we have $A/\fm_A^n=A$ and $\pi_n=\id_A$ for sufficiently large $n$.
	\end{proof}
\end{subfact}

\begin{subdefi}\label{RemarkSeparable}
	\emph{(\cite[page 198]{Matsumura}) Let $\varphi:K\rightarrow \mathcal{K}$ be a field extension. Then we say that $\mathcal{K}$ is \textit{separable} over   $K$, or $\varphi$ is a \textit{separable extension}, provided the tensor product ring $\mathcal{K}\otimes_K L$ is reduced for every field extension $L$ of $K$.}
\end{subdefi}


Sometimes in the literature, the notion of  $I$-smoothness    defined below  is  referred as formally smoothness  at $I$ (see e.g. \cite[Definition 2.2.1]{MajadasRodicioSmoothness}, or \cite[Tag 07EB]{Stacks} in conjunction with \cite[Tag 07EC(2)]{Stacks} and \cite[Tag 07NI]{Stacks}). 

\begin{subdefi}
  \emph{(\cite[page 213 and page 214]{Matsumura}) Let $A$ be a ring, $B$ an $A$-algebra and $I$ be an ideal of $B$. We say that $B$ is \textit{$I$-smooth}     over $A$ if given an $A$-algebra $C$, 	an ideal $N$ of $C$ satisfying  $N^2=0$  and an $A$-algebra homomorphism $u:B\rightarrow C/N$ which is continuous for the discrete topology of $C/N$ (that is, such that $u(I^v)=0$ some $v$), then there exists a   lifting  $\nu:B\rightarrow C$ of $u$ to $C$:}	 
  	  $$\begin{tikzcd}
  	      B \arrow[r,"u"] \arrow [dr,dotted,"\nu"] &C/N\\
  	      A \arrow[u]\arrow [r] &C\arrow[u,two heads,"nat.\ epi." swap]
  	  \end{tikzcd}$$ 
\end{subdefi}

\begin{subfact}\label{FactISmoothCompletion}\cite[Lemma 2.2.3]{MajadasRodicioSmoothness}
	Let $A\rightarrow B$ be a ring homomorphism, $J$ an ideal of $B$, $C$ an $A$-algebra and $T$ an ideal of $C$. Assume that $B$ is a $J$-smooth $A$-algebra and $C$ complete for the $T$-adic topology. 	
	 Let $f_1:B\rightarrow C/T$ be a continuous $A$-algebra homomorphism (i.e.  $f_1(J^n)=0$ for some $n$). Then  there exists a continuous  $A$-algebra homomorphism $f:B\rightarrow C$ inducing $f_1$ (i.e. $f_1=(C\overset{\text{nat. epi.}}{\twoheadrightarrow}C/T)\circ f$). 
\end{subfact}

We will refer to \cite{SpivakovskyANewProof} for (another proof of) the Popescu's result that a regular homomorphism of Noetherian rings is a filtered inductive limit  of smooth algebras of finite type, so for the notion of a  smooth algebra (in the sense used in \cite{SpivakovskyANewProof}) we refer the following definition. In view of \cite[page 382, first line]{SpivakovskyANewProof}, the author of \cite{SpivakovskyANewProof} has used the term ``smooth'' to describe the following stronger concept of $I$-smoothness.  


\begin{subdefi} \label{DefinitionSmooth}\emph{
		Let $B$ an $A$-algebra. Suppose that $B$ is $I$-smooth over $A$ for any ideal $I$ of $B$. Then we say that  \textit{$B$ is a smooth $A$-algebra in the sense of Spivakovsky} (cf. \cite[(28.D) DEFINITION]{MatsumuraCommutativeAlgebra}, see also \cite[(28.C) DEFINITION]{MatsumuraCommutativeAlgebra}). Equivalently, we say that  \textit{$B$ is a formally smooth $A$-algebra in the sense of stacks project} (see \cite[Tag 00TI and Tag  07EC(1)]{Stacks}).}
\end{subdefi}

 

\begin{subdefi}\label{DefinitionRegularHomomorphism}
	\emph{(see \cite[Definition, page 249]{MatsumuraCommutativeAlgebra}, or \cite[Tag 07BZ]{Stacks}) Let $A$ be a Noetherian ring containing a field $k$. We say that $A$ is \textit{geometrically regular } over $k$ if, for any finite extension $k'$ of $k$, the ring $A\otimes_kk'$ is regular (this equivalent to saying that $A\otimes_kk'$ is regular for every finitely generated field extension $k'$ of $k$, see \cite[Tag 0381]{Stacks}). We say that a homomorphism $\varphi:A\rightarrow B$ is \textit{regular} (or that $B$ is a regular $A$-algebra) if it is flat and for each $\fp\in \text{Spec}(A)$ the fiber $B\otimes_Ak(\fp)$ is geometrically regular over $k(\fp)$ (this is equivalent to saying that $B\otimes_A\mathcal{L}$ is regular for every field $\mathcal{L}$ that is an essentially of finite type $B$-algebra).}
\end{subdefi}

\begin{subdefi}\label{DefinitionLocallyCompleteIntersection}
	\emph{(\cite[Definition, page 458]{AvramovLocally}) A local homomorphism $\varphi:R\rightarrow S$ is said to be  \textit{complete intersection} at $\fm_S$ if in some Cohen factorization (see Definition and Remark \ref{DefinitionRemarkCohenFactorization}) $\Ker \varphi'$ is generated by a regular sequence (this property does not depend on the choice of Cohen factorization). A homomorphism 	of Noetherian rings $\varphi:R\rightarrow S$ is \textit{complete intersection at a prime ideal $\fq$ of $S$} if the induced local homomorphism $\varphi_{\fq}:R_{\fq\cap R}\rightarrow S_\fq$ is complete intersection at $\fq S_{\fq}$. A homomorphism that has this property for all $\fq\in \text{Spec}(S)$ is said to be \textit{locally complete intersection}.}
\end{subdefi}

For the definition  of the Andr\'e-Quillen homology $D_n(S|R;N)$, assigned to a ring homomorphism $\varphi:R\rightarrow S$ of commutative rings and an $S$-module $N$, we refer to \cite[5.8 Andr\'e-Quillen homology and chomology]{IyengarAndreQuillen}. In the sequel, we will use  Andr\'e-Quillen homology  (only) to   characterize    locally complete intersection homomorphisms and regular homomorphisms (as cited in the following fact).

\begin{subfact}\label{FactSmoothHomomorphismCollection}
  Let $A$ be a local ring and  $\varphi:A\rightarrow B$ 	 be a homomorphism of Noetherian rings.
  \begin{enumerate}
  	\item [(i)]  If $\varphi$  is essentially of finite type and formally smooth in the sense of stacks project (or equivalently, smooth in the sense of Spivakovsky), then $\varphi$ is
  	 regular. 
  	\item[(ii) ] Suppose that  $\varphi:A\rightarrow B$ is of finite type, and  formally smooth  in the sense of stacks project. Assume that $\fq\in \text{Spec}(B)$ contains  $\fm_AB$.  Then the induced local map $\varphi_{\fq}:A\rightarrow B_\fq$ is regular.
  	\item[(iii)]  If $A,B$ are complete local rings, $\varphi$ is a local homomorphism and $B$ is an $\fm_B$-smooth  $A$-algebra, then $\varphi$ is a regular homomorphism.
  	\item[(iv)]  Suppose that $\varphi$ is a regular homomorphism. Then $B$ is a filtered inductive limit of smooth $A$-algebras of finite type (smooth in the sense of Spivakovsky). Equivalently,  any commutative diagram 
  	\begin{center}\begin{tikzcd}
  		A \arrow[r,"\varphi"] \arrow[d,"\tau" swap] &B\\
  		C\arrow[ur,"\rho" swap]
  	\end{tikzcd}\end{center}
  where $C$ is a finitely generated $A$-aglebra, can be extended to a commutative diagram 
  \begin{center}\begin{tikzcd}
  	A \arrow[r,"\varphi"] \arrow[d,"\tau" swap] &B\\
  	C\arrow[ur,"\rho" swap] \arrow[r,"\phi" swap] &D \arrow[u,"\psi" swap]
  \end{tikzcd}\end{center}
  where $D$ is a smooth finitely generated $A$-algebra (smooth in the sense of Spivakovsky).
  	\item[(v)] $\varphi$ is regular if and only if $D_n(B|A;-)=0$ for all $n\ge 1$, or equivalently precisely when $D_1(B|A;-)=0$.
  	\item[(vi)] $\varphi$ is locally complete intersection if and only if $D_2(B|A;-)=0$.
  	  \end{enumerate}
  	\begin{proof}
  		(i) \cite[Corollary 2.9]{SpivakovskyANewProof}.
  		
  		
  		(ii) In view of  part (i),  it suffices to show that $A\rightarrow B_{\fq}$ is also formally smooth. So let's assume that for some ideal $IB_{\fq}$ of $B_{\fq}$ the  solid commutative diagram 
  		
  		\begin{center}\begin{tikzcd}
  			B \ar[r,"\mu"] \ar[drrrr,dotted,"v" swap] &B_\fq \ar[rrr,"u"] \ar[rrrd,dotted,"v_\fq"]&&&C/N \\
  			A \ar[r,equal]  \ar[u,"\varphi"]& A \ar[rrr] \ar[u,"\varphi_{\fq}"] &&&C \ar[u,"\pi:=\text{nat. epi.}" swap,twoheadrightarrow]
  		\end{tikzcd}\end{center}
  	
  	is given such that $u(I^nB_{\fq})=0$ for some $n\in \mathbb{N}$, the ideal $N$ of $C$ has square zero and $\mu$ is the localization map. Since $B$ is formally smooth over $A$, so there is a  homomorphism $v:B\rightarrow C$ such that  $\pi\circ v=u\circ \mu$. Let $s\in B\backslash \fq$.	 Then, $\pi\circ v(s)$ is invertible in $C/N$ because $\mu(s)$ is so in $B_\fq$. As $N$ is a nilpotent ideal, this implies that $v(s)$ is invertible in $C$. It follows that $v$ induces  $v_\fq:B_\fq\rightarrow C$ ($b/s\rightarrow v(b)v(s)^{-1}$). It is easily verified that then $\pi\circ v_\fq=u$, as required.
  		

  		(iii) \cite[Tag 07PM]{Stacks}.
  		
  		(iv) See \cite[Theorem 1.1]{SpivakovskyANewProof}, as well as the third paragraph of \cite[page 382]{SpivakovskyANewProof}.
  	
  		(v) \cite[Theorem 9.5]{IyengarAndreQuillen}.
  		
  		(vi) \cite[(1.2) Second Vanishing Theorem]{AvramovLocally}.
  	\end{proof}
\end{subfact}

\begin{subfact} \label{FactReductionToWeaklyUnramified}
	Suppose that $\varphi:R\rightarrow S$ is a weakly regular homomorphism of  local rings and that $T$ is  a test module for $R$. Let $\mathbf{x}$ be a sequence of elements of $S$ which forms a regular system of parameters for the regular local ring $S/\fm_RS$. Set $S':=S/\mathbf{x}S$.
	\begin{enumerate}
		\item [(i)] $R\overset{\varphi}{\rightarrow}S\overset{\text{nat. epi.}}{\twoheadrightarrow}S'$ is a flat weakly unramified local homomorphism.
		\item[(ii)] If $T\otimes_RS'$ is a test module for $S'$, then $T\otimes_RS$ is a test module over $S$.
	\end{enumerate}
  \begin{proof}
   The flatness mentioned in the statement of part (i) follows immediately from \cite[Corollary, page 177]{Matsumura}, while  \cite[Corollary, page 177]{Matsumura} also implies that $\mathbf{x}$ is a regular sequence on $S$.  
   
   (i) Since $\mathbf{x}$ forms a regular system of parameters for $S/\fm_RS$, so $S'/\fm_RS'\cong S/\big(\fm_RS+(\mathbf{x})\big)\cong K_S$ implying that $\fm_RS'=\fm_{S'}$, i.e. $S'$ is    weakly unramified over $R$.   Hence part (i) holds.

  	(ii) We assume that $T\otimes_R S'$ is a test module over $S'$ and we prove that $T\otimes_RS$ is a test module over $S$. Presume that $N$ is a finitely generated $S$-module with $\tor^S_{\gg}(N,T\otimes_R S)=0$. Replacing $N$ with its sufficiently high syzygy, we may assume that $\mathbf{x}$ is also a regular sequence on $N$ (see Fact \ref{SyzygyAndGrade}(ii)).  Set $N':=N/\mathbf{x}N$.   
  	From our hypothesis and Fact \ref{FactTorOfMAndTorOfSpecializationM} (by setting $a=1$ in its statement) we get $\tor^S_{\gg}(T\otimes_RS,N')=0$. Thus 
  	$$\tor^{S'}_{\gg}(T\otimes_RS',N')\overset{\text{Fact  \ref{FactTorOverBisTorOverAWhenBflatOverA}}}{\cong} \tor^{R}_{\gg}(T,N')\overset{\text{Fact  \ref{FactTorOverBisTorOverAWhenBflatOverA}}}{\cong}\tor^{S}_{\gg}(T\otimes_RS,N')=0.$$
  So  $\pd_{S'}(N')<\infty$ by our hypothesis. But $\pd_S(N)=\pd_{S'}(N')$ because tensoring to $S'$ a minimal projective resolution of $N$ over $S$ yields a minimal projective resolution  of $N'$ over $S'$ (by Fact \ref{ResolutionModuleRegularSequence}). This completes the proof.
  \end{proof}	
\end{subfact}

\begin{subfact}\label{FactKoszulComplexVanishingHomology} Suppose that $B$ is a local ring and that $M$ is a complex of $R$-modules with finitely generated homologies. Let $\mathbf{x}=x_1,\ldots,x_n$ be a sequence of elements of $\fm_B$ ($n\in\mathbb{N}$).
	\begin{enumerate}
 	\item[(i)] $H_i\big(K_\bullet(\mathbf{x};B)\otimes_BM\big)$ is a finitely generated $B$-module for each $i\in \mathbb{Z}$.
 	\item[(ii)] Let  $i\in \mathbb{Z}$.  If  $H_i\big(K_\bullet(\mathbf{x};B)\otimes_BM\big)=0$ then $H_i(M)=0$. 
	\end{enumerate}
   \begin{proof}
     Suppose that $n=1$. From \cite[Proposition 1.6.12(b)]{BrunsHerzogCohenMacaulay}, we get the exact sequence $H_i(M)\overset{\pm x_1}{\rightarrow} H_i(M)\rightarrow H_i\big(K_\bullet(x_1;B)\otimes_BM\big)\rightarrow H_{i-1}(M)$ for each $i\in I$,  which implies the statement of part (i).  It also implies the statement of part (ii), because   by Nakayama's lemma $H_i(M)=0$ if and only if the multiplicative map $H_i(M)\overset{\pm x_1}{\rightarrow}H_i(M)$ is surjective.   Now, assume that $n>1$ and the statements of part (i) and (ii) hold for  $n-1$. Set $\mathbf{x}':=x_1,\ldots,x_{n-1}$. Thus, the complex  $K_{\bullet}(\mathbf{x}';B)\otimes_BM$ has finitely generated homologies by our inductive hypothesis, so the base of the induction implies that      
      $K_\bullet(\mathbf{x};B)\otimes_BM\cong K_\bullet(x_n;B)\otimes_B\big(K_\bullet(\mathbf{x}';B)\otimes_BM)\big)$ also  has finitely generated homologies. This shows that part (i) holds, and for part (ii) we have $H_i\Big(K_\bullet(x_n;B)\otimes_B\big(K_\bullet(\mathbf{x}';B)\otimes_BM\big)\Big)\cong H_i\big(K_\bullet(\mathbf{x};B)\otimes_BM\big)=0$ by our hypothesis, therefore $H_i\big(K_\bullet(\mathbf{x}';B)\otimes_BM\big)=0$ by the base of the induction. Consequently, $H_i(M)=0$ by induction hypothesis.
   \end{proof}
\end{subfact}

\subsection{An alternative proof for the residually separable case}

The main idea  in our alternative proof of Corollary \ref{TestModuleAndFlatHomomorphism},   in the case where $\overline{\varphi}:K_R\rightarrow K_S$ is separable, is to apply the idea used in the proof of Linquan Ma's \cite[Lemma 5.1]{MaLech}: 

\begin{subprop}\label{LemmaNiceFactorization}
	Suppose that $\varphi:R\rightarrow S$ is a  local homomorphism of complete local  rings with Artinian closed fiber $S/\fm_RS$. Assume that the following coefficient rings containment holds: 
	
	there are coefficient rings $\lambda_R:C_R\rightarrow R,\ \lambda_S:C_S\rightarrow S$ and a  flat local  homomorphism $\sigma:C_R\rightarrow C_S$ such that  $\varphi\circ \lambda_R=\lambda_S\circ \sigma$. 
	
	Consider the complete tensor product $R\widehat{\otimes}_{C_R}C_S$  assigned to  $\lambda_R$ and $\sigma$. Then  $\varphi:R\rightarrow S$ factors into  two local homomorphisms  $$R\underset{\text{flat}}{\overset{\eta_{R,R\widehat{\otimes}_{C_R}C_S}}{\longrightarrow}} R\widehat{\otimes}_{C_R}C_S \underset{\text{}}{\overset{\varphi_{_{C_S}}}{\longrightarrow}} S$$ of complete local rings such that $\varphi_{_{C_S}}$ is a module-finite  homomorphism  (i.e. $\varphi=\varphi_{_{C_S}}\circ \eta_{R,R\widehat{\otimes}_{C_R}C_S}$ and $S\rceil_{\varphi_{_{C_S}}}$ is a finitely generated  $(R\widehat{\otimes}_{C_R}C_S)$-module). 
	
	 Moreover, \begin{enumerate}
	 	\item [(i)] if $\varphi$ is flat then $\varphi_{_{C_S}}$ is also flat.
	 	\item[(ii)] if $\varphi$ is flat and weakly unramified then $\varphi_{_{C_S}}$ is an isomorphism.
	 \end{enumerate}
	  \begin{proof}
	  	  Since $\sigma$ is flat, so  $C_R$ and $C_S$   have the same characteristic (without the flatness assumption $C_R$ can be a $p$-ring while $C_S$ is a field of prime characteristic).  So $\sigma$ is either a field extension or an extension of $p$-rings. 
	  	  Note that by Fact \ref{FactCompleteTensorProductFactCoefficientRingExtension}(ii), $(R\widehat{\otimes}_{C_R}C_S,\mathfrak{M}^e)$ is a complete local ring and $\eta_{R,R\widehat{\otimes}_{C_R}C_S}$ is a  flat weakly unramified local homomorphism. By Fact \ref{FactCompleteTensorProductFactCoefficientRingExtension}(i) $\mathfrak{M}\in \text{Max}(R\otimes_{C_R}C_S)$, so $R\widehat{\otimes}_{C_R}C_S$ is  $(R\otimes_{C_R}C_S)$-algebra isomorphic to the $\mathfrak{M}$-adic completion  
	  	  of $R\otimes_{C_R}C_S$ in view of Proposition \ref{CompleteTensorProductNoetherianness}(i).  Thus, we continue the proof of the existence of $\varphi_{_{C_R}}$ with the $\mathfrak{M}$-adic completion of $R\otimes_{C_R}C_S$ in place of $R\widehat{\otimes}_{C_R}C_S$ (this is to avoid considering different cases based on the characteristic).   In view of the equality $\varphi\circ \lambda_R=\lambda_S\circ \sigma$, there exists the natural  ring homomorphism $$R\otimes_{C_R}C_S\overset{}{\rightarrow} S,\ \ r\otimes c\mapsto \varphi(r)\lambda_S(c)$$ taking $\mathfrak{M}$ into $\fm_S$. So we get the  $(R\otimes_{C_R}C_S)$-algebra homomorphism  $$\lim\limits_{\underset{n\in \mathbb{N}}{\longleftarrow}}\big((R\otimes_{C_R}C_S)/\mathfrak{M}^n\big)\rightarrow \widehat{S}\overset{\cong}{\rightarrow}S$$
	  	 which  sends $(r\otimes 1+\mathfrak{M}^n)_{n\in \mathbb{N}}$ to $\varphi(r)$. This completes the proof of the existence of $\varphi_{_{C_S}}$ as well as the equality $\varphi_{_{C_S}}\circ \eta_{R,R\widehat{\otimes}_{C_R}C_S}=\varphi$.  If $\varphi$ is flat, then the flatness of $\varphi_{_{C_S}}$ follows from \cite[Tag 00MP]{Stacks} because  $(R\widehat{\otimes}_{C_R}C_S)/\fm_R(R\widehat{\otimes}_{C_R}C_S)$ is a field ($\eta_{R,R\widehat{\otimes}_{C_R}C_S}$ is weakly unramified). To show that $\varphi_{_{C_S}}$ is a module-finite homomorphism, we argue as follows.
	  	 
	  	 $S/\mathfrak{M}^eS=S/\fm_RS$ is an Artinian local ring by our hypothesis. Thus, if we can show that the field extension $$\overline{\varphi_{_{C_S}}}:(R\widehat{\otimes}_{C_R}C_S)/\mathfrak{M}^e\rightarrow K_S$$ is  an isomorphism, then by applying Fact  \ref{FactArtinianFiniteDimensiona} to $K_{R\widehat{\otimes}_{C_R}C_S}\rightarrow S/\mathfrak{M}^eS$ we can deduce that  $S/\mathfrak{M}^eS$ is finitely generated over $R\widehat{\otimes}_{C_R}C_S/\mathfrak{M}^e$.  This together with Fact \ref{FactCompleteNakayamsLemma} would imply that $S$ is a finitely generated $(R\widehat{\otimes}_{C_R}C_S)$-module. Hence  we should show that  $\overline{\varphi_{_{C_S}}}$ is an isomorphism.  
	  	 Let $\pi:S\rightarrow K_S$ be the natural epimorphism. It is easy to see that, $\varphi_{_{C_S}}\circ \eta_{C_S,R\widehat{\otimes}_{C_R}C_S}=\lambda_{S}$ by recalling  how we  constructed $\varphi_{_{C_S}}$ in the first paragraph of the proof. Consequently, $\pi\circ\varphi_{_{C_S}}\circ \eta_{C_S,R\widehat{\otimes}_{C_R}C_S}=\pi\circ\lambda_{S}$ implying that $\pi\circ\varphi_{_{C_S}}\circ \eta_{C_S,R\widehat{\otimes}_{C_R}C_S}$ is surjective as $\lambda_S$ is assumed to be a coefficient ring for $S$. From this we conclude that $\overline{\varphi_{_{C_S}}}$ is onto (i.e. an isomorphism), because $\pi\circ\varphi_{_{C_S}}\circ \eta_{C_S,R\widehat{\otimes}_{C_R}C_S}=\overline{\varphi_{_{C_S}}}\circ \overline{\eta_{C_S,R\widehat{\otimes}_{C_R}C_S}}\circ (C_S\overset{\text{nat. epi.}}{\twoheadrightarrow}K_{C_S})$.

	  	 Finally, under the condition (ii) ($\varphi$ being weakly unramified) we have $S/\mathfrak{M}^eS=S/\fm_RS=S/\fm_S$. Moreover, as we have seen in the previous paragraph $\overline{\varphi}:(R\widehat{\otimes}_{C_R}C_S)/\mathfrak{M}^e\rightarrow K_S$ is surjective.  These two facts imply that $R\widehat{\otimes}_{C_R}C_S/\mathfrak{M}^e\rightarrow S/\mathfrak{M}^eS=K_S$ is surjective. It follows that $S/\mathfrak{M}^eS$ is generated by $1_{S/\mathfrak{M}^eS}$ over $R\widehat{\otimes}_{C_R}C_S$, and consequently $S$ is generated by $1_S$ over $R\widehat{\otimes}_{C_R}C_S$ by virtue of Fact \ref{FactCompleteNakayamsLemma}. So $\varphi_{_{C_S}}$ is surjective, while it is also faithfully flat and injective by (i).
	  \end{proof}
\end{subprop}	

 The condition on the containment of coefficient fields in  Proposition \ref{LemmaNiceFactorization} holds if $\mathbb{Q}\subseteq R$ (as mentioned in the proof of \cite[Lemma 5.1]{MaLech}). But, in general it   holds    when the residue field extension $\overline{\varphi}:K_R\rightarrow K_S$ is separable:

\begin{subprop}\label{PropositionResiduallySeparableImpliesHomomorphismOfCoefficientRingsExist}
	Suppose that $\varphi:R\rightarrow S$ is a homomorphism of complete local rings such that the induced residue field extension $\overline{\varphi}:K_R\rightarrow K_S$ is separable. Presume that  either $R$ and $S$ are both equicharacteristic, or they both have mixed characteristic (this is the case where $\varphi$ is flat, or $R$ is equicharacteristic). Let $\lambda_R:C_R\rightarrow R$ be a coefficient ring for $R$. Then there is a  coefficient ring $\lambda_S:C_S\rightarrow S$ for  $S$   and a  (flat  local) ring homomorphism $\sigma:C_R\rightarrow C_S$ such that $\varphi\circ \lambda_R=\lambda_S\circ \sigma$.
	\begin{proof}

		\underline{Case 1:}   If $R$  (equivalently, $S$) has equicharacteristic then $C_R$ is a field. The field extension $$C_R\overset{\overset{\varphi\circ\lambda_R}{\cong}}{\rightarrow}(\varphi\circ\lambda_R)(C_R)\overset{\text{inclusion}}{\rightarrow}S\overset{\text{nat. epi.}}{\twoheadrightarrow}K_S$$ is  separable, because it coincides with the extension $$C_R\overset{\lambda_R}{\rightarrow}R\overset{\text{nat. epi.}}{\twoheadrightarrow}K_R\overset{\overline{\varphi}}{\rightarrow}K_S$$ while $C_R\overset{\lambda_R}{\rightarrow}R\overset{\text{nat. epi.}}{\twoheadrightarrow}K_R$ is an isomorphism and $\overline{\varphi}$ is separable by our assumption. Consequently, $(\varphi\circ\lambda_R)(C_R)\overset{\text{inclusion}}{\rightarrow}S\overset{\text{nat. epi.}}{\twoheadrightarrow}K_S$ is  also  a separable field extension and then the statement  holds by virtue of parts (iii) and (iv) of \cite[Theorem 28.3]{Matsumura}.
		
		\underline{Case 2:} Suppose that $R$ and $S$ have mixed characteristic  and $\text{Char}(K_R)=p>0$. Although it is possible to present an independent proof for this case, but in order to reduce the complexity,  we prove this case   by passing to the equicharacteristic case. 
		
		By our hypothesis, the induced ring homomorphism $\widetilde{\lambda_R}:C_R/pC_R\rightarrow R/pR$, by $\lambda_R$, is a coefficient field for $R/pR$. 
		Since $R/pR$ and $S/pS$ are complete local rings of prime characteristic (thus equicharacteristic), so in in view of  case 1 there is  a  coefficient field $\lambda_{S/pS}:C_{S/pS}\rightarrow S/pS$ and a field extension $\tau:C_{R}/pC_R\rightarrow C_{S/pS}$ such that $\widetilde{\varphi}\circ \widetilde{\lambda_R}=\lambda_{S/pS}\circ \tau$ where    $\widetilde{\varphi}:R/pR\rightarrow S/pS$ is the induced ring homomorphism by $\varphi$. In particular we have the commutative diagram,
		\begin{center}
			$\begin{CD}
			  C_{R}/pC_R @>\tau>> C_{S/pS}\\
			  @V(R/pR\overset{\text{nat. epi.}}{\twoheadrightarrow}K_R)\circ \widetilde{\lambda_{R}}V\cong V 
			  @V\cong V(S/pS\overset{\text{nat. epi.}}{\twoheadrightarrow}K_S)\circ \lambda_{S/pS}V\\
			  K_R @>>\overline{\varphi}> K_S
			\end{CD}$ 
		\end{center}
	which shows that the field extension $\tau$ is also a separable field extension (as so is $\overline{\varphi}$ and the vertical maps in the commutative diagram are isomorphisms). We can consider a $p$-ring $C_S$ with an isomorphism $\Psi:C_{S/pS}\rightarrow C_S/pC_S$ (by \cite[Theorem 29.1]{Matsumura} with $\mathbb{Z}_{p\mathbb{Z}}$ in place $A$ in its statement). Then applying \cite[Theorem 29.2]{Matsumura} to the $p$-ring $C_R$ and the filed extension $C_R/pC_R\overset{\tau}{\rightarrow}C_{S/pS}\overset{\Psi}{\rightarrow}C_S/pC_S$, we can deduce that there is local homomorphism $\sigma:C_R\rightarrow C_S$ such that $\overline{\sigma}:C_R/pC_R\rightarrow C_S/pC_S$ coincides with $\Psi\circ \tau$. By Remark \ref{RemarkExtensionOfUnramifiedDVRs}, $\sigma$ is  flat and local. Since $\tau$ is a separable extension and $\Psi$ is an isomorphism so $\Psi\circ \tau$ is also a separable field extension, thus it is $(0)$-smooth by \cite[Theorem 26.9]{Matsumura}. Consequently, $\sigma:C_R\rightarrow C_S$ (which induces $\Psi\circ \tau$ on residue fields) is $(p)$-smooth by \cite[Theorem 28.10]{Matsumura}. It is easily seen that the  diagram 
	\begin{center}\begin{tikzcd}
	   C_R \arrow[rr,"\lambda_R"] \arrow[d,"\sigma" swap] &&R \arrow[rr,"\varphi"] &&S \arrow[d,"\text{nat.\ epi.}",two heads]\\
	   C_S \arrow[r,two heads,"\text{nat. epi.}" swap] \ar[rrrru,dotted,"\lambda_S"]
	   &C_S/pC_S \arrow[r,"\Psi^{-1}" swap] &C_{S/pS}\arrow[r,"\lambda_{S/pS}" swap] &S/pS\arrow [r,two heads,"\text{nat. epi.}"] &K_S
	\end{tikzcd}\end{center}
	with solid arrows is commutative. Thus since $\sigma$ is $(p)$-smooth  and the ring homoorphism $C_S\rightarrow K_S$ in the bottom row sends $p$ to zero, so in view of Fact \ref{FactISmoothCompletion} there is an $C_R$-algebra homomorphism $\lambda_S:C_S\rightarrow S$ which makes the completed diagram commutative. 
	 Hence, $\lambda_S\circ \sigma=\varphi\circ \lambda_R$  (the commutativity of the upper triangle) and  $\lambda_S:C_S\rightarrow S$ is a coefficient ring for $S$ (by  the commutativity of lower  triangle as well as the surjectivity of $\Psi^{-1}$ in conjunction with the fact that $C_{S/pS}$ is a coefficient field for $S/pS$). So the proof is complete.
	\end{proof}
\end{subprop}

We deduce the following corollary from the previous proposition, which will be used in the proof of Lemma \ref{LemmaCoefficientPRingBaseChangeIsOKWhenSeparableResidueFieldExtension}.

\begin{subcor}\label{CorollaryMixedCharacteristicResiduallyFinitelyGeneratedSeparableWeaklyRegularHomomorphism}
Suppose that $R$ is a complete local ring of mixed characteristic. Let $\varphi:R\rightarrow S$ be a weakly regular local homomorphism of local rings such that $\overline{\varphi}:K_R\rightarrow K_S$ is   separable and finitely generated. Let $T$ be a  test module over $R$. If $K_R$ is uncountable or $\text{CI-dim}_R(T)<\infty$,  then $T\otimes_RS$ is a test module over $S$.
  \begin{proof}
  	 Since the natural map $S\rightarrow \widehat{S}$ induces isomorphism on residue fields and $\widehat{S}/\fm_R\widehat{S}\cong \widehat{S/\fm_RS}$ is regular, so without loss of generality we can assume that $S$ is complete (in view of Remark \ref{RemarkTestModulePropertyDesentFromFaithfullyFlatExtension}). Moreover, by Fact \ref{FactReductionToWeaklyUnramified}, without loss of generality, we can assume that $\varphi$ is weakly unramified as well. 
  	 
  	  Let   $\lambda_R:C_R\rightarrow R$ be a coefficient ring for $R$. By Proposition \ref{PropositionResiduallySeparableImpliesHomomorphismOfCoefficientRingsExist}, there is a coefficient ring $\lambda_S:C_S\rightarrow S$, of $S$, and a (flat local) extension $\sigma:C_R\rightarrow C_S$ of $p$-rings such that $\lambda_S\circ \sigma=\varphi \circ \lambda_R$. Thus by Proposition \ref{LemmaNiceFactorization}(ii), without loss of generality, we may assume that $S=R\widehat{\otimes}_{C_R}C_S$ and $\varphi=\eta_{R,R\widehat{\otimes}_{C_R}C_S}$. So the statement follows from Lemma \ref{FinitelyGeneratedExension} (the residue field extension $\overline{\sigma}:C_R/pC_R\rightarrow C_S/pC_S$ is also finitely generated as required in the statement of Lemma \ref{FinitelyGeneratedExension}, because   $\overline{\varphi}$ is a finitely generated field extension  while  $\overline{\lambda_R}:C_R/pC_R\rightarrow K_R,\ \overline{\lambda_{S}}:C_S/pC_S\rightarrow K_S$ are both isomorphisms and $\overline{\lambda_S}\circ \overline{\sigma}=\overline{\varphi}\circ\overline{\lambda_R}$).
  \end{proof}
\end{subcor}

In order to prove Corollary \ref{CorollaryAlternativeProofCWQuestionSeparableCase}, we need the following lemma which is the mixed characteristic version of Corollary \ref{BasedFieldExtensionLemma}(i) given with extra assumption that  the residue field extension is separable.

\begin{sublem}\label{LemmaCoefficientPRingBaseChangeIsOKWhenSeparableResidueFieldExtension}
	Suppose that $R$  is a complete local ring of mixed characteristic and that $\lambda_R:C_R\rightarrow R$ is a coefficient ring for $R$. Let $\sigma:C_R\rightarrow C'$ be an  extension of $p$-rings such that $\overline{\sigma}:K_{C_R}\rightarrow K_{C'}$ is a separable field extension. 
	 Let $T$ be a test module for $R$. If $\text{CI-dim}_R(T)<\infty$ or $K_R$ is uncountable, then $T\otimes_R(R\widehat{\otimes}_{C_R}C')$ is a test module for $R\widehat{\otimes}_{C_R}C'$.
	\begin{proof}
	The residue field extension $\overline{\eta_{R,R\widehat{\otimes}_{C_R}C'}}:K_R\rightarrow K_{R\widehat{\otimes}_{C_R}C'}$, which is a separable field extension in view of Fact \ref{FactCompleteTensorProductFactCoefficientRingExtension}(iii), is $0$-smooth by \cite[Theorem 26.9]{Matsumura}.	By Fact \ref{FactCompleteTensorProductFactCoefficientRingExtension}(ii) 
	$\eta_{R,R\widehat{\otimes}_{C_R}C'}$ is a weakly unramified homomorphism. Consequently, $\eta_{R,R\widehat{\otimes}_{C_R}C'}$ is $\mathfrak{M}^e$-smooth in view of \cite[Theorem 28.10]{Matsumura}. It follows that then $\eta_{R,R\widehat{\otimes}_{C_R}C'}$  is a regular homomorphism (Fact \ref{FactSmoothHomomorphismCollection}(iii)). Consequently,  by  Fact \ref{FactSmoothHomomorphismCollection}(iv) $R\widehat{\otimes}_{C_R}C'$
	is a filtered inductive limit of smooth $R$-algebras  of finite type
	  (smooth in the sense of Spivakovsky, or formally smooth in the sense of stacks project). Using this fact,  similarly as in the case of Lemma \ref{LemmaTestModuleDirectLimit}, we will   deduce the statement. However, the   transition maps in the statement of Lemma \ref{LemmaTestModuleDirectLimit}  are assumed to be flat local homomorphisms, that is not the case here necessarily (note that, moreover, here we have a filtered inductive limit rather than a direct limit). To remedy this, we argue as follows. 
	   
	
	  
	  Let $N$ be a finitely generated $(R\widehat{\otimes}_{C_R}C')$-module such that 
	  \begin{equation}
	  \label{EquationLastTorVanishing?}
	  \tor^{R\widehat{\otimes}_{C_R}C'}_{\gg}\big(T\otimes_R(R\widehat{\otimes}_{C_R}C'),N\big)=0.
	  \end{equation}
	  Let $\mathcal{F}_{\bullet}:=0\rightarrow F_2\overset{[a_{l,k}]}{\rightarrow}F_1\overset{[a'_{l,k}]}\rightarrow F_0\rightarrow 0$ be a truncation of a minimal free resolution of $N$ (thus $H_1(\mathcal{F}_{\bullet})=0$ and $H_0(\mathcal{F}_{\bullet})=N$). Let $C$ be the finitely generated $R$-subalgebra  of  $R\widehat{\otimes}_{C_R}C'$ generated, over $R$, by all entries $a_{l,k}$ and $a'_{l,k}$ of the matrices of the differentials of $\mathcal{F}_\bullet$. We denote by $\tau:R\rightarrow C$ the natrual map $r\mapsto \eta_{R,R\widehat{\otimes}_{C_R}C'}(r)\in C$, and by $\rho:C\rightarrow R\widehat{\otimes}_{C_R}C'$ the inclusion map. Then by virtue of Fact \ref{FactSmoothHomomorphismCollection}(iv), there is a commutative diagram  with solid rows
	  \begin{center}\begin{tikzcd}
	  	  	R \arrow[d,"\tau" swap] \arrow[rrr,"\eta_{R,R\widehat{\otimes}_{C_R}C'}"] &&& R\widehat{\otimes}_{C_R}C'\\
	  	  	C\arrow[rrr,"\phi" swap] \arrow[rrru,"\rho"] &&& S \arrow[u,"\psi"] \arrow[r,"\mu" swap,dotted] & S_\fq \arrow[lu,"\psi_{\fq}" swap,dotted]
	  \end{tikzcd}\end{center}
	  of ring homomorphisms where $S$ is a  formally smooth (in the sense of stacks project) finitely generated $R$-algebra.  Setting  $\fq=\psi^{-1}(\mathfrak{M}^e)$, we get the dotted arrows in the above commutative diagram where $\mu$ stands for the localization map and $\psi_{\fq}$ is induced  by $\psi$ ($\psi_{\fq}(s/t)=\psi(s)\psi(t)^{-1}$). 	  
	  We can construct a length $2$ (perfect) complex  $_{S_{\fq}}\mathcal{F}_{\bullet}$ of finite free $S_{\fq}$-modules whose matrices of differentials are $[\mu\circ \phi(\underset{\in C}{\underbrace{a'_{l,k}}})]$ and $[\mu\circ \phi(\underset{\in C}{\underbrace{a_{l,k}}})]$, so  it is straightforward to check that 
	  \begin{equation}
	  \label{EquationAurora}
	  _{S_{\fq}}\mathcal{F}_{\bullet}\otimes_{S_{\fq}}(R\widehat{\otimes}_{C_R}C')\cong \mathcal{F}_{\bullet}.
	  \end{equation}



 The $R$-algebra homomorphism $\psi_{\fq}:S_{\fq}\rightarrow R\widehat{\otimes}_{C_R}C'$ yields the exact sequence of Andr\'e-Quillen homologies 
\begin{equation}
\label{EquationAndreQuillen}
D_2\big(R\widehat{\otimes}_{C_R}C'|R;M\big)\rightarrow D_2\big(R\widehat{\otimes}_{C_R}C'|S_{\fq};M\big)\rightarrow D_1\big(S_{\fq}|R;M\big)
\end{equation}
 for each $(R\widehat{\otimes}_{C_R}C')$-module $M$ (see \cite[6.7 Transivity]{IyengarAndreQuillen}).
 
  By Fact \ref{FactSmoothHomomorphismCollection}(ii), $S_{\fq}$ is a regular $R$-algebra. So, from Fact \ref{FactSmoothHomomorphismCollection}(v) we get $D_1(S_\fq|R;-)=0$.  
   Moreover, $D_2\Big(R\widehat{\otimes}_{C_R}C'|R;-\big)=0$ again by Fact \ref{FactSmoothHomomorphismCollection}(v). It turns out that $\psi_{\fq}:S_{\fq}\rightarrow R\widehat{\otimes}_{C_R}C'$ is a locally complete intersection homomorphism in view of (\ref{EquationAndreQuillen})  and by virtue of Fact \ref{FactSmoothHomomorphismCollection}(vi). Hence  there is some  Cohen factorization $$S_{\fq}\overset{\delta}{\rightarrow} B\overset{\pi}{\twoheadrightarrow}R\widehat{\otimes}_{C_R}C'$$ of $\psi_{\fq}:S_{\fq}\rightarrow R\widehat{\otimes}_{C_R}C'$,  where $\delta$ is weakly regular, $B$ is a complete local ring and $\pi$ is a surjection whose kernel is generated by a regular sequence, say $\mathbf{x}$. In particular, $K_\bullet(\mathbf{x};B)$ is a flat resolution of $R\widehat{\otimes}_{C_R}C'$ over $S_{\fq}$.
   
    So $$_{S_{\fq}}\mathcal{F}_{\bullet}\otimes_{S_{\fq}}K_\bullet(\mathbf{x};B)\overset{\text{\cite[(A.4.1)]{ChrisensenGorenstein}}}{\simeq} {_{S_{\fq}}\mathcal{F}_{\bullet}}\otimes_{S_{\fq}}(R\widehat{\otimes}_{C_R}C')\overset{(\ref{EquationAurora})}{\cong} \mathcal{F}_{\bullet}$$ implying that $$H_1\big(({_{S_{\fq}}}\mathcal{F}_{\bullet}\otimes_{S_{\fq}}B)\otimes_B K_\bullet(\mathbf{x};B)\big)\cong H_1\big({_{S_{\fq}}}\mathcal{F}_{\bullet}\otimes_{S_{\fq}}K_\bullet(\mathbf{x};B)\big)\cong H_1(\mathcal{F}_{\bullet})=0.$$ So from Fact \ref{FactKoszulComplexVanishingHomology}(ii) we get $H_1({_{S_{\fq}}}\mathcal{F}_{\bullet}\otimes_{S_{\fq}}B)=0$, i.e. $H_1({_{S_{\fq}}}\mathcal{F}_{\bullet})=0$ as   $\delta$ is faithfully flat. It follows that ${_{S_{\fq}}}\mathcal{F}_{\bullet}$ is a truncation of  a free resolution, $P_{\bullet}$, of $N_{S_{\fq}}:=H_0({_{S_{\fq}}}\mathcal{F}_{\bullet})$ (over $S_\fq$). Note that $N_{S_{\fq}}\otimes_{S_\fq}(R\widehat{\otimes}_{C_R}C')\cong N$ by (\ref{EquationAurora}). We have $$\tor^{S_{\fq}}_{1}(N_{S_\fq},R\widehat{\otimes}_{C_R}C')=H_1\big(P_\bullet\otimes_{S_\fq}(R\widehat{\otimes}_{C_R}C')\big)=H_1\big({_{S_{\fq}}}\mathcal{F}_{\bullet}\otimes_{S_{\fq}}(R\widehat{\otimes}_{C_R}C')\big)\overset{\text{(\ref{EquationAurora})}}{\cong}H_1(\mathcal{F}_{\bullet})=0,$$ in other perspective $$H_1\big(N_{S_{\fq}}\otimes_{S_{\fq}}K_\bullet(\mathbf{x};B)\big)=H_1(\mathbf{x};N_{S_{\fq}}\otimes_{S_{\fq}}B)=0$$ ($K_\bullet(\mathbf{x};B)$ is a flat resolution of $R\widehat{\otimes}_{C_R}C'$ over $S_{\fq}$). Consequently, the rigidity  of the Koszul complex (\cite[Corollary 1.6.19]{BrunsHerzogCohenMacaulay}) implies that \begin{equation}
 \label{EquationTorIndependence}
 \tor^{S_{\fq}}_{i}(N_{S_{\fq}},R\widehat{\otimes}_{C_R}C')\cong H_i(\mathbf{x};N_{S_{\fq}}\otimes_{S_{\fq}}B)=0,\ \ \forall\ i\ge 1.
 \end{equation} It follows that $P_\bullet\otimes_{S_{\fq}}(R\widehat{\otimes}_{C_R}C')$ is a free resolution of $N$, and  so considering a minimal free resolution $G_\bullet$ of $T$ over $R$ we get 
 \begin{align*}
   (G_\bullet\otimes_RP_\bullet)\otimes_{S_{\fq}}K_\bullet(\mathbf{x};B)
   &\overset{\text{\cite[(A.4.1)]{ChrisensenGorenstein}}}{\simeq} (G_\bullet\otimes_RP_\bullet)\otimes_{S_{\fq}}(R\widehat{\otimes}_{C_R}C')
   &\\&\overset{\text{}}{\simeq}
   G_\bullet \otimes_RN && (\text{\cite[(A.4.1)]{ChrisensenGorenstein}})
   &\\&\cong 
   \big(G_\bullet\otimes_R (R\widehat{\otimes}_{C_R}C')\big){\otimes_{R\widehat{\otimes}_{C_R}C'}}N
   &\\&\simeq 
 \big(T\otimes_R(R\widehat{\otimes}_{C_R}C')\big)\otimes_{R\widehat{\otimes}_{C_R}C'}^{\mathbf{L}}N 
&\\&
 \overset{\text{(\ref{EquationLastTorVanishing?})}}{\in} \mathcal{D}_b\big(S_{\fq}\big).
 \end{align*}
 From this in conjunction with  Fact \ref{FactKoszulComplexVanishingHomology}(ii),
   by the same way as we used above to prove that $H_1({_{S_{\fq}}}\mathcal{F}_{\bullet})=0$, 
 we can conclude that   $\big(G_\bullet \otimes_R S_{\fq}\big)\otimes_{S_{\fq}}P_\bullet\cong G_\bullet\otimes_RP_\bullet\in \mathcal{D}_b\big(S_{\fq}\big)$, in other words
 
 \begin{equation}
   \label{EquationTorVanishignOverSq}
   \tor^{S_{\fq}}_{\gg}\big(T\otimes_RS_{\fq},N_{S_{\fq}}\big)=0.
 \end{equation}
   
Again by Fact \ref{FactSmoothHomomorphismCollection}(ii),  $S_{\fq}$ is a regular $R$-algebra, thus it is indeed weakly regular over $R$. On the other hand,   as $S$ is  of finite type over $R$ so the induced residue field extension 
$$K_R\rightarrow K_{S_{\fq}}\big(=S_{\fq}/\fq S_{\fq}\overset{\text{as\ } K_R\text{-algebras}}{\cong} \text{Frac}(S/\fq)\big)$$	  
is a finitely generated field  extension, while it is also a separable field extension as it is a subextension of the separable field extension $\overline{\varphi}:K_R\rightarrow K_{R\widehat{\otimes}_{C_R}C'}$. Consequently,  Corollary \ref{CorollaryMixedCharacteristicResiduallyFinitelyGeneratedSeparableWeaklyRegularHomomorphism} implies that $T\otimes_RS_\fq$ is a test module over $S_\fq$. From this fact and  (\ref{EquationTorVanishignOverSq})  we conclude that
 $\pd_{S_{\fq}}(N_{S_\fq})<\infty$.  Therefore from  (\ref{EquationTorIndependence}) we deduce that $\pd_{R\widehat{\otimes}_{C_R}C'}(N)<\infty$, as was to be proved  ($N=N_{S_{\fq}}\otimes_{S_\fq}(R\widehat{\otimes}_{C_R}C')$). 	
 	\end{proof}
\end{sublem}

Now we  can present our alternative proof of Corollary \ref{TestModuleAndFlatHomomorphism} that is not  using \cite[Theorem 4.8]{WagstaffAscent}, for residually separable weakly regular homomorphisms.

\begin{subcor}\label{CorollaryAlternativeProofCWQuestionSeparableCase} Suppose that $\varphi:R\rightarrow S$ is a flat local homomorphism   with regular closed fiber $S/\fm_RS$ such that $\overline{\varphi}:K_R\rightarrow K_S$ is a separable field extension. Let $T$ be a test $R$-module.  If  $K_R$ is uncountable or $\text{CI-dim}_R(T)<\infty$, then $T\otimes_R S$ is a test module for $S$.
\begin{proof}
   
  We can assume that $R$ and $S$ are complete local rings, as we argued in the first paragraph of the proof of  Corollary \ref{TestModuleAndFlatHomomorphism} (without loosing the separability of the induced residue field extension).  Afterwards, in view of  Fact \ref{FactReductionToWeaklyUnramified} we may, and we do, assume that $\varphi$ is weakly unramified. 	  

Considering a coefficient ring $\lambda_R:C_R\rightarrow R$ of $R$, by Proposition \ref{PropositionResiduallySeparableImpliesHomomorphismOfCoefficientRingsExist} and our hypothesis there exists a coefficient ring $\lambda_S:C_S\rightarrow S$,  of $S$,  and $\sigma:C_R\rightarrow C_S$ such that $\lambda_S\circ \sigma=\varphi\circ \lambda_R$. Thus by Proposition \ref{LemmaNiceFactorization}(ii) we may, and we do, assume that $S=R\widehat{\otimes}_{C_R}C_S$ and $\varphi=\eta_{R,R\widehat{\otimes}_{C_R}C_S}$. Note that $\overline{\sigma}:K_{C_R}\rightarrow K_{C_S}$ is also a separable field extension, as $\overline{\sigma}\circ \overline{\lambda_R}^{\ -1}$ is $K_R$-algebra isomorphic to $\overline{\eta_{R,R\widehat{\otimes}_{C_R}C_S}}$ by Fact \ref{FactCompleteTensorProductFactCoefficientRingExtension}(iii) (and $\overline{\lambda_R}$ is an isomorphism).
 So our statement follows from Lemma \ref{LemmaCoefficientPRingBaseChangeIsOKWhenSeparableResidueFieldExtension} (mixed characteristic case), and Corollary \ref{BasedFieldExtensionLemma}(i) (equicharacteristic case).
 \end{proof}
\end{subcor}  

\subsection{An alternative proof for  the Cohen-Macaulay and equicharacteristic case} In this subsection all rings have equicharacteristic. 
In the Cohen-Macaulay case  ($R$ being Cohen-Macaulay),  we relax the separability condition    in the statement of      Corollary \ref{CorollaryAlternativeProofCWQuestionSeparableCase}, by applying the following  Proposition.  The following proposition is analogous to Proposition \ref{LemmaNiceFactorization}, but the coefficient ring containment hypothesis in Proposition \ref{LemmaNiceFactorization} is now replaced with another condition.
Since only residually inseparable  case (thus prime characteristic case, as $R$ and $S$ have equicharacteristic)   of the next proposition shall be used in the sequel, so   the following proposition has been stated only for prime characteristic rings.

\begin{subprop}\label{ArtinianFactorization} Suppose that $\varphi:R\rightarrow S$ is a flat local homomorphism of Artinian local rings of prime characteristic $p>0$. Then there is a  flat local homomorphisms   $\widetilde{\varphi}:\widetilde{R}\rightarrow \widetilde{S}$ of Artinian local rings as well as weakly unramified flat local homomorphisms $\nu_{\widetilde{R}}:\widetilde{R}\overset{\text{}}{\longrightarrow} R,\ \nu_{\widetilde{S}}:\widetilde{S}\overset{\text{}}{\rightarrow} S$  such that they fit into the  commutative diagram
	\begin{equation}
	\label{EquationCommutativeDiagram}
	\begin{CD}R@>\varphi>>S\\@A\nu_{\widetilde{R}}AA@AA\nu_{\widetilde{S}}A\\\widetilde{R}@>>\widetilde{\varphi}>\widetilde{S}\end{CD}
	\end{equation}
	and such that $\widetilde{\varphi}$ admits a factorization $(\widetilde{\varphi}=\widetilde{\varphi}_{\mathcal{L}}\circ \eta)$, $$\widetilde{R}\overset{\eta}{\rightarrow} \widetilde{R}_\mathcal{L}:=\widetilde{R}\otimes_{F}\mathcal{L}\overset{\widetilde{\varphi}_{\mathcal{L}}}{\longrightarrow} \widetilde{S}$$  into a  coefficient field base change $\eta:\widetilde{R}\overset{r\mapsto r\otimes 1}{\longrightarrow} \widetilde{R}_\mathcal{L}$  followed by a  module-finite flat local homomorphism $\widetilde{\varphi}_{\mathcal{L}}:\widetilde{R}_\mathcal{L}\rightarrow \widetilde{S}$ ($F$  is a coefficient field of $\widetilde{R}$ and $F\rightarrow \mathcal{L}$ is a field extension).
	
	Additionally, let $T\in \text{mod}(R)$ and $N\in \text{mod}(S)$. Then the factorization can be chosen  so that  there exist $\widetilde{T}\in \text{mod}(\widetilde{R})$ and  $\widetilde{N}\in\text{mod}(\widetilde{S})$ with  $$\widetilde{T}\otimes_{\widetilde{R}}R=T,\  \widetilde{N}\otimes_{\widetilde{S}}S=N$$ and $\widetilde{T}$ is a test module for $\widetilde{R}$ provided $T$ is a test module for $R$.
	
	\begin{proof}
		By Fact \ref{FactArtinianIsPolynomialQuotient}	we can present $R$ and $S$ as quotients of polynomial rings, say $$(R,\fm_R)=\big(C_R[\mathbf{X}:=X_1,\ldots,X_n]/\maa,(\mathbf{X})\big),\ \ \ \ (S,\fm_S)=\big(C_S[\mathbf{Y}:=Y_1,\ldots,Y_m]/\mab,(\mathbf{Y})\big)$$ where $C_R\subseteq R$ (resp. $C_S\subseteq S$) is a coefficient field  of $R$ (resp. of $S$).  
		
		Consider the prime subfield  $\mathbb{F}_p=\mathbb{Z}/p\mathbb{Z}$ of $C_R$ and of $C_S$, and the  natural inclusion maps from $\mathbb{F}_p$ to $C_R$ and $C_S$ and thus to each of $C_R[\mathbf{X}]$, $C_S[\mathbf{Y}]$, $R$ and $S$. 
		
		Let $R^u\overset{f}{\rightarrow} R^h\twoheadrightarrow T\rightarrow 0$ be a presentation of $T$ over $R$ and    $H$ be the   matrix of the map $f:R^u\rightarrow R^v$ of finite free modules. Afterwards, let
		$\underline{\zeta}\subseteq C_R$ be the set of   coefficients appearing in the entries of the matrix $H$ as well as the  coefficients appearing in a minimal generating set of‌ the ideal $\maa$ of $R$. Namely, suppose that $\fa$ is generated by $t$-elements and  suppose that $j$-th generator of $\fa$ is 
		\begin{equation}
		\label{EquationGeneartorsOfa}
		_jc_{\underline{i}_1}X_1^{i_{1,1}}\ldots X_n^{i_{1,n}} +\cdots+\ _jc_{\underline{i}_w}X_1^{i_{w,1}}\ldots X_n^{i_{w,n}},\ \ \  \underline{i}_1=(i_{1,1},\ldots,i_{1,n}),\ldots,\underline{i}_w=(i_{w,1},\ldots,i_{w,n})\in \mathbb{N}_0^n
		\end{equation}  
		and the $(k,l)$-th entry of $H$ modulo $\fa$ is 
		\begin{equation}
		\label{EquationEntryOfMatrix}
		_{k,l}c'_{\underline{i}_1}X_1^{i_{1,1}}\ldots X_n^{i_{1,n}} +\cdots+\ _{k,l}c'_{\underline{i}_w}X_1^{i_{w,1}}\ldots X_n^{i_{w,n}},\ \ \  \underline{i}_1=(i_{1,1},\ldots,i_{1,n}),\ldots,\underline{i}_w=(i_{w,1},\ldots,i_{w,n})\in \mathbb{N}_0^n
		\end{equation}
		(since there are finitely many elements appearing in a minimal generating set of $\fa$ and in the entries of $H$ we may, and we do, fix $w\in \mathbb{N}$ and $\underline{i}_1,\ldots,\underline{i}_w\in \mathbb{N}^n_0$ in the above two displays working for all generators of $\fa$ and all entries of $H$,  at the cost of allowing some coefficients to be zero). Set  $$\underline{\zeta}:=\{_1c_{\underline{i}_1},\ldots,{_1c_{\underline{i}_w}},\ldots,{_tc_{\underline{i}_w}},\ldots,{_tc_{\underline{i}_w}}\} \cup \{_{1,1}c'_{\underline{i}_1},\ldots,_{1,1}c'_{\underline{i}_w},\ldots,{_{1,u}c'_{\underline{i}_1}},\ldots,{_{1,u}c'_{\underline{i}_w}},\ldots,{_{h,u}c'_{\underline{i}_1}},\ldots,{_{h,u}c'_{\underline{i}_w}}\}.$$

		The subfield $\mathbb{F}_p(\underline{\zeta})$ of $C_R$ generated by $\underline{\zeta}$ contains all data needed to  simulate  $R$ and the $R$-module $N$,  such that the simulation has coefficient field   $\mathbb{F}_p(\underline{\zeta})$. 
		Namely, let $\widetilde{\fa}$  be the ideal of $\mathbb{F}_p(\underline{\zeta})[\mathbf{X}]$ whose $j$-th generator is   (\ref{EquationGeneartorsOfa})   ($1\le j\le t$) and set $$\widetilde{R}:=\mathbb{F}_p(\underline{\zeta})[\mathbf{X}]/\widetilde{\fa}.$$ Also let $\widetilde{H}$ to be the $h\times u$ matrix with entries in $\widetilde{R}$  whose $(k,l)$-th entry is (the image in  $\mathbb{F}_{p}(\underline{\zeta})[\mathbf{X}]/\widetilde{a}$  of) (\ref{EquationEntryOfMatrix}) and set $\widetilde{T}=\text{Coker}(\widetilde{H})$.
		Note that $\mathbb{F}_p(\underline{\zeta})[\mathbf{X}]\overset{\text{inclusion}}{\rightarrow} C_R[\mathbf{X}]$ is a flat extension (Fact \ref{FaithfullyFlat}(ii)). Then,
		\begin{itemize}	
			\item [(i)] The inclusion $\mathbb{F}_{p}(\underline{\zeta})\overset{\text{}}{\rightarrow}C_R$ and the natural map $\mathbb{F}_{p}(\underline{\zeta})\overset{\text{incl.}}{\rightarrow}\mathbb{F}_{p}(\underline{\zeta})[\mathbf{X}]\overset{\text{nat. epi.}}{\twoheadrightarrow}\widetilde{R}$ yields the tensor product ring  $\widetilde{R}\otimes_{\mathbb{F}_p(\underline{\zeta})}C_R$ and the isomorphism $$\widetilde{R}\otimes_{\mathbb{F}_p(\underline{\zeta})}C_R\cong C_R[\mathbf{X}]/\widetilde{\fa}C_R[\mathbf{X}]=C_R[\mathbf{X}]/\fa=R.$$ So we have  the  weakly unramified flat local homomorphism $$\nu_{\widetilde{R}}:\widetilde{R}\overset{\text{flat coefficient field base change}}{\longrightarrow} \widetilde{R}\otimes_{\mathbb{F}_p(\underline{\zeta})}C_R\overset{\cong}{\rightarrow} R,\ \ f(\mathbf{X})+\widetilde{\fa}\mapsto f(\mathbf{X})+\fa.$$
			\item[(ii)] $\widetilde{T}\otimes_{\widetilde{R}}R=\text{Coker}(\widetilde{H})\otimes_{\widetilde{R}}R\cong \text{Coker}\big(\nu_{\widetilde{R}}(\widetilde{H})\big)=\text{Coker}(H)=T$. Moreover, $\widetilde{T}$ is a test $\widetilde{R}$-module provided $T$ is a test $R$-module by Remark \ref{RemarkTestModulePropertyDesentFromFaithfullyFlatExtension}.
		\end{itemize}

		In the same vein, we consider  the set $\underline{\xi}_1$ consisting of the   coefficients appearing in a minimal generating set of the defining ideal $\mab$ of $S$, as well as the set $\underline{\xi}_2$ consisting of coefficients appearing in the matrix of  a presentation of $N$ over $S$. The data stored in $\underline{\xi}_1\cup \underline{\xi}_2$,  is sufficient for constructing a  simulation $\widetilde{S}$ of $S$ whose coefficient field is a finitely generated field extension  of $\mathbb{F}_p$, and is enough for constructing a finitely generated $\widetilde{S}$-module $\widetilde{N}$ such that $\widetilde{S}$ and $\widetilde{N}$ satisfy the  properties similar to (i) and (ii) above (precisely as we did in the above, but by using $\underline{\xi}_1\cup \underline{\xi}_2, S,\widetilde{S}, \widetilde{N}$ and $N$ in place of $\underline{\zeta}$, $R$, $\widetilde{R},\widetilde{T}$ and $T$ respectively). However here, we also add the following (finitely many) extra coefficients to the  previously mentioned coefficients stored in $\mathbb{\xi}_1\cup\underline{\xi}_2$:
		
		Assuming for each $\zeta\in \underline{\zeta}$ and each $X_i$ that 
		\begin{equation}
		\label{EquationCoefficientsRequiredForConstructingHomomorphism}
		\varphi(\zeta)=(\sum\limits_{k=1}^v\xi''_{\zeta,\underline{j}_k}Y_1^{j_{k,1}}\ldots Y_m^{j_{k,m}})+\mab,\ \ \ \ 
		\varphi(X_i)=(\sum\limits_{k=1}^v\xi''_{X_i,\underline{j}_k}Y_1^{j_{k,1}}\ldots Y_m^{j_{k,m}})+\mathfrak{b},
		\end{equation}
		we set $$\underline{\xi}=\{\xi''_{\zeta,\underline{j}_k}:\ \zeta\in \underline{\zeta},\ 1\le k\le v\}\cup \{\xi''_{X_i,\underline{j}_k}:\ 1\le i\le n,\ 1\le k\le v\}\cup \underline{\xi}_1\cup \underline{\xi}_2.$$
		
		Here  perhaps it is worth to stress that $\varphi(\zeta)$ is an element of $S$ but not necessarily an element of $C_S$, because $\varphi$ does not necessarily map the coefficient ring $C_R$ to a coefficient field of $S$. This extra data stored in $\underline{\xi}$ is required to construct a ring homomorphism $\widetilde{\varphi}:\widetilde{R}\rightarrow \widetilde{S}$ such that $\widetilde{\varphi}$ fits into the  commutative diagram (\ref{EquationCommutativeDiagram}).

		Again, $\mathbb{F}_p(\underline{\xi})[\mathbf{Y}]$ has an ideal $\widetilde{\mab}$  such that $\widetilde{S}:=\mathbb{F}_p(\underline{\xi})[\mathbf{Y}]/\widetilde{\fb}$ admits a module $\widetilde{N}$ with 
		
		\begin{itemize}	
			\item [(i$'$)] The inclusion $\mathbb{F}_{p}(\underline{\xi})\overset{\text{}}{\rightarrow}C_S$ and the natural map $\mathbb{F}_{p}(\underline{\xi})\overset{\text{incl.}}{\rightarrow}\mathbb{F}_{p}(\underline{\xi})[\mathbf{Y}]\overset{\text{nat. epi.}}{\twoheadrightarrow}\widetilde{S}$ yields the the tensor product ring  $\widetilde{S}\otimes_{\mathbb{F}_p(\underline{\xi})}C_S$  and the isomorphism $$\widetilde{S}\otimes_{\mathbb{F}_p(\underline{\xi})}C_S\cong C_S[\mathbf{Y}]/\widetilde{\fb}C_S[\mathbf{Y}]=C_S[\mathbf{Y}]/\fb=S$$ and we have  the  weakly unramified flat local homomorphism $$\nu_{\widetilde{S}}:\widetilde{S}\overset{\text{flat coefficient field base change}}{\longrightarrow} \widetilde{S}\otimes_{\mathbb{F}_p(\underline{\xi})}C_S\overset{\cong}{\rightarrow} S,\ \ f(\mathbf{Y})+\widetilde{\fb}\mapsto f(\mathbf{Y})+\fb.$$
			\item[(ii$''$)] $\widetilde{N}\otimes_{\widetilde{S}}S\cong N$.
		\end{itemize} 
		
		Note that $\widetilde{S}$ and $\widetilde{R}$ are also Artinian (thus complete) local rings, as $\nu_{\widetilde{S}}:\widetilde{S}\rightarrow S$ and $\nu_{\widetilde{R}}:\widetilde{R}\rightarrow R$ are   flat weakly unramified local homomorphisms and $S,R$ are Artinian (see \cite[Theorem A.11]{BrunsHerzogCohenMacaulay}). Later, we will show that there is a ring homomorphism $\widetilde{\varphi}:\widetilde{R}\rightarrow \widetilde{S}$ (which has to be flat necessarily, as we will see) such that it fits into the commutative diagram (\ref{EquationCommutativeDiagram}). But,     we first assume the existence of $\widetilde{\varphi}$ and we prove the statement of the proposition.

		Let $\mathcal{L}$ be a maximal subfield  of $\widetilde{S}$  containing the subfield  $\widetilde{\varphi}\big(\mathbb{F}_p(\underline{\zeta})\big)$ of  $\widetilde{S}$. ($\mathcal{L}$ exists in view of the Zorn's lemma). By Fact \ref{FactMaximalSubfieldAlgebraicExtension}  the extension of fields   $\mathcal{L}\rightarrow \widetilde{S}\twoheadrightarrow \widetilde{S}/(\mathbf{Y})$ is algebraic. 
		We consider the tensor product ring $\widetilde{R}_\mathcal{L}:=\widetilde{R}\otimes_{\mathbb{F}_p(\underline{\zeta})}\mathcal{L}$ assigned to $\widetilde{\varphi}:\mathbb{F}_p(\underline{\zeta})\rightarrow \mathcal{L}$ and $\mathbb{F}_{p}(\underline{\zeta})\overset{\text{incl.}}{\rightarrow}\mathbb{F}_{p}(\underline{\zeta})[\mathbf{X}]\overset{\text{nat. epi.}}{\twoheadrightarrow}\widetilde{R}$. Note that since $\widetilde{R}$ is  of finite type over $\mathbb{F}_p(\underline{\zeta})$, so $\widetilde{R}_{\mathcal{L}}$ is Noetherian. Then we have the ring homomorphism $$\widetilde{\varphi}_{\mathcal{L}}:\widetilde{R}_\mathcal{L}\rightarrow \widetilde{S},\ \ \ \widetilde{\varphi}_{\mathcal{L}}(r\otimes \kappa)=\widetilde{\varphi}(r)\kappa\  \ \ (\kappa\in \mathcal{L},\ r\in \widetilde{R}).$$ It is clear, by the definition,    that $\widetilde{\varphi}$ is factored into $\widetilde{R}\overset{r\mapsto r \otimes 1}{\rightarrow}\widetilde{R}_{\mathcal{L}}$ and $\widetilde{\varphi}_{\mathcal{L}}:\widetilde{R}_{\mathcal{L}}\rightarrow S$. To see also that $\widetilde{\varphi}_{\mathcal{L}}$ is module finite, note that  the induced residue field extension $\overline{\widetilde{\varphi}_{\mathcal{L}}}$ fits into the commutative diagram of field extensions
		\begin{center}
			$\begin{CD}
			\mathcal{L} @>\gamma:=\text{inclusion}>> \widetilde{S}\\
			@V\kappa\mapsto (1\otimes \kappa)+(\mathbf{X})V\cong V @VV\mu:=\text{natural epimorphism}V\\
			K_{\widetilde{R}_\mathcal{L}}=\widetilde{R}_{\mathcal{L}}/(\mathbf{X}) @>>\overline{\widetilde{\varphi}_{\mathcal{L}}}>  \widetilde{S}/(\mathbf{Y})=K_{\widetilde{S}}
			\end{CD}$
		\end{center}
		thus $\overline{\widetilde{\varphi}_{\mathcal{L}}}$ is an algebraic  field extension in view of the above commutative diagram, because   $\mu\circ \gamma$ is algebraic as mentioned before (it is easily seen that the left vertical map in the diagram is isomorphism by the same reason as in the proof of Fact \ref{FactCompleteTensorProductFactCoefficientRingExtension}(iv)). On the other hand $\overline{\widetilde{\varphi}_{\mathcal{L}}}$  is also a finitely generated field extension, because  $ K_{\widetilde{S}}=\mathbb{F}_{p}(\underline{\xi})$ is  finitely generated even over  $\mathbb{F}_p$.
		It follows that $\overline{\widetilde{\varphi}_{\mathcal{L}}}$  is both of finitely generated and algebraic, thus a finite extension of fields. 
		So, by  applying Fact \ref{FactArtinianFiniteDimensiona} to $K_{\widetilde{R}_{\mathcal{L}}}\rightarrow \widetilde{S}/\fm_{\widetilde{R}_{\mathcal{L}}}\widetilde{S}$ we can see that $\widetilde{S}/\fm_{\widetilde{R}_{\mathcal{L}}}\widetilde{S}$ is  a finitely generated $K_{\widetilde{R}_{\mathcal{L}}}$-module (and so a finitely generated $\widetilde{R}_{\mathcal{L}}$-module).
 Hence from Fact \ref{FactCompleteNakayamsLemma} we conclude that $\widetilde{S}$ is a finitely generated $\widetilde{R}_{\mathcal{L}}$-module (via $\widetilde{\varphi}_{\mathcal{L}}$). Finally,  applying \cite[Tag 00MP]{Stacks} to $\widetilde{R}\rightarrow \widetilde{R}_{\mathcal{L}}\rightarrow \widetilde{S}$ we  can deduce the flatness of $\widetilde{\varphi}_{\mathcal{L}}$ from the flatness of $\widetilde{\varphi}$. However, $S\rceil_{\nu_{\widetilde{S}}\circ \widetilde{\varphi}}\cong \widetilde{S}\otimes_{\widetilde{S}}S$ is faithfully flat over $\widetilde{R}$, because of the commutative diagram (\ref{EquationCommutativeDiagram})  and faithfully flatness of $\nu_{\widetilde{R}}$ and $\varphi$. So $-\otimes_R(\widetilde{S}\otimes_{\widetilde{S}}S)\cong (-\otimes_R\widetilde{S})\otimes_{\widetilde{S}}S$ is an exact functor on the category of $R$-modules.    This implies that $\widetilde{\varphi}:\widetilde{R}\rightarrow\widetilde{S}$ is flat,  as $\nu_{\widetilde{S}}:\widetilde{S}\rightarrow S$ is faithfully flat.

		It remains to deal with the existence of $\widetilde{\varphi}$, as promised before. Let us consider two finite sets of indeterminates $\Theta=\{\Theta_{\zeta}:\zeta\in \underline{\zeta}\}$  and $\Xi=\{\Xi_{\xi}:\xi\in \underline{\xi}\}$   and the presentations 
		\begin{align*}
		\mu:\Frac(\mathbb{F}_p[\Theta]/\mathfrak{k}_1)\overset{\cong}{\rightarrow}\mathbb{F}_p(\underline{\zeta}),\ \ \  \mathfrak{k}_1=\Ker\big(\mathbb{F}_p[\Theta]\overset{\Theta_\zeta\mapsto \zeta}{\longrightarrow} \mathbb{F}_p(\underline{\zeta})\big)\\
		\mu':\Frac(\mathbb{F}_p[\Xi]/\mathfrak{k}_2)\overset{\cong}{\rightarrow}\mathbb{F}_p(\underline{\xi}),\ \ \  \mathfrak{k}_2=\Ker\big(\mathbb{F}_p[\Xi]\overset{\Xi_{\xi}\mapsto \xi}{\rightarrow} \mathbb{F}_p(\underline{\xi})\big).
		\end{align*}
		Let $$\mu_{\mathbf{X}}:\Frac(\mathbb{F}_p[\Theta]/\mathfrak{k}_1)[\mathbf{X}]\rightarrow \mathbb{F}_p(\underline{\zeta})[\mathbf{X}]\ \ \text{(resp. \ \ }\mu'_{\mathbf{Y}}:\Frac(\mathbb{F}_p[\Xi]/\mathfrak{k}_2)[\mathbf{Y}]\rightarrow \mathbb{F}_p(\underline{\xi})[\mathbf{Y}]\text{)}$$ be the natural isomorphism extending $\mu$ (resp. $\mu'$) by the rule $X_i\mapsto X_i$ (resp. $Y_i\mapsto Y_i$).  There are ideals $\widetilde{\fa}_1$ and $\widetilde{\fb}_1$ of $\Frac(\mathbb{F}_p[\Theta]/\mathfrak{k}_1)[\mathbf{X}]$   and  $\Frac(\mathbb{F}_p[\Xi]/\mathfrak{k}_2)[\mathbf{Y}]$  respectively such that $\mu_{\mathbf{X}}(\widetilde{\fa}_1)=\widetilde{\fa}$ and $\mu'_{\mathbf{Y}}(\widetilde{\fb}_1)=\widetilde{\fb}$, thus we have isomorphisms $$\overline{\mu_{\mathbf{X}}}:\Frac(\mathbb{F}_p[\Theta]/\mathfrak{k}_1)[\mathbf{X}]/\widetilde{\fa}_1\overset{\cong}{\rightarrow} \widetilde{R},\ \ \ \ \ \  \overline{\mu'_{\mathbf{Y}}}:\Frac(\mathbb{F}_p[\Xi]/\mathfrak{k}_2)[\mathbf{Y}]/\widetilde{\fb}_1\overset{\cong}{\rightarrow} \widetilde{S}.$$
		
		First we have the following commutative diagram:
		\begin{center}
			$\begin{CD}
			\mathbb{F}_p[\Theta] @>\Theta_{\zeta}\mapsto \zeta  > > \widetilde{R}=\mathbb{F}_p(\underline{\zeta})[\mathbf{X}]/\widetilde{\fa}   @>\nu_{\widetilde{R}} >>   R=C_R[\mathbf{X}]/\maa \\
			@V\alpha VV                          @.                                                             @V\varphi VV\\
			\Frac(\mathbb{F}_p[\Xi]/\mathfrak{k}_2)[\mathbf{Y}]/\widetilde{\mab}_1  @>\overline{\mu'_{\mathbf{Y}}} >\cong>  \widetilde{S}=\mathbb{F}_p(\underline{\xi})[\mathbf{Y}]/\widetilde{\mab}  @>\nu_{\widetilde{S}}>1-1> S=C_S[\mathbf{Y}]/\mab,
			\end{CD}$ 
		\end{center}
		where $\alpha$‌ is  is given by the rule $$\alpha(\Theta_\zeta)= \sum\limits_{k=1}^v(\Xi_{\xi''_{\zeta,\underline{j}_k}}+\mathfrak{k}_2)Y_1^{j_{k,1}}\ldots Y_m^{j_{k,m}}+\widetilde{\mathfrak{b}}_1$$ so that the diagram is commutative  in view of (\ref{EquationCoefficientsRequiredForConstructingHomomorphism}).

		The injectivity of $\nu_{\widetilde{S}}$ and $\overline{\mu'_{\mathbf{Y}}}$ and the commutativity of the diagram shows that $\alpha(\mathfrak{k}_1)=0$, hence we may replace $\mathbb{F}_p[\Theta]$ with $\mathbb{F}_p[\Theta]/\mathfrak{k}_1$ in the diagram while keeping the diagram commutative. Similarly, since $\overline{\mu'_{\mathbf{Y}}}$ and $\nu_{\widetilde{S}}$ are local homomorphisms and the diagram is commutative so any $0\neq f\in \mathbb{F}_p[\Theta]/\mathfrak{k}_1$ is mapped to an invertible element by $\alpha$, because  $\varphi\circ \nu_{\widetilde{R}}\circ \mu (f)=\nu_{\widetilde{S}}\circ\overline{\mu'_{\mathbf{Y}}}\circ \alpha(f)$ is an invertible element of $S$ ($\mu(f)$ is invertible). Therefore, $\alpha$ extends to $\Frac(\mathbb{F}_p[\Theta]/\mathfrak{k}_1)\rightarrow\Frac(\mathbb{F}_p[\Xi]/\mathfrak{k}_2)[\mathbf{Y}]/\widetilde{\mab}_1$, and we can replace $\mathbb{F}_p[\Theta]/\mathfrak{k}_1$ with $\Frac(\mathbb{F}_p[\Theta]/\mathfrak{k}_1)$ in the diagram.
		Afterwards, we  define the ring homomorphism  $$\beta:\Frac(\mathbb{F}_p[\Theta]/\mathfrak{k}_1)[\mathbf{X}]\rightarrow \Frac(\mathbb{F}_p[\Xi]/\mathfrak{k}_2)[\mathbf{Y}]/\widetilde{\mab}_1,\ \ \ X_i\mapsto \sum\limits_{k=1}^v(\Xi_{\xi''_{X_i,\underline{j}_k}}+\mathfrak{k}_2)Y_1^{j_{k,1}}\ldots Y_m^{j_k,m}+\widetilde{\mathfrak{b}}_1$$ 
		extending $\alpha:\Frac(\mathbb{F}_p[\Theta]/\mathfrak{k}_1)\rightarrow \Frac(\mathbb{F}_p[\Xi]/\mathfrak{k}_2)[\mathbf{Y}]/\widetilde{\mab}_1$.  Then, in view of (\ref{EquationCoefficientsRequiredForConstructingHomomorphism}),  if we replace $\alpha$ and its source in the diagram with $\beta$ and its source,   the  diagram remains commutative (with the homomorphism $\Frac(\mathbb{F}_p[\Theta]/\mathfrak{k}_1)[\mathbf{X}]\overset{\mu_{\mathbf{X}}}{\rightarrow}\mathbb{F}_p(\underline{\zeta})[\mathbf{X}]\overset{\text{nat. epi.}}{\twoheadrightarrow}\widetilde{R}$,  in the top row).  Finally, again the commutativity of the diagram as well as the injectivity of $\nu_{\widetilde{S}}\circ \overline{\mu'_{\mathbf{Y}}}$ implies that $\widetilde{\fa}_1\subseteq \ker(\beta)$. Now, the induced map $\overline{\beta}:\Frac(\mathbb{F}_p[\Theta]/\mathfrak{k}_1)[\mathbf{X}]/\widetilde{\fa}_1\rightarrow \Frac(\mathbb{F}_p[\Xi]/\mathfrak{k}_2)[\mathbf{Y}]/\widetilde{\mab}_1$ in conjunction with the isomorphisms $\overline{\mu_{\mathbf{X}}}$ and $\overline{\mu'_{\mathbf{Y}}}$ shows us that $$\widetilde{\varphi}:=\overline{\mu'_{\mathbf{Y}}}\circ \overline{\beta} \circ \overline{\mu_{\mathbf{X}}}^{-1}:\widetilde{R}\rightarrow \widetilde{S}$$ is our desired ring homomorphism.
	\end{proof}
\end{subprop}

Now, we present our alternative proof for Corollary \ref{TestModuleAndFlatHomomorphism} when $R$ is a  Cohen-Macaulay ring containing a field.

\begin{subcor} \label{CorollaryAlternativeProofCohenMacaulayCase}
		Suppose that $\varphi:R\rightarrow S$ is a weakly regular homomorphism of local rings  and that $T$ is a test module over $R$. Assume that $R$ is Cohen-Macaulay and $R$ contains a field. If $K_R$ is uncountable, or $\text{CI-dim}_R(T)<\infty$ then $T\otimes_RS$ is a test module over $S$.
\begin{proof} If $R$ contains $\mathbb{Q}$, then $K_R$ has characteristic zero so $\overline{\varphi}:K_R\rightarrow K_S$ is a separable field extension and the statement follows from Corollary \ref{CorollaryAlternativeProofCWQuestionSeparableCase}. So we assume that $R$ has prime characteristic. 
 By Fact \ref{FactReductionToWeaklyUnramified} we may, and we do, assume that $\varphi$ is a weakly unramified homomorphism. We consider some $N\in \text{mod}(S)$ such that 
 \begin{equation}
    \label{EquationIsThisLastOne}
    \tor^S_{\gg}(T\otimes_RS,N)=0
 \end{equation} 
 and we prove that $\pd_S(N)<\infty$.
 
    Considering,   sufficiently high syzygies of $N$ and $T$,   we may assume that $T$‌ and $N$ are both maximal Cohen-Macaulay modules (since $R$ is Cohen-Macaulay and $\varphi$ is weakly regular, so $S$ is also  Cohen-Macaulay by \cite[Theorem 2.1.7]{BrunsHerzogCohenMacaulay}). If $R$ is not Artinian, then in view of  $\fm_R\nsubseteq \fm_R^2\bigcup (\bigcup\limits_{\fp\in \Ass(R)}\fp)$ in conjunction with the Prime Avoidance Lemma  there is some regular element $x_1\in \fm_R\backslash \fm_R^2$. Proceeding in this way, we  can pick a regular sequence $\mathbf{x}$ of $R$ of length $\depth(R)(=\dim(R))$ such that $x_{i+1}\notin (x_1,\ldots,x_i)+\mam_R^2$.  Since $\varphi$ is flat, so $\mathbf{x}$ is a regular sequence on $S$, and also on $T\otimes_RS$ and $N$ ($N$ and $T\otimes_RS$ are high syzygies and so maximal Cohen-Macaulay modules).  Thus, $$\tor^{S/\mathbf{x}S}_{i}\big((T\otimes_RS)/\mathbf{x}(T\otimes_RS),N/\mathbf{x}N\big)\overset{\text{Fact \ref{FactTorOverRTorOverRxR}}}{\cong}
    \tor^{S}_{i}\big((T\otimes_RS)/\mathbf{x}(T\otimes_RS),N\big)$$ for each $i\in \mathbb{N}_0$. Consequently, from Fact \ref{FactTorOfMAndTorOfSpecializationM} (by setting $a=1$ in its statement) and $\tor^S_{\gg}(T\otimes_RS,N)=0$, we conclude that  
    \begin{align}
    \label{EquationTorVanihsingOverArtinianRings}
    \tor^{S/\mathbf{x}S}_{\gg}\big(\ \underset{\cong(T/\mathbf{x}T)\otimes_{R/\mathbf{x}R}(S/\mathbf{x}S)}{\underbrace{{(T\otimes_RS)/\mathbf{x}(T\otimes_RS)}}},\ N/\mathbf{x}N\big)=0.
    \end{align}
    
     But $T/\mathbf{x}T$ is a test module over the Artinian ring $R/\mathbf{x}R$ by \cite[Proposition 2.2(ii)]{Celikbas}, while the induced homomorphism $R/\mathbf{x}R\rightarrow S/\mathbf{x}S$ by $\varphi$ is weakly regular. Note that since $\varphi$ is weakly unramified so $\dim(S)=\dim(R)$ and $S/\mathbf{x}S$ is also Artinian.	 Also $\text{CI-dim}_{R/\mathbf{x}R}(T/\mathbf{x}T)<\infty$ provided $\text{CI-dim}_R(T)<\infty$ by \cite[(1.12) Proposition]{Avramov} (as $\mathbf{x}$ is also a regular sequence over $T$). Hence,  if the statement holds for  weakly regular homomorphism on Artinian rings, then from (\ref{EquationTorVanihsingOverArtinianRings})  we deduce that $\pd_S(N)\overset{\text{Fact \ref{ResolutionModuleRegularSequence}}}{=}\pd_{S/\mathbf{x}S}(N/\mathbf{x}N)<\infty$, as was to be proved.
    
Consequently, without loss of generality, we can assume that $R$ and $S$ are Artinian rings.

By Proposition \ref{ArtinianFactorization}  there is a  flat local homomorphism of Artinian rings $\widetilde{\varphi}:\widetilde{R}\rightarrow \widetilde{S}$ and weakly unramified flat local homomorphisms $\nu_{\widetilde{S}}:\widetilde{S}\rightarrow S$, $\nu_{\widetilde{R}}:\widetilde{R}\rightarrow R$ fitting in the commutative diagram (\ref{EquationCommutativeDiagram}), such that $\widetilde{\varphi}$ admits a factorization $$\widetilde{R}\overset{\eta:\widetilde{r}\mapsto \widetilde{r}\otimes 1}{\longrightarrow} \widetilde{R}_{\mathcal{L}}:=\widetilde{R}\otimes_{F}\mathcal{L}\overset{\widetilde{\varphi}_{\mathcal{L}}}{\longrightarrow }\widetilde{S}$$ where $\widetilde{\varphi}_{\mathcal{L}}$ is a module-finite flat  extension ($\widetilde{\varphi}=\widetilde{\varphi}_{\mathcal{L}}\circ \eta$, and $F$ is a coefficient field of $\widetilde{R}$). Moreover, there are  $\widetilde{N}\in \text{mod}(\widetilde{S})$ and  $\widetilde{R}$-test module $\widetilde{T}$ such that $\widetilde{N}\otimes_{\widetilde{S}}S\cong N$ and $\widetilde{T}\otimes_{\widetilde{R}}R\cong T$.

Let $F_\bullet$ be a free resolution of $\widetilde{T}$. Then (see Notation \ref{NotationRestrictionOfScalars}),  
\begin{align*}
	\big((F_\bullet\otimes_{\widetilde{R}}\widetilde{S})\otimes_{\widetilde{S}}\widetilde{N}\big)\otimes_{\widetilde{S}}S
	\cong 
	(F_\bullet\otimes_{\widetilde{R}}\widetilde{N}\rceil_{\widetilde{\varphi}})\otimes_{\widetilde{S}}S
	&\cong 
	F_\bullet\otimes_{\widetilde{R}}N\rceil_{\underset{=\varphi\circ \nu_{\widetilde{R}}}{\underbrace{\nu_{\widetilde{S}}\circ \widetilde{\varphi}}}}
	&\\&\cong F_\bullet\otimes_{\widetilde{R}}(R\otimes_RN\rceil_{\varphi})
	&\\&\cong (F_\bullet\otimes_{\widetilde{R}}R)\otimes_RN\rceil_{\varphi}
	&\\&\cong (F_\bullet\otimes_{\widetilde{R}}R)\otimes_R(S\otimes_SN)
	&\\&\cong \big((F_\bullet\otimes_{\widetilde{R}}R)\otimes_RS\big)\otimes_SN 
	&\\&\simeq (T\otimes_RS)\otimes_S^{\mathbf{L}}N && (\widetilde{T}\otimes_{\widetilde{R}}R\cong T)
&\\&\in \mathcal{D}_b(S)  && (\ref{EquationIsThisLastOne}).
\end{align*}
	
	It follows that $(F_\bullet\otimes_{\widetilde{R}}\widetilde{S})\otimes_{\widetilde{S}}\widetilde{N}\in \mathcal{D}_b(\widetilde{S})$, i.e.  $\tor^{\widetilde{S}}_{\gg}(\widetilde{T}\otimes_{\widetilde{R}}\widetilde{S},\widetilde{N})=0$,  because of the faithfully flatness of $S$ over $\widetilde{S}$. Therefore considering $\widetilde{N}$ as $\widetilde{R}_{\mathcal{L}}$-module and $\widetilde{R}$-module via restricting the scalars,
		 \begin{align}
	  \label{EquationVanishingOfTorHopefullyLastOne}
	  \tor^{\widetilde{R}_{\mathcal{L}}}_{\gg}(\widetilde{T}\otimes_{\widetilde{R}}\widetilde{R}_{\mathcal{L}},\widetilde{N}\rceil_{\widetilde{\varphi}_{\mathcal{L}}})\overset{\text{Fact \ref{FactTorOverBisTorOverAWhenBflatOverA}}}{\cong} \tor^{\widetilde{R}}_{\gg}(\widetilde{T},\widetilde{N}\rceil_{\underset{=\widetilde{\varphi}}{\underbrace{\widetilde{\varphi}_{\mathcal{L}}\circ \eta}}})\overset{\text{Fact \ref{FactTorOverBisTorOverAWhenBflatOverA}}}{\cong} \tor^{\widetilde{S}}_{\gg}(T\otimes_{\widetilde{R}}\widetilde{S},\widetilde{N})=0.
	\end{align}

Since $\widetilde{R}$ is Artinian and $\mathcal{L}$ is a field, so $\widetilde{R}_{\mathcal{L}}=\widetilde{R}\otimes_{F}\mathcal{L}\cong \widetilde{R}\widehat{\otimes}_{F}\mathcal{L}$ (as $R$-algebras). Thus   $\widetilde{T}\otimes_{\widetilde{R}}\widetilde{R}_\mathcal{L}$  is a test module over $\widetilde{R}_\mathcal{L}$  by  Corollary \ref{BasedFieldExtensionLemma}(i), because: 

In case $\text{CI-dim}_R(T)=\infty$,    $\nu_{\widetilde{R}}:\widetilde{R}\rightarrow R$ fulfills the required  condition on the existence of some  flat local map with target $R'$ such that $K_{R'}$ is uncountable and $\widetilde{T}\otimes_{\widetilde{R}}R'$ is a test module for $R'$, as in the statement of Corollary \ref{BasedFieldExtensionLemma} and   condition (i) of Corollary \ref{CorollaryTranscendentalExtension}. If otherwise  $\text{CI-dim}_R(T)<\infty$, let $\pd_{Q}(T\otimes_{R}R')<\infty$ for some quasi-deformation $R\overset{g}{\rightarrow} R'\twoheadleftarrow Q$. Then $\widetilde{R}\overset{g\circ \nu_{\widetilde{R}}}{\hookrightarrow}R'\twoheadleftarrow Q$ is a quasi-deformation and $\pd_Q(\widetilde{T}\otimes_{\widetilde{R}}R')=\pd_Q\big(\widetilde{T}\otimes_{\widetilde{R}}(R\otimes_RR')\big)=\pd_Q(T\otimes_{R}R')<\infty$. So $\text{CI-dim}_{\widetilde{R}}(\widetilde{T})<\infty$.

Hence from (\ref{EquationVanishingOfTorHopefullyLastOne}) we get $\pd_{\widetilde{R}_{\mathcal{L}}}(\widetilde{N}\rceil_{\widetilde{\varphi}_\mathcal{L}})<\infty$, because $\widetilde{N}$‌ is a finite $\widetilde{R}_\mathcal{L}$-module (as $\widetilde{S}$ is finite over $\widetilde{R}_{\mathcal{L}}$).  Therefore 
\begin{align}
\label{EquationNoThisWasLastTorVanishing}
\tor^{\widetilde{S}}_{\gg}(\widetilde{S}/\mam_{\widetilde{R}}\widetilde{S},\ \widetilde{N})
\overset{\text{Fact \ref{FactTorOverBisTorOverAWhenBflatOverA}}}{=}
\tor^{\widetilde{R}}_{\gg}(\widetilde{R}/\mam_{\widetilde{R}},\ \widetilde{N}\rceil_{\widetilde{\varphi}})
\overset{\text{Fact \ref{FactTorOverBisTorOverAWhenBflatOverA}}}{=}
\tor^{\widetilde{R}_{\mathcal{L}}}_{\gg}( \widetilde{R}_\mathcal{L}/\mam_{\widetilde{R}}\widetilde{R}_\mathcal{L},\ \widetilde{N}\rceil_{\widetilde{\varphi}_{\mathcal{L}}})
\overset{\pd_{\widetilde{R}_{\mathcal{L}}}(\widetilde{N}\rceil_{\widetilde{\varphi}_{\mathcal{L}}})<\infty}{=}0.
\end{align}

Since $\varphi, \nu_{\widetilde{R}}$  and $\nu_{\widetilde{S}}$ are all weakly unramified local homomorphisms and $\nu_{\widetilde{S}}\circ \widetilde{\varphi}=\varphi\circ \nu_{\widetilde{R}}$ so $$(\fm_{\widetilde{R}}\widetilde{S})S=(\fm_{\widetilde{R}}R)S=\fm_RS=\fm_S=\fm_{\widetilde{S}}S.$$ Thus $\fm_{\widetilde{R}}\widetilde{S}=\fm_{\widetilde{S}}$ ($\nu_{\widetilde{S}}$ is faithfully flat). This equality and (\ref{EquationNoThisWasLastTorVanishing}) yield $\text{Tor}^{\widetilde{S}}_{\gg}(K_{\widetilde{S}},\widetilde{N})=0$  i.e. $\pd_{\widetilde{S}}(\widetilde{N})<\infty$. Consequently, $N\cong \widetilde{N}\otimes_{\widetilde{S}}S$ has finite projective dimension, as was to be proved.
\end{proof}
\end{subcor}

The proof of the above corollary provides  an alternative proof for Corollary \ref{MainResult}  that is avoiding the use of \cite{WagstaffAscent} (for equicharacteristic rings).

\begin{subcor}\label{CorollaryAlternativeProofCDTQuestion} Suppose that $R$ is a local ring containing a field. Assume that $R$ admits a test module $T$ with $\text{CI-dim}_R(T)<\infty$. Then $R$ is a complete intersection.
	\begin{proof}
		This corollary follows from Corollary \ref{CorollaryAlternativeProofCohenMacaulayCase}, by exactly the same arguments as in the proof of Corollary \ref{MainResult}.
	\end{proof}
\end{subcor}

  \section*{Acknowledgement}

We are greatly thankful to Kamran Divaani-Aazar for suggesting us to investigate Question \ref{CDTQuestion} and for   his   comments on   earlier versions of this paper.  We are also grateful to Rasoul Ahangari Maleki for a useful discussion which helped us to find the reference \cite{AvramovSunCohomology}. Finally, we thank Simon H\"aberli, Massoud Tousi,  Olgur Celikbas and  Mohsen Gheibi for their  comments. 


\end{document}